\documentclass[leqno,10pt]{amsart}
\usepackage{amsmath,amssymb,amscd,latexsym,amsthm}
\usepackage{amsfonts}
\usepackage[mathscr]{eucal}
\usepackage{enumerate,verbatim,calc}
\usepackage{epsfig}
\usepackage[all]{xy}
\xyoption{curve}
\SelectTips{eu}{}
\CompileMatrices
\usepackage{graphics}
\usepackage{tikz}
\usepackage[colorlinks,pagebackref]{hyperref}

\newtheorem{slem}[equation]{Lemma}
\newtheorem{sprop}[equation]{Proposition}
\newtheorem{scor}[equation]{Corollary}

\theoremstyle{definition}

 \newtheorem{sdef}[equation]{Definition}
 \newtheorem{srem}[equation]{Remark}
 \newtheorem{srems}[equation]{Remarks}
 \newtheorem{sscholium}[equation]{Scholium}
 \newtheorem{exam}[equation]{Example}
 \newtheorem{exams}[equation]{Examples}

 \newtheorem{subcosa}[equation]{}

\newcommand{\Cref}[1]{Corollary~\textup{\ref{#1}}}
\newcommand{\Dref}[1]{Definition~\textup{\ref{#1}}}
\newcommand{\Eref}[1]{Example~\textup{\ref{#1}}}
\newcommand{\Lref}[1]{Lemma~\textup{\ref{#1}}}
\newcommand{\Pref}[1]{Proposition~\textup{\ref{#1}}}
\newcommand{\Rref}[1]{Remark~\textup{\ref{#1}}}
\newcommand{\Sref}[1]{Section~\textup{\ref{#1}}}

\def\bilap#1{\hbox to 0pt{\hss#1\hss}}
 \def\Rarrow#1{\bilap{\hbox to#1{\rightarrowfill}}}
 \def\Larrow#1{\bilap{\hbox to#1{\leftarrowfill}}}
\def\Equals#1{\bilap
                  {\hbox{\rule[3.5pt]{#1} {.5pt}}
                   \kern-#1
                   \hbox{\rule[1pt]{#1}{.5pt}}%
                 }}

\def\UnderElement#1#2#3#4{\vbox to 0pt{
\hbox{$
\llap{$\scriptstyle#1$}
\left#2\vbox to #3{}\right.
\rlap{$\scriptstyle#4$}
     $}
\vss}}
\def\lift#1,#2,{\vbox to 0pt{\vskip-#1 ex\hbox{$\scriptstyle #2$}\vss}}

\newcommand{\under}[2]
{\vbox to 0pt{\vskip-#1 ex\hbox{$\scriptstyle #2$}\vss}}

\newcommand{\EQAL}[1]%
{\,\begin{picture}(#1,0)%
\put(0,3){\line(1,0){#1}}%
\put(0,1){\line(1,0){#1}}%
\end{picture}\,}%

\newcommand{\vlto}[1]%
{\,\begin{picture}(#1,3)%
\put(0,2){\vector(1,0){#1}}%
\end{picture}\,}%

\newcommand{\vllarrow}[1]%
{\,\begin{picture}(#1,3)%
\put(#1,2){\vector(-1,0){#1}}%
\end{picture}\,}%

\newcommand{\dirlm}[1]%
  {
     {\lim\hskip-1.58em\lower.65ex
       \hbox{$
                {}_{\stackrel{\lower1ex\hbox
                                        {$\scriptstyle -\!\!\!\longrightarrow$}
                                      }{\vbox to0pt{\vss\vskip1ex
                                            \hbox{$\scriptstyle{}^{#1}$}\vss}}
                   }
            $}
     }
\:}

\newcommand{\subdirlm}[1]%
  {
     {\lim\hskip-1.5em\lower.6ex
       \hbox{$
                   {}_{\stackrel{\lower1ex\hbox
                                           {$\scriptstyle\longrightarrow$}
                                }{ ^{#1} }
                      }
             $}
     }
\:}

\newcommand{\inlm}[1]%
   {
      {\lim\hskip-1.58em\lower.6ex
        \hbox{$
                 {}_{\stackrel{\lower1ex\hbox
                                        {$\scriptstyle \longleftarrow\!\!\!-$}
                              }{\vbox to0pt{\vss\vskip.6ex
                                            \hbox{$\scriptstyle{}^{#1}$}\vss}}
                    }
             $}
      }
\:}

\def\halfsize#1{{\scalebox{.6}{#1}}}
\def\hz#1{{\hbox to 0pt{#1}}}
\def\Iso{\vbox to 0pt{\vss\hbox{$\widetilde{\phantom{nn}}$}\vskip-7pt}}

\def\>{\mspace {1mu}}
\def\<{\mspace{-1mu}}
\renewcommand{\(}{{\textup(}}
\renewcommand{\)}{{\textup)}}
\newcommand{\vs}[1]{\vspace{#1pt}}

\newcommand{\T}{{\mathscr T}}

\newcommand{\I}{{\mathscr I}}
\newcommand{\J}{{\mathscr J}}

\newcommand{\Spec}{{\mathrm {Spec}}}
\newcommand{\Spf}{{\mathrm {Spf}}}
\newcommand{\Supp}{{\mathrm {Supp}}}
\newcommand{\A}{{\mathcal A}}

\newcommand{\Hr}{{\mathrm H}}

\newcommand{\cH}{{\mathcal H}}

\newcommand{\cK}{{\mathcal K}}
\newcommand{\cO}{{\mathcal O}}
\newcommand{\OX}{{\mathcal O_{\<\<X}}}

\newcommand{\OU}{{\mathcal O_U}}
\newcommand{\OV}{{\mathcal O_V}}
\newcommand{\OW}{{\mathcal O_W}}

\newcommand{\D}{{\mathbf D}}
\newcommand{\K}{{\mathbf K}}
\newcommand{\vc}{{\vec{\mathrm{c}}}}

\newcommand{\Dqc}{\D_{\mkern-1.5mu\mathrm {qc}}}
\newcommand{\wDqc}{ \widetilde
         {\vbox to6.5pt{\vss\hbox{$\mathbf D$}}}
   _{\mkern-1.5mu\mathrm {qc}} }
\newcommand{\wDqcp}{\wDqc^{\lower.5ex\hbox{$\scriptstyle+$}}}
\newcommand{\Dvc}{\D_{\<\vc}}
\newcommand{\Dqct}{\D_{\mkern-1.5mu\mathrm{qct}}}

\newcommand\Dpl{\D^{\lift.95,\text{\cmt\char'053},}}
\font\cmt=cmtex10

\newcommand{\qc}{{\mathrm{qc}}}

\newcommand{\bj}{{\boldsymbol j}}

\newcommand{\bt}{\textup{\bf t}}

\newcommand{\R}{{\mathbf R}}

\newcommand{\bL}{{\mathbf L}}
\newcommand{\Hom}{{\mathrm {Hom}}}

\newcommand{\Homb}{{\mathrm {Hom}}^{\bullet}}
\newcommand{\KK}{{\mathcal K_{\<\<\scriptscriptstyle\infty}^\bullet}\<\<}
\newcommand{\Ab}{\mathfrak{Ab}}

\newcommand{\Aqc}{\A_{\>\qc}}
\newcommand{\Avc}{\A_{\vec {\mathrm c}}}
\newcommand{\Aqct}{\A_{\mathrm {qct}}\<}

\newcommand{\ush}[1]{{#1^{\textup{\texttt\#}}}}

\newcommand{\sHom}{\cH{om}}

\newcommand{\iGp}[1]{{\Gamma_{\<\!#1}'}}
\newcommand{\iG}[1]{{\varGamma_{\<\!#1}^{\phantom\prime}}}

\newcommand{\set}{\!:=}
\newcommand{\lra}{\longrightarrow}
\newcommand{\iso}%
{{\mkern8mu\longrightarrow \mkern-25.5mu{}^\sim\mkern17mu}}
\newcommand{\osi}%
{{\mkern8mu\longleftarrow \mkern-24.5mu{}^\sim\mkern16mu}}
\newcommand{\Otimes}{\underset
  {\vbox to 0pt {\vskip-1ex\hbox{$\scriptscriptstyle=$}\vss}}
    {\otimes}\vadjust{\kern.4pt}} 
\newcommand{\BG}{\boldsymbol{\varGamma}}

\newcommand{\BL}{{\boldsymbol\Lambda}}

\newcommand{\smcirc}%
  {{\raise.15ex\hbox to.7em{$\hss \scriptstyle\circ\hss$}}} 
\newcommand{\xto}{\xrightarrow}

\newcommand{\upl}{{\mkern-1.3mu\raise.2ex\hbox{$\sst +$}}} 

\newcommand{\fst}{{f^{}_{\<\<*}}}

\newcommand\steq{\stepcounter{equation}}

\newcommand{\OOtimes}{\underset {\vbox to 0pt
  {\vskip-.78ex\hbox{$\scriptscriptstyle=$}\vss}}
  {\otimes}\vadjust{\kern.4pt}}

\newcommand{\kf}{\kern.5pt}

\renewcommand{\iG}{{\Gamma'}}
\newcommand{\vG}[1]{\varGamma_{\<\!#1}^{}}

\newcommand{\GGt}[1]{\iGp{#1}\lower.01ex\hbox{$^{{\textup t}}$}}

\newcommand{\sX}{_{\<\<X}}

\newcommand{\ma}{{\mathsf a}}
\newcommand{\ml}{{\mathsf l}}

\newcommand{\mr}{{\mathsf r}}
\newcommand{\ms}{{\mathsf s}}

\newcommand{\mU}{{\mathfrak U}}
\newcommand{\mV}{{\mathfrak V}}

\newcommand{\mfs}{{\mathfrak s}}
\def\b1{{\mathbf1}}

\newcommand{\boh}{{\boldsymbol h}}

\newcommand{\lle}{{\scriptscriptstyle \le}}

\newcommand\OT[2]{\Otimes_{\<\<\!#1}\,#2}
\newcommand{\ot}{\otimes}

\newcommand{\smdm}%
  {{\<\raise.15ex\hbox to.7em{\hss{$\scriptscriptstyle\diamond$}\hss}\<}} 

\newcommand{\supp}{\textup{supp}}
\newcommand\sst[1]{\scriptscriptstyle #1}

\newcommand{\circled}[1]{\textcircled{\scriptsize{#1}}}

\newcommand{\va}[1]{\vspace{#1pt}}
\DeclareMathOperator{\via}{{\textup{via}}}

\newcommand{\bpic}{\begin{tikzpicture}}
\newcommand{\epic}{\end{tikzpicture}}

\author[ Joseph Lipman]{Joseph Lipman}
\address{142 Ranch Ln. \\
              Santa Barbara CA 93111, USA}
\email {lipman@math.purdue.edu}
\urladdr{www.math.purdue.edu/\~{}lipman/}

\numberwithin{equation}{thm}
\usepackage{relsize}
\usepackage{xcolor}
\definecolor{antique}{cmyk}{0, 0.06, 0.12, 0}

\makeatletter
\def\@tocline#1#2#3#4#5#6#7{\relax
  \ifnum #1>\c@tocdepth 
  \else
    \par \addpenalty\@secpenalty\addvspace{#2}%
    \begingroup \hyphenpenalty\@M
    \@ifempty{#4}{%
      \@tempdima\csname r@tocindent\number#1\endcsname\relax
    }{%
      \@tempdima#4\relax
    }%
    \parindent\z@ \leftskip#3\relax \advance\leftskip\@tempdima\relax
    \rightskip\@pnumwidth plus4em \parfillskip-\@pnumwidth
    #5\leavevmode\hskip-\@tempdima
      \ifcase #1
       \or\or \hskip 1em \or \hskip 2em \else \hskip 3em \fi%
      #6\nobreak\relax
    \dotfill\hbox to\@pnumwidth{\@tocpagenum{#7}}\par
    \nobreak
    \endgroup
  \fi}
\makeatother

\begin{document}

\title[Cohomology with supports; idempotent pairs]%
{Cohomology with supports; idempotent pairs}

\begin{abstract} 
This chapter sets out preliminaries for the duality theory in later chapters.
An underlying idea is that \emph{local cohomology functors  are higher derived functors of 
colocalizations \(\kern-.7pt a.k.a.~coreflections\kf\).}

 Predominantly well-known facts about cohomology with supports---often 
 under ``\kf finitary" conditions that obtain, e.g., under noetherian hypotheses---and its 
 local and global interactions with quasi-coherence and with colimits, are reviewed from both the topological 
 and scheme\kf-theoretic perspectives. Some refinements of standard results are needed to accommodate certain features involving unbounded complexes and general systems of supports.

An important attribute of such cohomology 
 is   \mbox{``$\otimes$-coreflectiveness"}, in its avatar---ultimately in the context of closed 
 categories---as ``idempotent pair," a~notion which plays an important role in the sequel. 
 
Some basic facts about  linearly topologized noetherian rings and their maps, related to cohomology with supports, and subsumed under properties of idempotent pairs, are brought forth; and similarly for the less-familiar context of formal schemes.
 
\end{abstract}

\maketitle
\setcounter{tocdepth}{2}
\tableofcontents

 \newpage
\section{\bf \large Cohomology with supports; idempotent pairs\va2}

This chapter sets out some preliminaries for the duality theory in later chapters.
An underlying idea is that \emph{local cohomology functors  are higher derived functors of colocalizations (a.k.a.~coreflections\kf).}

 The first three sections review, from both the topological 
 and scheme\kf-theoretic perspectives (connected, as in \ref{ex:supports} and \ref{2RGams}), rudimentary facts about cohomology with supports---often under ``\kf finitary" conditions that obtain, in particular, under noetherian hypotheses---and its local and global interactions with quasi-coherence and with colimits 
 (see e.g.,  \ref{Dqc to itself} and~\ref{RGam and colim2}). A basic attribute of such cohomology is  
 \mbox{``$\otimes$-coreflectiveness"} (see~\ref{coreflections}). This property is  elaborated on 
 in the context of monoidal categories, as is its avatar ``idempotent pair," a~notion which plays
 an important role in the sequel (see Sections~\ref{tensorcompat}--\ref{idemcat}). 
 
 For instance,  if $X$ is a locally noetherian scheme and\/ $Z\subset X$ is closed, then  
$\R\vG{\<Z}\OX$ and the natural map $\iota\colon\R\vG{\<Z}\OX\to\OX$ form an idempotent $\D(X)$-pair,\va{.5} i.e.,  \smash{$\iota\Otimes{}\b1$} and~\smash{$\b1\Otimes{}\iota$} are \emph{equal isomorphisms} from\va{.5} \smash{$\R\vG{\<Z}\OX\<\Otimes{}\R\vG{\<Z}\OX$}\va1 to $\R\vG{\<Z}\OX$; 
and the corresponding $\otimes$-coreflection is
given by the functor\va1 \smash{$\R\vG{\<Z}(-)\set \R\vG{\<Z}\OX\OT{}(-)$} together\va1 with the map
\va1 \smash{$\iota\OT{} \b1\colon \R\vG{\<Z}(-)\to(-)$}. 
(See \ref{idpt3}.) 
 
 The idempotent pairs in a monoidal category~$\D$ are the objects of a strictly full monoidal subcategory $\mathbf I_\D$ 
 of the slice category $\D/\cO$ ($\cO\set$ unit object of $\D$);
$\mathbf I_\D$~is preordered, and the functor induced by the canonical functor $\D/\cO\to\D$ 
 is final in the category of all strong monoidal functors from preordered monoidal categories to $\D$
 (see Remark~\ref{I universal}).

Some basic facts about  linearly 
topologized noetherian rings and their maps, related to cohomology with supports, are subsumed under properties of idempotent pairs (Sections~\ref{ALC} and \ref{top+idem}); and similarly for 
formal schemes (Section~\ref{formal schemes0}).
In the latter case, if $\Dqct$ is the full subcategory of the derived category spanned by complexes with quasi-coherent torsion homology, then sending an idempotent pair in $\Dqct$ to its support gives an equivalence of $\mathbf I_{\Dqct}$ (modulo isomorphism) with the category of
inclusion maps of specialization-stable subsets (see~\ref{idem and specstable}).

This material is predominantly well-known (cf.~e.g.,~\cite[Expos\'es I, II\kf]{SGA2}), \cite{Hg}, \cite{AJS2}; but some refinements of the standard results are needed to accommodate certain features involving
unbounded complexes and general systems of supports. It is recommended to skim through these preliminaries,
referring back as needed in the subsequent duality theory. 

\setcounter{subsection}{-1}
\subsection{Terminology and notation}\label{notation}

\stepcounter{thm}
Let $\A$ be an abelian category. 

An $\A$-\emph{complex}
\mbox{$\>\>C=(C^\bullet\<,\>d^{\>\bullet})$} is  a sequence of
$\A$-maps 
$$
\cdots \xrightarrow{d^{\>i-2}} C^{\>i-1}
\xrightarrow{d^{\>i-1}} C^{\>i}
\xrightarrow{\ d^{\>i}\ } C^{\>i+1}
\xrightarrow{d^{\>i+1}}\cdots
\qquad(i\in\mathbb Z)\\
$$  
such that  $d^{\>i}d^{\>i-1}=0$ for all $i\>$. Homotopy equivalence of
maps of
$\A$-complexes is defined as usual  \cite[p.\,25]{H}.  The \emph {$i$-th cohomology}
$\>\Hr^iC\set\ker(d^{\>i})/{\rm im}(d^{\>i-1})$%
\footnote{Implicit here and elsewhere is the assumption that a \emph{specific choice} has been made in 
$\A$ of the kernel and cokernel of each $\A$-map, of a 0\kf-\kf object, of a direct sum for any two objects,\kern1pt\dots} is the object part of a natural
$\A$-valued functor on the category $\mathbf C(\A)$  of $\A$-complexes, or on
the homotopy category~$\K(\A)$ whose objects are $\A$-complexes and whose
morphisms are homotopy-equivalence classes  of maps of $\A$-complexes, or on
the derived category~$\D(\A)$ of~$\K(\A)$. (See e.g., \cite[\S\S1.1, 1.2]{Dercat}.) 

A \emph{quasi-isomorphism} in
$\mathbf C(\A)$ (resp.~$\K(\A)$) is a map of $\A$-complexes
$\phi\colon C\to C'$ which induces isomorphisms $\Hr^iC\iso\Hr^iC'$ for
all~$i$ (resp.~the homotopy equivalence class~$\bar\phi$ of such a map); or
equivalently, with
\mbox{$q_{\A}^{}\colon \K(\A)\to\D(\A)$} the canonical functor, such that 
$q_{\A}^{}\bar\phi$ is an isomorphism.
When the context dictates what is meant, there will usually be no
notational distinction among a map in $\mathbf C(\A)$,
its homotopy class in $\K(\A)$, and the image of that class under
$q^{}_{\A}\mkern.5mu$.\vs2

With reference to maps or diagrams in $\A$, $\K(\A)$ or
$\D(\A)$, ``natural'' means that unless otherwise specified,
the maps involved are the obvious ones.\vs2

Any  functor $\Gamma$ between triangulated categories (such as $\K(\A)$ or
$\D(\A)$) is understood to be \emph{additive} and
\emph{triangle-preserving}\va1 (a $\Delta$-functor, for short), i.e.,  equipped
with a functorial isomorphism\vs1
$\theta(E\>)\colon\Gamma(E[1])\iso(\Gamma E\>)[1]$ such that\vspace{1pt} for any
triangle
\smash{$E\underset{\under{3.2}{u\ }}\to F\underset{\under{3.2}{v\;}}\to
G\underset{\under{3.2}{w\ }}\to E[1]$,  the sequence\va{.5} 
$\Gamma E\underset{\under{3.5}{\Gamma u\ }}\lra 
\Gamma F\underset{\under{3.5}{\Gamma v\; }}\lra 
\Gamma G\xrightarrow{\raisebox{-1ex}{$\scriptstyle\<\theta\smcirc\Gamma
w\>$}}(\Gamma E\>)[1]$}  is also a triangle. In each instance the natural
definition of
$\theta$ is left to the reader.  (Sometimes there are sign considerations,
see, e.g., \cite[\S1.5]{Dercat} for more details, and for examples involving
$\otimes$ and Hom.) 
By definition, maps of $\Delta$-functors commute with the associated
$\theta\>$s.\vs1
\vspace{1pt}

A \emph{plump} subcategory
(or \emph{weak Serre subcategory}) $\ush\A\subset\A$ is a full subcategory containing 0 and such that for any exact $\A$-sequence $\>M_1\to M_2\to M\to M_3\to M_4,$  if 
$M_i\in\ush\A$ for $i=1,2,3,4$, then $M\in\ush\A$. The kernel and cokernel (in $\A$) of a map in such an $\ush\A$ both lie in $\ush\A$; so $\ush\A$ is abelian, and any object of $\A$ isomorphic to one in $\ush\A$ is itself in $\ush\A\<$.\va2

An $\A$-complex $I$ is \emph{K-injective} (\emph{q-injective} in the terminology of \cite{Dercat}, with%
``q" connoting ``quasi-isomorphism") 
if any quasi-isomorphism $\psi\colon I\to I'$ has a
left homotopy-inverse, that is, there exists an $\A$-homomorphism
$\psi'\colon I'\to I$ such that $\psi'\psi$ is homotopic to the identity map
of~$C$. Numerous equivalent\vspace{1pt} conditions can be found in
\cite[p.\,129, Prop.\,1.5]{Sp} and in
\cite[\S2.3]{Dercat}. One such is that the
functor~$\Homb(-,I\>)\colon\mathbf C(\A)\to\mathbf C(\A)$ preserves
quasi-isomorphism. Another is that for every $\A$-complex $F$, the natural
map is an isomorphism \va{-1}
\[
\Hom_{\K(\A)}(F,I\>)\iso\Hom_{\D(A)}(q_{\A}^{}F,q_{\A}^{}I\>).
\]
\vskip-1pt

Any bounded-below injective (in every degree) complex is K-injective.

A  \emph{K-injective resolution} of $E\in\A$ 
is a quasi-isomorphism $\sigma_{\!E}^{}\colon E\to I_{\<E}\>$
where $I_{\<E}$ is K-injective and also injective.\va2

\emph{All rings will be commutative.}\va2

For a  commutative ring $S$,  $\A(S)$ is the (abelian) category of small $S$-modules.
We write $q_S^{}\colon \K(S)\to \D(S)$ for
$q_{\A(S)}^{}\colon \K(\A(S))\to\D(\A(S))$. Each $S$-complex~$E$
admits a K-injective resolution $\sigma_{\!E}^{}\colon E\to I_{\<E}\>$,
see \cite[p.\,133, Prop.\,3.11]{Sp}. 

Similar considerations hold for any ringed space $(X,\cO\sX)$ 
($X$ a topological space and $\cO\sX$ a sheaf of commutative rings on~$X$), with $\A(X)$ the category
of \mbox{$\cO\sX$-modules,} with $q_X^{}\colon \K(X)\to \D(X)$ signifying
$q_{\A(X)}^{}\colon \K(\A(X))\to\D(\A(X))$,\vs{-1.5} etc.
(See \cite[p.\,138, Thm.\,4.5]{Sp}).%
\footnote{\kf Such assertions
hold in any Grothendieck category, see \cite[p.\,243, Thm.\,5.4]{AJS1}, or \cite[Propositions 1.3.5.3 
and 1.3.5.6]{Lu}, noting that K-injective $\Leftrightarrow$ homotopically equivalent to fibrant.}
$\Dpl(X)\subset\D(X)$ is the full subcategory spanned by the locally cohomologically bounded-below $\OX$-complexes (those $C\in\D(X)$ 
for which there is an open cover $(X_\alpha)_{\alpha\in A}$ of $X$ and  for each $\alpha$ 
an integer $n_\alpha$ such that the restriction 
$(H^iC)|_{X_\alpha}\<\<$ vanishes for all $i<n_\alpha\>$).
For such $(X,\cO\sX)$, restriction to open subsets preserves K-injectivity of $\OX$-complexes \cite[Lemma 2.4.5.2]{Dercat}.

For example, any topological space $X$ can be regarded as a ringed space, with $\OX$ the sheaf $\>\mathbb Z_X\<\<$ of locally constant functions from $X$ to $\mathbb Z\>$; and then $\A(X)$ is just the category~$\Ab (X)$ of sheaves of abelian groups. \vs2 

When $(X,\OX)$ is a scheme, $\Aqc(X)\subset\A(X)$ is the full subcategory of quasi-coherent $\OX$-modules, and $\Dqc(X)\subset\D(X)$ is the full subcategory whose objects are the complexes with quasi-coherent homology. \va2

\subsection{Finitary supports}\label{Prelims}

\begin{subcosa}\label{supports}
A \emph{system of supports} (s.o.s.) (a.k.a.~\emph{family of supports}) in a topological space~$X$ 
is a nonempty set $\Phi$ of closed subsets of~$\>X$
such that any closed subset of any finite union of members of $\Phi$ is a member of $\Phi$. \vs1

For instance, if $\>Y\<\subset X$ then the set~$\Phi_Y$ consisting of all subsets of~$Y$ that are closed in $X\<$ is an s.o.s.  An s.o.s.\kern2.24pt has this form if and only if it contains 
every \mbox{$X\<$-closed} subset of 
the union of all its members. In fact, there~is a one\kf-one correspondence between such s.o.s.\kern2pt and
\emph{specialization-stable} $Y\subset X\<$ (i.e.,  $Y$ contains the $X$-closure of each of its points, or equivalently, $Y$ is a union of closed subsets of $X$): 
to such an s.o.s.\kern2pt$\Phi$ corresponds the  union of its
members, and to a specialization-stable~$Y\subset X\<$ corresponds $\Phi_Y$.

If $X$ is \emph{noetherian,} i.e., every open subset is
quasi-compact \cite[II, \S4.2]{Bou},%
\footnote{A (possibly non-Hausdorff) topological space is \emph{quasi-compact} if every open cover has a finite subcover.}
then every closed subset of~$X$
is a finite union of irreducible closed subsets; and if,
furthermore, every irreducible closed subset of $X$ is the closure of one of its points 
(for instance, if $X$ is the underlying space of a noetherian scheme), then
every s.o.s.\ in~$X$ is $\Phi_Y\<$ for a unique specialization-stable~$Y\<$.\vs1

An s.o.s.~$\Phi$ in~$X$ is \emph{finitary} if each member of $\Phi$ is contained in a member $Z$ such that  $X\setminus Z$ is \emph{retrocompact} in~$X\<,$ 
i.e., for every quasi-compact open $U\subset X\<,$ the open subset
$U\setminus Z$ is quasi-compact.\vs1

For example, the s.o.s.~$\Phi_{\mkern-2.5mu X}$ consisting of all closed subsets of~$X$ is finitary.\va1

One checks that \emph{every s.o.s.~in $X$ is finitary $\Leftrightarrow$ every quasi-compact open subset of~$\>X\<$ is noetherian $\Leftrightarrow$ every open subset of~$X$ is retrocompact in} $X\<.$ 
(To see this, consider the s.o.s.~$\Phi_Y$ for an arbitrary closed $Y\subset X$\dots). 

These conditions on $X$ have no substance if the only quasi-compact open subset of $X$ is the empty one. 
More noteworthy is the situation where $X$ is a union of quasi-compact open subsets 
(for instance, the underlying space of a scheme): then the conditions
hold if and only if $X$~is \emph{locally noetherian,}  i.e., every point of $X$ has a noetherian neighborhood.

\end{subcosa}

\begin{slem}\label{finite intersection}
Let\/ $X$ be a quasi-compact topological space and let\/ $\Phi$ be a finitary
s.o.s.\,\,in\/~$X\<$. If\/ $(Z_\delta)_{\delta\in\mathscr D}$ is
a family of closed subsets of\/ $X$ such that\/ 
$\raisebox{2pt}{$\bigcap_{\>\delta\in\mathscr D\>}$}Z_\delta\in\Phi$ then there
is a finite subset\/ $\mathscr D_0\subset \mathscr D$ such that\/ 
$\raisebox{2pt}{$\bigcap_{\>\delta\in\mathscr D_0}\>$}Z_\delta\in\Phi$.
\end{slem}

\begin{proof} Fix
$Z\supset\raisebox{2pt}{$\bigcap_{\>\delta\in\mathscr D}\>$}Z_\delta$
such that $Z\in\Phi$ and $X\setminus Z$ is retrocompact in~$X\<,$ hence quasi-compact.
 The family $(Z_\delta\setminus Z)_{\delta\in\mathscr D}$ of closed subsets of $X\setminus Z$
has empty intersection, whence there is a finite 
$\mathscr D_0\subset \mathscr D$ such that 
$\raisebox{2pt}{$\bigcap_{\>\delta\in\mathscr D_0}\>$}(Z_\delta\setminus Z)$ is empty,  i.e.,~
$\raisebox{2pt}{$\bigcap_{\>\delta\in\mathscr D_0}\>$}Z_\delta\subset Z$, so that 
$\raisebox{2pt}{$\bigcap_{\>\delta\in\mathscr D_0}\>$}Z_\delta\in\Phi$.
\end{proof}

\begin{subcosa}\label{PhitoW} (Inverse image of an s.o.s.)
Let $f\colon W\to X$ be a continuous map of topological spaces,  and $\Phi$ an s.o.s.\kern2pt in~$X\<$. Set
\[
\begin{aligned}\label{Phi_f}
\Phi_{\!f}&\set\{\,V\textup{ closed in }W\mid \textup{the closure of } f(V)\textup{\ belongs to\ }\Phi\,\}\\
&\,=\{\,V\textup{ closed in }W\mid V\subset f^{-\<1}Z\textup{ for some }Z\in\Phi\,\}\\
&\,=\bigcup_{Z\in\Phi}\Phi_{f^{-\<1}Z}\>,
\end{aligned}\tag{\ref{PhitoW}.1}
\]
the smallest s.o.s.\ in $W$ that contains $f^{-\<1}Z$ for all $Z\in\Phi$.

For example, if $Y\subset X$ then $(\Phi_Y\mkern-1.5mu)_{\<f}\subset\Phi_{\!f^{-\<1}Y}$, with equality if
$Y$ is closed or if $Y$~is specialization-stable, $W$ is noetherian and every irreducible closed subset 
of~$W$ is the closure of one of its points.\va1

For another example, if $f$ is the inclusion map of a subspace $W\subset X$ then  
\begin{equation*}\label{Phi_fsub}
\Phi_{\!f}=\Phi|_W\set\{\,Z\cap W\mid Z\in\Phi\,\}.
\tag{\ref{PhitoW}.2}
\end{equation*}
Moreover, every s.o.s.~$\Phi_0$ in $W$ has the form $\Phi_{\!f}\colon$ let $\Phi$ consist~of 
all closed subsets of~$X$ whose intersection with $W$ is in $\Phi_0$.  \va2

The pairs $(X,\Phi)$ with $X$ a topological space and $\Phi$ an s.o.s.~in $X$ are the objects of a category in which a morphism $(W,\Psi)\to(X,\Phi)$ is a continuous map $f\colon W\to X$ 
such that $\Psi\subset\Phi_{\!f}$.
Such a morphism is called \emph{strict} if  $\Psi=\Phi_{\!f}$.

\begin{srem}\label{finitary restricts}
Let $W\subset X$ be open. If $\Phi$ is finitary then so is
$\Phi|_W$.
\end{srem}

\begin{srem}\label{soslocal}
Suppose that $Z$ is \emph{locally} in $\Phi$,\va{.6} i.e., 
$Z\subset \cup_{\alpha\in A}\>U_\alpha$ with each~$U_\alpha$  an open subset of $X\<,$ \va1 
such that $Z\cap U_\alpha\in \Phi|_{U_{\<\<\alpha}}$ (i.e., $\overline{Z\cap U_\alpha}\in\Phi$). 
If $A$ is finite, or if $\>\Phi=\Phi_Y$ ($Y\subset X$), then 
$Z\in\Phi$.\va3
\end{srem}
\end{subcosa}

\centerline{* * * * *}

\begin{subcosa}\label{basedef}
Let  $(X,\OX)$ be a scheme.  An $\OX$-\emph{base} is a  nonempty set~$\I$ of   
\emph{quasi-coherent} $\OX$-ideals such that: \vs2

(i) if $I\in\I$ and if $J$ is a quasi-coherent $\OX$-ideal such that\va1 $\sqrt{J}\supset I$, 
then $J\in\I\>$, and

(ii) if  $I\in\I$ and  $J\in \I$  then $I\cap J\in \I\>$.\va4

\noindent Since $\sqrt{IJ}\supset (I\cap J)\supset IJ$, therefore if  (i) holds\va1 then (ii) is equivalent to:\va2

(ii)$'$ if  $I\in\I$ and  $J\in \I$  then $IJ\in \I\>$.\va4

For example, if $\>\mathbf I$ is a nonempty set of $\OX$-ideals, 
and $f\colon W\to X$ is a map of schemes, then the smallest $\mathcal O_W$-base 
containing $I\mathcal O_W$ for all $I\in\mathbf I$ is\va2
\[
\begin{aligned}\label{I_f}
\mathbf I_{\<f}\set\{\,&\textup{quasi-coherent $\mathcal O_W$-ideals }J\mid\\
& \sqrt{J}\supset I_1\cdots I_n\mathcal O_W \textup{ for some integer } n\ge 0\textup{ and }I_1,\dots,I_n\in\mathbf I\,\}.
\end{aligned}
\tag{\ref{basedef}.1}
\]

When $f$ is the inclusion map of a subspace $W\subset X\<,$ $\mathbf I_{\<f}$ is 
denoted $\mathbf I|_W$.\va3

If $\I$ and $\J$ are $\OX$-bases, then so is $\I\cap\>\J=\{\,I+J\mid I\in\I,\:J\in\J\,\}$.\va2

An $\OX$-base $\I$ is \emph{finitary} if $X$ is covered by open subsets~$U$ such that each member 
of~$\I|_U$ contains a \emph{finite-type} member of~$\I|_U$.\va4
\end{subcosa}

\vskip3pt
\centerline{* * * * *}

\begin{subcosa}\label{suppdef}
Again, let $\<(X,\OX)\<$ be a scheme. 
\vspace{1pt}

The \emph{support of an $\OX\<$-\kf module~$M$} is\looseness=-1
\begin{equation*}\label{Suppdef}
\Supp(M\>)\set\{\,x\in X\mid M_x\ne (0)\,\}.
\tag{\ref{suppdef}.1}
\end{equation*}

For example, let $u\colon U\hookrightarrow X$ be the inclusion map of an open subscheme,\va{.6} and let $N$ be an $\OU$-module, with support (in $U$) $\Supp_U(N\>)$. Any point\va1 $x\in X$ lying outside the 
$X\<$-\kf closure $\overline{\Supp_U(N\>)}$ has a neighborhood\va1 in which $u_*N$ vanishes, 
and so\vs1 \mbox{$\>\Supp(u_*N\>)\subset \overline{\Supp_U(N\>)}$.} \va2

For any $s\in\Gamma(X,M\>)$, the \emph{support of} $s$ is
the closed set
\begin{equation*}\label{exc}
\textup{supp}(s)=\textup{supp}\sX\<(s)\set\{\,x\in X\mid s_x\ne 0\,\}=\Supp(s\OX)\subset\Supp(M\>).
\tag{\ref{suppdef}.2}
\end{equation*}

If $M$ is of finite type then $\Supp(M\>)$ is locally the union of the supports of members of a finite generating set, and so $\Supp(M\>)$ is a closed subset of $X$.\va 3

For an $\OX$-ideal $I$, the \emph{zero\kf -set of} $I$ is the closed set
\[
Z(I\>)\set\Supp(\OX/I\>).
\]

Every closed subset of $X$ is  $Z(I\>)$ for
some quasi-coherent $\OX$-ideal $I\<$. \va2

Clearly, $Z(I\>)=Z(\sqrt{I}\,)$  and $Z(I_1I_2)=Z(I_1)\cup Z(I_2)$. \va1

If $I_1$ and $I_2$ are quasi-coherent, 
one checks locally that
\begin{equation*}\label{radsup}
Z(I_1)\supset Z(I_2) \!\iff\! \sqrt{I_1}\subset\sqrt{I_2}.
\tag{\ref{suppdef}.3}
\end{equation*}

For $s\in\Gamma(X,M\>)$,
\begin{equation*}\label{supp=Zann}
\textup{supp}(s)=Z(\textup{ann}(s))
\end{equation*}
where ann$(s)$, the \emph{annihilator  of} $s$,  is the kernel of the $\OX$-homomorphism 
$\OX\to M$ taking $1\in\Gamma(X,\OX)$ to $s$.  \va2

\end{subcosa}

\begin{sprop} \label{ex:supports}
There is an inclusion-preserving bijection\/ $\mathcal S$
from the set of\/ $\OX\<$-\kf bases onto the set of systems of supports in a scheme\/~$X,$ 
such that for any quasi-coherent\/ $\OX\<$-\kf ideal\/~$I,$ $\OX\<$-\kf base~$\I$ 
and s.o.s.\/~$\Phi,$
\[\label{PhisubI}
 I\in \I \!\iff\! Z(I\>)\in\Phi_\I\set\mathcal S(\I\>),
 \tag{\ref{ex:supports}.1}
\] 
or equivalently, 
\[\label{IsubPhi}
 I\in \I_\Phi\set\mathcal S^{-1}(\Phi)\!\iff\! Z(I\>)\in\Phi.
\tag{\ref{ex:supports}.2}
\]

\end{sprop}

\begin{proof} 
Left to the reader.
\end{proof}

\begin{exam}\label{example2}
Let $f\colon W\to X$ be a scheme\kf-map, $\Phi$ an s.o.s.~in $X\<,$ $\Phi_{\!f}$ as in~\ref{PhitoW}, $\>\I_\Phi$ and $\I_{\Phi_{\!f}}$
as in (\ref{ex:supports}.2), and $\>(\I_\Phi\<)_{\<\<f}$~as in (\ref{I_f}).
For $I\in\I_\Phi$ and $J$ a quasi-coherent $\mathcal O_W$-ideal,
\eqref{radsup}~gives\va{-2}
\[
\{Z(J)\subset f^{-\<1}Z(I\>)=Z(I\OW)\}\ \!\iff\!\  
\{\sqrt{J}\supset I\mathcal O_W\},
\]
whence $\I_{\Phi_{\!f}}=(\I_\Phi\<)_{\<\<f}\>$.\va4
\end{exam}

\begin{scor}\label{baselocal}
Let\/ $\I$ be an\/ $\OX\<$-\kf base,\va{.5} and\/ $I$ an\/ $\OX\<$-\kf ideal
locally in\/~$\I,$ i.e., 
$X$~has an open covering $(U_\alpha)_{\alpha\in A}$ such that for each $\alpha$, 
$I\cO_{U_\alpha}\in \I|_{\<U_\alpha}.$ If\/~$A$ is finite, or if\/ $\Phi_\I=\Phi_Y$ for some\/ $Y\subset X$ \(see~\textup{\ref{PhisubI}, \ref{supports}),} 
then\/ $I\in\I$.
\end{scor}

\begin{proof} 
For all $\alpha\in A,$ the $\cO_{U_\alpha}\<\<$-ideal $I\cO_{U_\alpha}$ is quasi-coherent; so the $\OX$-ideal~$I$ is quasi-coherent.
Also, with $\Phi\set\Phi_\I,$ so that, by~\ref{ex:supports}, $\I=\I_\Phi\>$,
\[
I\cO_{U_\alpha}\in\I|_{\<U_\alpha}
\underset{\textup{\ref{example2}}}\implies I\cO_{U_\alpha}\in\I_{(\Phi|_{U_\alpha}\<)}
\underset{\textup{\ref{ex:supports}}}\implies
Z(I\cO_{U_\alpha})=Z(I\>)\cap U_\alpha\in \Phi|_{U_\alpha}.
\]
Remark~\ref{soslocal} ensures then that $Z(I\>)\in\Phi$, that is, $I\in\I$.
\end{proof}

Recall that the scheme $X$ is \emph{quasi-separated} if the intersection of any two quasi-compact open subsets is quasi-compact, see \cite[p.\,296, (6.1.12)]{GD}.

\begin{slem} \label{global finitary}
Let\/ $X$ be a quasi-compact quasi-separated scheme, and\/ $\I$ a finitary 
$\OX$-base. Every member of\/ $\I$ contains a finite-type member of\/ $\I$.
\end{slem}

\begin{proof}
 Let $I\in\I$. As $\I$ is finitary and $X$ quasi-compact, there exists a covering 
 $(U_i)\ (1\le i\le n)$ of~$X$ by a finite family of affine open subsets, 
 and for each $i$, a finite\kf-type  $J_i\in\I|_{U_i}$ with 
$J_i\subset I\cO_{U_i}$. Let $\bar J_i\subset I$ be a finite\kf-type $\OX$-ideal whose restriction\va{.6} to~$U_i$ is $J_i$ (see
 \cite[p.\,318, Thm.\,(6.9.7)]{GD}). Then $\bar J\set\sum_{i=1}^n \bar J_i\subset I$ is a quasi-coherent finite-type $\OX$-ideal whose restriction to each~$U_i$ contains $J_i$\kf, hence 
 lies in~$\I|_{U_i}$\kf; so by~\ref{baselocal}, $\bar J\in\I$. 
 \end{proof}
 
\begin{sprop} \label{finitary either way}
Let\/ $\I$ be an\/ $\OX\<$-\kf base. If\/ $\Phi_\I$ is finitary then\/ $\I$ is finitary. The converse holds
if $X$ is quasi-compact and quasi-separated.
\end{sprop}

\begin{proof}
Suppose $\Phi_\I$ finitary. For any open
$U\subset X$, $\Phi|_U$ is finitary. 
Hence to show that $\I$ is finitary,\va{.6}
one may assume that $X$ is affine, say $X=\Spec(R)$. Let $I\in\I$. Then 
$Z(I\>)\subset Z(\bar I\>)$ for some $\bar I\in\I$ such that $X\setminus Z(\bar I\>)$ is quasi-compact and so covered by finitely many open subsets 
$X\setminus Z(f_iR)$ with $f_i\in\Gamma(X,\bar I\>)\ (i=1,2,\dots,n)$. 
Since $\bar I\subset\sqrt{I}$ (see \eqref{radsup}),\va{.6} one can, upon replacing each $f_i$ by a suitable power, assume that every $f_i$ is in $\Gamma(X, I\>)$;\va1 and then by
\ref{basedef}(i), the ideal $(f_1,f_2,\dots,f_n)R$, whose radical contains 
$\Gamma(X,\bar I\>)$,\va{.6} sheafifies to a finite\kf-type ideal in~$\I$ that is contained in $I$. Thus $\I$ is indeed finitary.\va1

 For the converse, suppose $\I$  finitary and $X\va{.5}$ quasi-compact and quasi-separated.
 Let $Z\in\Phi_\I$, say $Z=Z(I\>)$ ($I\in\I$\kf). Let $(U_i)\ (1\le i\le n)$\va{1} and 
$\bar J$ be as in the proof of~\ref{global finitary},
so that $Z\subset Z(\bar J\>)\in\Phi_\I$.  $U_i$ being affine,\va{.5}  $\bar J\cO_{U_i}$ is generated by finitely many of its sections over $U_i\>$; so $U_i\setminus Z(\bar J\>)$, being an intersection of finitely many quasi-compact open sets, is quasi-compact, whence
$
X\setminus Z(\bar J)= \bigcup_{i=1}^{\>n} \big(U_i\setminus Z(\bar J\>)\big)
$ 
 is quasi-compact, hence retrocompact in $X\<$. 
 Thus $\Phi_\I$~is finitary.
\end{proof}

\medskip
\centerline{* * * * *}

\begin{subcosa}\label{gamph} 
Let  $(X,\OX)$ be a ringed space, and $M\in\A(X)$.
The support supp$(s)$ of $s\in\Gamma(X,M\>)$
is closed in~$X\<$ (see \eqref{exc}). 

For any s.o.s.~$\Phi$ in $X\<,$  and open $U\subset X$,
one has the $\Gamma(U\<,\OX\<)$-module
\[
\Gamma^{}_{\<\!\Phi}(U,M\>)\set\{\,s\in \Gamma(U,M\>)\mid\supp_U(s)\in\Phi|_U\}.
\]

Let $\vG{\Phi}$ be
the \emph{left-exact subfunctor of the identity functor on\/~$\A(X)$}\va{.6} such that
for any  $M\in\A(X)$, $\vG{\Phi}(M\>)$ is the sheaf associated to the presheaf
$\,U\mapsto \Gamma^{}_{\<\!\Phi}(U,M\>)$ $(U\textup{ open in }X),$ that is,
the sheaf of sections of~$M$\va{.6} whose support is locally in~$\Phi$. (See, apropos, Remark~\ref{soslocal}.)\va1

Following \cite[Expos\'e I, \S1]{SGA2}, for closed $Z\subset X$ we set $\Gamma^{}_{\<\<\!Z}\set\Gamma^{}_{\<\!\Phi_{\mkern-2.5mu Z}}$ and 
$\vG{\<Z}\set\vG{\Phi_{\mkern-2.5mu Z}}$ ($\Phi_{\!Z}$ consisting, as in \ref{supports}, of 
all closed subsets of $Z$).

Clearly, for any s.o.s.\ $\Phi$ and for $U\<$, $M$ as above,
\[
\begin{aligned}\label{Global Gam as lim}
\Gamma^{}_{\<\!\Phi}(U,M\>)=&\bigcup_{Z\in\Phi}\<\Gamma^{}_{\<\<\!Z}(U,M\>)
=\dirlm{\lift .25,{\scriptstyle Z\<\in\<\Phi^{}},\,}\<\<\Gamma^{}_{\<\<\!Z}(U,M\>),
\\
\vG{\Phi}(M\>)=&\bigcup_{Z\in\Phi}\<\vG{\<Z}(M\>)
=\dirlm{\lift .25,{\scriptstyle Z\<\in\<\Phi},\,}\<\<\vG{\<Z}(M\>).\\[4pt]
\end{aligned}
\tag{\ref{gamph}.1}
\]
\vspace{-9pt}

As in \cite[Expos\'e I, 1.6]{SGA2}, the functor $\Gamma^{}_{\<\<\!Z}(U,M\>)$ (denoted there by $\Gamma^{}_{\<\<\!Z\cap\>U}(M\>)$) is naturally isomorphic to  
$\Hom_{\Ab(U)}\<\<(\mathbb Z_{Z\cap \>U\!,\>U},M|_U)$,
where  
$\mathbb Z_{Z\cap \>U\!,\>U}$~is the abelian sheaf on $U$ which restricts over $Z\cap U$ to the locally constant sheaf of integers $\mathbb Z$ and vanishes elsewhere.
Hence there~is a functorial isomorphism of $\OX$-modules 
\begin{equation*}\label{Gam and Hom}
\vG{\<Z}(M\>)\cong \sHom_{\>\Ab\<(X)}\<\<(\mathbb Z_{Z,X},M\>). 
\end{equation*}

\vskip2pt
The functor $\vG{\Phi}$ is \emph{idempotent}: $\vG{\Phi}\vG{\Phi}=\vG{\Phi}\>$.
In fact, if each of $\Phi$ and $\Psi$ is an s.o.s. in $X$ then so is $\Phi\cap\>\Psi$, and one checks that
\begin{equation*}\label{intersect vG}
\vG{\Phi}\vG{\Psi}=\vG{\Phi}\cap\>\vG{\Psi}=\vG{\Phi\>\cap\>\Psi}\>.
\tag{\ref{gamph}.2}
\end{equation*}
And if $U$ is  an open subset of~$X$ such that every member of $\Psi|_{\<U}$ is quasi-compact (for instance, if $\>U$ itself is quasi-compact), 
or if $\>\Psi=\Phi_Y$ for some $Y\subset X\<,$ then using Remark~\ref{soslocal} one checks that
\begin{equation*}\label{intersect sos}
\Gamma^{}_{\!\<\Phi}(U,\vG{\Psi} M\>) = \Gamma^{}_{\<\!\Phi\>\cap\>\Psi}(U,M\>).
\tag{\ref{gamph}.3}
\end{equation*}

\vskip2pt

Let $f\colon W\to X$ be a continuous map of topological spaces,
$\Phi$ an s.o.s.\kern2pt in~$X\<,$ and as in \ref{PhitoW},\va{-1} 
$\Phi_{\<\!f}\set\{\,V\textup{ closed in }W\mid  \textup{the closure of } f(V)\textup{\ belongs to\ }\Phi\,\}.$
In particular, if $Y\subset X$ is closed and $\Phi=\Phi_Y$  then $\Phi_{\<\!f}=\Phi_{\!f^{-\<1}Y}$.

It is straightforward to see                                                                                                                                                                                                                  that for any $N\in\A(W)$, the support of a global section of $\fst N$ is the closure of the image under $f$
of the support of the corresponding global section of $N\<$. It follows that
\begin{equation*}\label{Gam and f}
\mkern14mu \Gamma^{}_{\<\!\Phi_{\<\!f}}\<\<(W,N)=\Gamma^{}_{\<\!\Phi}(X,\fst N)
\qquad(N\in\A(W)).
\tag{\ref{gamph}.4}
\end{equation*}

If  $X$ has a base of open sets $U$ such that
$f^{-\<1}U$ is quasi-compact, or if \mbox{$\Phi=\Phi_Y$} with $Y\subset X$ closed, 
then~\eqref{intersect sos} (with $(\Phi,U,\Psi)$ replaced by $(\Phi_W, f^{-\<1}U,\Phi_{\<\!f})$) and~\eqref{Gam and f} (with $X$ replaced by an arbitrary $U$ and $W$ by $f^{-\<1}U$) give\begin{equation*}\label{vGam and f}
\fst\vG{\Phi_{\!f}}= \vG{\Phi}\fst\> .
\tag{\ref{gamph}.5}
\end{equation*} 

\end{subcosa}

\medskip
\centerline{* * * * *}

\begin{subcosa}\label{gambase}
Let $X$ be a scheme, $U\subset X$ open, $\I$ an $\OX$-base, $M$ an $\OX$-module and 
\[
\Gamma^{}_{\<\!\I}(U,M\>)\set
\dirlm{\lift.01,{{\lift.9,\halfsize{$I$},\in\>}\lift.95,\halfsize{$\I$},},\,}\<\Hom_{\OU}\<\<\big(\OU/(I|_{U}),M|_{U}\big).
\]
\vskip4pt
There is a natural isomorphism
\[
\Gamma^{}_{\<\!\I}(U,M\>)\iso\{\,s\in\Gamma(U,M\>)\mid \textup{ann}_{\>U}(s)\supset I|_U
\textup{ for some }I\in\I\,\}.
\]

There is an obvious presheaf $U\mapsto\Gamma^{}_{\<\!\I}(U,M\>)$. The associated sheaf is
\begin{equation*}\label{defvg}
\vG{\I}(M\>)\set\dirlm{\lift.01,{{\halfsize{$I$}\in\>}\lift.95,\halfsize{$\I$},},\,}\< \sHom_{\OX}\<\<(\OX/I,M\>)
\subset M.
\tag{\ref{gambase}.1}
\end{equation*}
\vskip4pt
\noindent There results a \emph{left-exact subfunctor\/ $\vG{\I}\colon\A(X)\to\A(X)$ of the identity functor}.\va1

For a quasi-coherent $\OX$-ideal~$I\<$,  let $\I_{\<\<I}$
be the $\OX$-base consisting of all quasi-coherent $\OX$-ideals whose radical contains $I$, i.e., the smallest $\OX$-base containing~$I$. (According to \eqref{I_f}, this is $\{I\}_{\>\mathbf 1_{X}}$.)\looseness=-1

Set $\Gamma^{}_{\!\<\<I}\set \Gamma^{}_{\<\!\I_{\!I\<}}$ and 
$\vG{I}\set \vG{\I_{\!I\<}}$. Then for any $\OX$-base~$\I$,\va{-3}
 \begin{equation*}\label{Global Gam as lim1}
\Gamma^{}_{\<\!\I}(U,M\>)=\bigcup_{I\lift1,\in,\>\I}\<\Gamma^{}_{\<\<\!I}(U,M\>)=
\dirlm{\lift .01,{{\halfsize{$I$}\lift.8,{\scriptscriptstyle\in},\>}\lift1,\halfsize{$\I$},},\,}\<\Gamma^{}_{\<\<\!I}(U,M\>),\\[-3pt]
\tag{\ref{gambase}.2}
\end{equation*}
\vskip-3pt
\noindent whence\va{-2}
\begin{equation*}\label{Gam as lim1}
\vG{\I}(M\>)
=\dirlm{\lift .01,{{\halfsize{$I$}\lift.8,{\scriptscriptstyle\in},\>}\lift1,\halfsize{$\I$},},\,}
\<\vG{I}(M\>).
\tag{\ref{gambase}.3}
\end{equation*}

\medskip
If the open $U\subset X$  is \emph{quasi-compact} then for any $\OX$-bases $\I$ and $\J$,\va{-1}
\begin{equation*}\label{intersect bases}
\Gamma^{}_{\!\<\I}(U,\vG{\J} M\>) = \Gamma^{}_{\<\!\I\cap\>\J}(U,M\>).
\tag{\ref{gambase}.4}
\end{equation*}
\vskip-1pt\noindent
Indeed,  for any $s\in \Gamma^{}_{\!\<\I}(U,\vG{\J} M\>)\subset \Gamma(U,M\>)$ there is a finite open cover $U=\cup_{i=1}^nU_i$ such
that for~each~$i$, the restriction $s|_{U_i}$ is annihilated by (the restriction of) some $J_{i}\in\J$; and then for some $I\in\I$, $s$ is annihilated by  $I+J_1J_2\cdots J_n\in\I\cap\>\J$. Thus
$\Gamma^{}_{\!\<\I}(U,\vG{\J} M\>) \subset \Gamma^{}_{\<\!\I\cap\>\J}(U,M\>)$;\va1 and the opposite inclusion is clear. 

As $X$ has a base of quasi-compact open sets, sheafifying shows then that\va{-2}
\begin{equation*}\label{local intersect sos}
\vG{\I}\vG{\J}=\vG{\I\cap\>\J}\>.
\tag{\ref{gambase}.5}
\end{equation*}
In particular (set $\J\set\I$), the functor $\vG{\I}$ is \emph{idempotent}.\va2
\end{subcosa}

\begin{subcosa}
Let $f\colon (W,\OW)\to (X,\OX)$ be a map of schemes,\va{.5} and $N\in\A(W)$. A section  
\mbox{$s\in\Gamma(X,\fst N)=\Gamma(W,N)$} can be regarded\va1 as the $\OX$-homomorphism
$\mathbf s\colon\OX\to\fst N$ that takes $1\in\Gamma(X,\OX)$ to $s$, 
or as the natural composite $\OW$-homomorphism
$\bar{\mathbf s}\colon\OW=f^*\OX\xto{\<f^*\mathbf s\>}f^*\<\<\fst N\lra N$ 
\va1(taking $1\in\Gamma(W,\OW)$ to~$s$).\va2

Let $\I$ be an $\OX$-base. Complying with \ref{I_f}, set
\[
\I_{\!f}\set\<\{\textup{\kf\kf quasi-coherent $\mathcal O_W$-ideals }J\mid \sqrt{J}\supset  I\mathcal O_W \textup{ for some }I\in\I\,\}.
\] 
If  $I\subset\ker(\mathbf s)=\textup{ann}_{X}(s)$, then 
$I\OW\subset\ker(\bar{\mathbf s})=\textup{ann}_{W}(s)$, whence
\[
\Gamma^{}_{\!\I}(X,\fst N)\subset\Gamma^{}_{\mkern-3mu\I_{\!f}}\<\<(W,N).
\]
Furthermore, \emph{if~$X$ is quasi-compact and quasi-separated and $\I$ is finitary} then by~\ref{global finitary}, 
one can assume that in the definition of $\I_{\!f}$,  the ideal $I$ is of finite\- type,                                                                                                                                                                                                                                                                                                                                                                                                                                                                                                                                                                                                                                                                                                                                
so one can replace $\sqrt{J}$ by $J$. Thus if $s\in\Gamma^{}_{\mkern-3mu\I_{\!f}}\<\<(W,N)$,
then there is an
$I\in\I$ such that $ I\OW\subset \ker(\bar{\mathbf s})$, whence the top row in the natural commutative diagram\va{-2}
\[
\def\1{$\fst f^*I$}
\def\2{$\fst \OW$}
\def\3{$\fst N$}
\def\4{$I$}
\def\5{$\OX$}
 \bpic[xscale=1.75, yscale=1.3]
  \node(11) at (1,-1){\1} ;
  \node(12) at (2,-1){\2} ;
  \node(13) at (3,-1){\3} ;
  \node(21) at (1,-2){\4} ;
  \node(22) at (2,-2){\5} ;
  \node(23) at (3,-2){\3} ;
 
   \draw[->] (11) -- (12) ;
   \draw[->] (12) -- (13) node[above, midway, scale=.75]{$\fst\bar{\mathbf s}$}  ;
   
   \draw[->] (21) -- (22) ;
   \draw[->] (22) -- (23) node[below=1, midway, scale=.75]{$\mathbf s$}  ;
   
   \draw[->] (21) -- (11) ;   
   \draw[->] (22) -- (12) ;
   \draw[double distance=2] (23) -- (13);
    
 \epic
\]
\vskip-3pt\noindent  
composes to 0, whence so does the bottom row, and so $s\in\Gamma^{}_{\!\I}(X,\fst N).$ 
Hence 
\[
\Gamma^{}_{\!\I}(X,\fst N)=\Gamma^{}_{\mkern-3mu\I_{\!f}}\<\<(W,N).
\]

From this plus \eqref{intersect bases},
it follows---without $X$ having to be quasi-compact and quasi-separated---that 
\emph{if the map $f$ is quasi-compact and $\I$ is finitary then}\va{-1} 
\[
\vG{\I}\fst=\fst\vG{\I_{\!\<f}}.
\]
\vskip1pt
\end{subcosa}

\begin{sprop}\label{Gam and qc} 
Let\/ $X$ be a scheme, $\I$ a finitary\/ $\OX\<$-\kf base,
and\/ $M$ a quasi-coherent\/ $\OX\<$-\kf module. Then the $\OX\<$-\kf module\/ $\vG{\I}M$ is quasi-coherent. 
\end{sprop}

\begin{proof}
The assertion being local (see \ref{finitary restricts}), $X$ can be assumed affine, so that 
every member of $\I$ 
contains a finite\kf-type member (see \ref{global finitary}). The assertion follows then from
\eqref{defvg} and \cite[p.\,217, (2.2.2)]{GD}. 
\end{proof}

\begin{sprop}\label{2Gammas} 
Let\/ $X$ be a scheme,  $M$ an $\OX\<$-\kf module, $\Phi$ an s.o.s.~in $X$
and\/ $\I\set\I_\Phi$ \textup{(see \ref{ex:supports})}.
Then $\Gamma^{}_{\<\!\I}(X,M\>)\subset\Gamma^{}_{\!\<\Phi}(X,M\>)$
and $\vG{\I}M\subset\vG{\Phi}M,$ with equality \(in~either case\) if\/ $M$ is quasi-coherent.\footnotemark
\footnotetext{\kf For examples of inequality, with $X$ noetherian and $M$ injective, see the proof of \ref{Gam preserves injective}.}
\end{sprop}

\begin{proof} Any $s\in\Gamma^{}_{\<\!\I}(X,M\>)$ is annihilated by an $I\in\I$, whence supp$(s)\subset Z(I)\in\Phi$, that is, $s\in\Gamma^{}_{\!\<\Phi}(X,M\>)$.  
Thus $\Gamma^{}_{\<\!\I}(X,M\>)\subset\Gamma^{}_{\!\<\Phi}(X,M\>).$\va2

If, moreover,  $M$~is quasi-coherent, then so is ann$(s)$ for any $s\in\Gamma(X,M\>)$, and 
$
s\in\Gamma^{}_{\<\!\Phi}(X,M\>)\<\<\!\iff\!\<\< \supp(s)=Z(\textup{ann}(s))\in\Phi
\<\<\!\iff\!\<\< \textup{ann}(s)\in\I\<\<\!\iff\!\<\< s\in\Gamma^{}_{\<\!\I}(X,M\>),
$ 
so that $\Gamma^{}_{\!\<\I}(X,M\>)=\Gamma^{}_{\!\<\Phi}(X,M\>)$. \va2

Replacing $X$ by an arbitrary open subset, one gets inclusion (resp.~equality) for the resulting presheaves,
and sheafification gives inclusion (resp.~equality) for~$\vG{}\<$.
\end{proof}

The next result is immediate from  \ref{Gam and qc} and \ref{2Gammas}. (See  also \ref{Dqc to itself} below for an essentially well-known generalization.)
\begin{scor}\label{Cor2Gammas}
Let\/ $X$ be a scheme and\/ $\Phi$ an s.o.s.~in\/~$X\<$.
If\/ $M$~is a quasi-coherent\/ $\OX\<$-\kf module then so is\/ $\vG{\Phi}M$. 
\end{scor}
 \vs3

\noindent\emph{Remark.} In \cite[p.\,2293]{GS} there is an example in which $X$ is 
the spectrum of a
polynomial ring in countably many variables over a field, $I$ is the sheafification of the ideal
generated by the variables, and $M$ is a certain quasi-coherent $\OX$-module such that  $\vG{I}(M\>)$ is not quasi-coherent. (There, of course, $\I_{\<\<I}$ is not finitary.)

\vskip13pt
\centerline{* * * * *}

\begin{sprop} \label{Gam and lim} 

\textup{(i)}  Suppose that the topological space\/ $X$ has a base of quasi-compact open sets, and that the s.o.s.\kern2.4pt$\Phi$ in\/ $X$ is finitary. Then 
$\vG{\Phi}$~commutes with small filtered colimits, hence
with small direct sums. 
                                                                            
More exactly, if~$A$ is a small filtered category\/ \textup{\cite[p.\,211]{Ma2}} 
and\/  \mbox{$\mathscr M\colon A\to\Ab(X)$}  is a~functor, then
the natural map is an isomorphism
\[
\lambda\colon\dirlm{A}\<\<(\vG{\Phi}\<\smcirc \mathscr M\>)\iso 
\vG{\Phi}(\>\dirlm{A}\<\< \mathscr M\>).
\]
\vskip2pt
\textup{(ii)} Let\/ $X$ be a scheme and\/ $\I$ a finitary\/ $\OX\<$-\kf base. 
Then $\vG{\I}$~commutes with small filtered colimits \(as in \textup{(i)}\)$,$ hence with small direct sums. 
\end{sprop}

\begin{proof} (i). Since the composite map\va{-2} 
\[
\dirlm{A}\<\<(\vG{\Phi}\<\smcirc \mathscr M\>)
\xto{\;\lambda\;}
\vG{\Phi}(\dirlm{A}\<\< \mathscr M\>)
\xto{\<\<\textup{natural}\>}\dirlm{A}\<\< \mathscr M
\]
is the natural injection, therefore $\lambda$ is injective.

Surjectivity can be checked stalkwise. Fix $x\in X\<$. Any element of  
$(\vG{\Phi}(\dirlm{} \<\<\mathscr M\>))_x$\vs{-1} is the germ $\sigma_{\<\<x}$ of a section $\sigma$ of $\dirlm{} \<\<\mathscr M$ over a quasi-compact open neighborhood~$V$ of $x$, such that \va{.6} $\sigma$  is the natural image of a section 
$\sigma_{\<\<a}\in\Gamma(V\<\<,\>\mathscr Ma)$ for some $a\in A$, and   
$\supp(\sigma)  \in\Phi|_V$. 

\pagebreak[3]

For each $A$-morphism $\alpha\colon a\to b$, let $\sigma_{\<\<\alpha}$ be the image of $\sigma_{\<\<a}$
under the induced map 
$\Gamma(V\<\>,\>\mathscr M\alpha)\to \Gamma(V,\mathscr Mb)$.
Then $\sigma$ is  the natural image of~$\sigma_{\<\<\alpha}\>$; and for all $y\in V\<\<$, 
$\sigma_{\<\<y}\ne0 \!\iff\!  (\sigma_{\<\<\alpha})_y\ne0\textup{ for all }\alpha,$
i.e., $\raisebox{2pt}{$\scriptstyle\bigcap$}^{}_\alpha \>\>\supp(\sigma_{\<\<\alpha})
=\supp(\sigma)\in\Phi|_V.$
Since $V$ is quasi-compact and $\Phi|_V$~is finitary, and since  $A$~is filtered, 
 \Lref{finite intersection} implies that there exists a single~$\alpha\colon a\to b$ with~ 
$\supp(\sigma_{\<\<\alpha}) \in\Phi|_V$. For such an $\alpha$,  
$(\sigma_{\<\<\alpha})_x$ is an element of $(\vG{\Phi}(\mathscr Mb))_x$ whose natural image
in $(\dirlm{}\!(\vG{\Phi}\<\smcirc\mathscr M\>))_x$ is taken by $\lambda_x$ to $\sigma_x$.
Thus $\lambda_x$ is surjective for any $x\in X$, that is, $\lambda$ is surjective.\va1
 
The passage from filtered direct limits to direct sums is standard (cf.~the last part of the proof of \Pref{RGam and colim1} below).
 \vs2

(ii) The assertion being locally verifiable, one can assume $X$ affine. Lemma~\ref{global finitary} shows then that every member of the 
finitary $\OX$-base contains a finite\kf-type one.  Hence the assertion is given by the natural isomorphisms, with $\I_0$ consisting of all finite\kf-type $I\in\I,$\va{-2}
\begin{align*}
\dirlm{\lift.2,\sst A,}\<\<(\vG{\I}\<\smcirc \mathscr M\>)
=\:
&\dirlm{\lift.2,\sst A,}\dirlm{\lift .3,{{\halfsize{$I$}\lift.85,{\scriptscriptstyle\in},\>}\lift1,\halfsize{$\I_0$},},\,}\<\<\big(\>\sHom_{\OX}\mkern-1.5mu(\OX/I, -)\smcirc \mathscr M\big)\\[4pt]
\iso
&\dirlm{\lift .3,{{\halfsize{$I$}\lift.85,{\scriptscriptstyle\in},\>}\lift1,\halfsize{$\I_0$},},\,}
\dirlm{\lift.2,\sst A,}\<\<\big(\sHom_{\OX}\mkern-1.5mu(\OX/I, -)\smcirc \mathscr M\big)\\[4pt]
\iso 
&\dirlm{\lift .3,{{\halfsize{$I$}\lift.85,{\scriptscriptstyle\in},\>}\lift1,\halfsize{$\I_0$},},\,}\<\sHom_{\OX}\mkern-1.5mu\big(\OX/I, \>\> \dirlm{\lift.15,\sst A,}\<\mathscr M\big)
\iso
\vG{\I}(\>\dirlm{\lift.15,\sst A,}\< \mathscr M).\\[-34pt]
\end{align*}
\end{proof}

As in \cite[p. 640]{Kf}, a \emph{quasi-noetherian} topological space is one that is quasi-compact and 
has a base of quasi-compact open subsets any two of which have quasi-compact intersection.

For example, the underlying space of a quasi-compact quasi-separated scheme or formal scheme is quasi-noetherian.

\begin{scor}\label{lims of Gams} 

\textup{(i)} Let\/ $X$ be a quasi-noetherian topological space, 
 $\Phi$ a finitary s.o.s.~in\/ $X\<,$ and\/ $\mathscr M$ as in \textup{\ref{Gam and lim}}. The natural map is an isomorphism
$$
\dirlm{A}\<\big(\Gamma^{}_{\!\<\<\Phi}(X, -)\smcirc\mathscr M\big)\iso 
\Gamma^{}_{\!\<\<\Phi}(X,\>\> \dirlm{A} \mathscr M\>).
$$
In particular, $\Gamma^{}_{\!\<\<\Phi}(X, -)$ commutes with direct sums.\va1

\textup{(ii)} Let\/ $X$ be a quasi-compact quasi-separated scheme,  
$\I$ a finitary\/ $\OX\<$-\kf base, and\/~$\mathscr M$ as in \textup{\ref{Gam and lim}}.
The natural map is an isomorphism
$$
\dirlm{A}\<\big(\Gamma^{}_{\!\<\<\I}(X, -)\smcirc\mathscr M\big)\iso 
\Gamma^{}_{\!\<\<\I}(X,\>\> \dirlm{A} \mathscr M\>).
$$
In particular, $\Gamma^{}_{\!\<\<\I}(X, -)$ commutes with direct sums.
\end{scor}

\begin{proof}
Let $\bullet$ denote one of $\Phi$ and $\I$. As $\,\dirlm{}\<$ commutes\va{-.5} with
$\Gamma(X,-)$ (see~\cite[p.\,641, Prop.\,6]{Kf}), one~gets, by~setting, in \eqref{intersect sos},
$\Phi \!\set\,$\{all closed subsets of~$X$\kf\},\va{.5}  or by switching, in~\eqref{intersect bases},
 $\I$ and~$\J$  and then setting $\J\set$\{\kf all quasi-coherent $\OX$-ideals\kf\}, natural iso\kf\-morphisms
\begin{flalign*}
\qquad\dirlm{A}\<\big(\Gamma_{\!\<\<\bullet}(X, -)\smcirc\mathscr M\big)\>
&{\iso}\,
\dirlm{A}\<\big(\Gamma(X,-)\smcirc\varGamma_{\!\bullet}\smcirc\mathscr M\big)
\\[3pt]
&\iso
\Gamma(X,-)\smcirc\>\dirlm{A}\<\<(\varGamma_{\!\bullet}\smcirc\mathscr M\>)
\\[3pt]
&\underset{\ref{Gam and lim}^{\mathstrut}}\iso
\Gamma(X,-)\smcirc\>\varGamma_{\!\bullet}(\>\dirlm{A}\< \mathscr M\>)\>
\iso
\Gamma_{\!\<\<\bullet}(X,\>\> \dirlm{A}\< \mathscr M\>), 
\end{flalign*}
whose composition is the map in question.
\end{proof}
\vskip1pt

\subsection{Cohomology with supports: topological spaces and schemes}\label{CohSupp}

\stepcounter{thm}
 Next, the derived functors of those just considered. Notation remains as in~\Sref{notation}.

Let $(X,\OX)$ be a ringed space. 
An additive functor ${\mathcal G}\colon\A(X)\to\A$
extends naturally\va{.4} to a functor $\overline{\mathcal G}\colon \K(X)\to\K(\A)$.
Given, for each $\OX$-\kf complex~$E$,  a \mbox{K-injective} resolution\va{.2}
$\sigma_{\!E}^{}\colon E\to I_{\<E}\>$, with homotopy class~$\tilde\sigma_{\!E}^{}\>$, 
there~exists a right-derived functor
\mbox{$\R{\mathcal G}\colon\D(X)\to\D(\A)$} and a~functorial map
$\zeta^{}_{\mathcal G}\colon q^{}_{\<\A}\overline{\mathcal G}  \to\R{\mathcal G}q^{}_{\<X}\>$  
such that for all~$E$,\va{.5} 
$\R{\mathcal G}q^{}_{\<X}E=q^{}_{\<\A}\overline{\mathcal G}\< I_{\<E}$  and 
$\zeta^{}_{\mathcal G}(E\>)=q^{}_{\<\A}\overline{\mathcal G}\tilde\sigma_{\!E}^{}\>$. 
(See, e.g., \cite[\S2.3]{Dercat}.) 
To a functorial map $\lambda\colon{\mathcal G}\to\mathcal G'\<$,\va1 with natural extension 
$\bar\lambda\colon\overline{\mathcal G}\to\>\>\overline{\mathcal G'}$,  there is associated
a unique functorial map\va{.5} $\R\lambda\colon\R{\mathcal G}\to\R{\mathcal G}'$ such that 
$\R\lambda\smcirc\zeta^{}_{\mathcal G}=\zeta^{}_{\mathcal G'}\<\smcirc\bar\lambda\>$.
\va1

For instance, with $\Phi$ an s.o.s.\ in $X$ and $U\subset X$ open, let ${\mathcal G}$ be the functor
\[
\Gamma^{}_{\!\<\<\Phi}(U,-)\colon\A(X)\to\A(\Hr^0(U,\OU)).
\]
Let $u\colon U\hookrightarrow X$ be the inclusion. Then $u^*\!$ takes
K-injective resolutions to K-injective resolutions, and therefore one has, with 
$\Phi|_U$ as in \eqref{Phi_fsub}, a natural isomorphism of functors (from $\D(X)$ to $\D\big(\Hr^0(U,\OU)$): 
$\R\Gamma^{}_{\<\<\!\Phi}(U,-)\iso\R\Gamma^{}_{\<\<\!\Phi|_U}\<\<(U,-)\smcirc u^*\<$.

Likewise, if $(X,\OX)$ is a scheme and $\I$ an $\OX$-base then one has, with $\I|_U$ as in 
the line following \eqref{I_f}, a natural isomorphism 
$\R\Gamma^{}_{\<\<\!\I}(U,-)\iso\R\Gamma^{}_{\<\<\!\I|_U}\<\<(U,-)\smcirc u^*\<$.\va1

\enlargethispage*{10pt}
\emph{The ordered set\/ $\Phi$ \(resp.~$\I\>$\) will always be regarded as a filtered category,
with inclusions \(resp.~containments\) as morphisms.}

\begin{subcosa}\label{derived}
With preceding notation, set $H^n_{\<\Phi}\set H^n\R\vG{\Phi}\ (n\in\mathbb Z),$ and
$H^n_{\<\<Z}\set H^n\R\vG{\<\Phi_{\!Z}}$. One has natural functorial isomorphisms
\begin{equation*}\label{HGam as lim.1}
H^n_{\<{\Phi}} E
\cong H^n\<\vG{\Phi}\>I_{\<E}
\underset{\eqref{Global Gam as lim}}\cong 
H^n\>\>\dirlm{\lift .15,{\!\!\lift.9,\halfsize{$Z$},\in\lift1,\Phi,},\!\!}\<\!\vG{\<Z}I_{\<E}
\cong\dirlm{\lift .15,{\!\!\lift.9,\halfsize{$Z$},\in\lift1,\Phi,},\!\!}\!H^n\<\vG{\<Z}I_{\<E}
\cong\dirlm{\lift .15,{\!\!\lift.9,\halfsize{$Z$},\in\lift1,\Phi,},\!\!}\!H^n_{\!Z}E\>.\tag{\ref{derived}.1}
 \end{equation*}
 \vskip2pt
Similarly, if $(X,\OX)$ is a scheme and $\I$ an $\OX$-base, then with
\[
H^n_{\<\I}\<E\set H^n\R\vG{\I}E=\dirlm{\lift .25,{\lift.9,\halfsize{$I$},\in{\lift1,\I,}},\,}
  \<\<\mathcal{E}xt^n(\OX/I,E)
\qquad (n\in\mathbb Z)
\]
\vskip3pt\noindent
(set $M\set I_{E}$ in \eqref{defvg}), with $\vG{I}$ as in the lines preceding~\eqref{Global Gam as lim1}, 
and \mbox{$H^n_{\<I}E\set H^n\R\vG{I}E,$}
one has natural functorial isomorphisms 
\begin{equation*}\label{HGam as lim.2}
H^n_{\<\I} \mkern-1.5mu E
\cong H^n\<\vG{\I}I_{E}
\underset{\eqref{Gam as lim1}}\cong 
H^n\>\>\dirlm{\lift.25,{\lift.9,\halfsize{$I$},\in{\lift1,\I\>,}},\,}\!\vG{\<I}I_{\<E}
\cong\dirlm{\lift .25,{\lift.9,\halfsize{$I$},\in{\lift1,\I\>,}},\,}
  \<\<H^n\<\vG{\<I}I_{\<E}
\cong\dirlm{\lift .25,{\lift.9,\halfsize{$I$},\in{\lift1,\I\>,}},\,}
  \<\<H^n_{\<I}\<E.
  \tag{\ref{derived}.2}
 \end{equation*}
 \vskip1pt\noindent
 
 Ditto,  via  \eqref{Global Gam as lim} or \eqref{Global Gam as lim1}, with 
$\bullet\set\Phi$ or $\I$, for $\Hr^n_{\>\bullet}(U,E\>)\set \Hr^n\R\Gamma_{\!\<\bullet}(U,E\>)$\quad ($U$ open in $X$).\vs2
\end{subcosa}

\begin{sprop}\label{Dqc to itself} 
If\/ $\Phi$ is a finitary s.o.s.\kern2.4pt in a scheme\/ $X\<,$ then 
\[
\R\vG{\Phi}\Dqc(X)\subset\Dqc(X).
\]
\end{sprop}

\begin{proof}  
In view of \eqref{HGam as lim.1}, one may assume $\Phi=\Phi_{\<\<Z}$ with $Z\subset X$ closed and the inclusion map
$i\colon (X\setminus Z)\hookrightarrow X$ quasi-compact. The assertion is given then by  
\cite[p.\kern1.4pt 25, (3.2.5)(iii)]{AJL}.
\end{proof}

\begin{sprop}\label{2RGams}
Let\/ $X$ be a locally noetherian scheme,  $E\in\Dqc(X),$ $\Phi$~an s.o.s.~in~$X$ and\/ $\I\set\I_\Phi$ \textup{(see \ref{ex:supports}\kf)}.  
Deriving the inclusion 
$\vG{\I}\hookrightarrow\vG{\Phi}$ from \textup{\ref{2Gammas}}
gives an \emph{isomorphism} $\R\vG{\I}E\iso\R\vG{\Phi}E$.
\end{sprop}

\begin{proof}
One needs the natural maps $H^n_{\<\I}\<E\iso H^n_{\<\Phi}E\ (n\in\mathbb Z)$ to be isomorphisms. 
By \eqref{IsubPhi},  \eqref{HGam as lim.1} and \eqref{HGam as lim.2}, and the fact that 
for $I\in\I$, \kf $\I_{\<\Phi_{\mkern-2.5mu Z\<(I\>)}}\!=\>\I_{\<\<I}$ (see~\eqref{radsup} with $I_1\set I$), one reduces to where 
$\Phi=\Phi_{\<Z(I\>)}$ for some quasi-coherent $\OX$-ideal~$I$, and $\I=\I_{\<\<I}$.
As $Z(I\>)$ is proregularly embedded in~$X$  (\cite[p.\,16, Example (a)]{AJL} and 
the lines before it), the assertion is given  by  \cite[p.\,25, (3.2.4)]{AJL}.\looseness=-1
\end{proof}

From \ref{Dqc to itself} and \ref{2RGams} (or from  {\cite[p.\kern1.4pt 21, (3.1.4)(iii)]{AJL}})
one gets:
\begin{sprop}\label{RGam respects Dqc}
For any locally noetherian scheme\/ $X$ and $\OX$-base $\I,$
\[
\mkern235mu\R\vG{\I}\>\Dqc(X)\subset\Dqc(X). \makebox{$\mkern218mu$}\square
\]
\end{sprop}

\vskip3pt

\centerline{* * * * *\va2}

Next, the stage is set for subsequent propositions.

\begin{slem}\label{Gam preserves injective}
{\hskip-1.4pt}Let\/ $X$ be a locally noetherian scheme, $\J$ an\/ $\OX$-base, \kern-.4pt$\Psi$ \kern-.7pt an s.o.s~in~$X$ and\/ $E$ a $(\<$degreewise\kf\) injective\/  $\OX$-complex.
Then both\/ $\vG{\J}E$ and\/~$\vG{\Psi}E$ are injective.
\end{slem}

\begin{proof} 
Using the results about injective $\OX$-modules on p.~127 of \cite{H}, and the fact that $\vG{\J}$ commutes with direct sums (see \ref{Gam and lim}), one reduces to
checking that if~$x\in X$ specializes to $x'\in X$ and $J(x,x')$ is the direct image on~$X$ of the constant sheaf on the closure~$\overline{x'}$ of $x'$ whose stalk at $x'$ is the injective hull $J_x$ of the residue field of $\cO_{\<\<X\<\<, \>x}\>$, then 
\[
\vG{\J}J(x,x') 
=
\begin{cases}
  J(x,x') &\textup{ if}\mkern8mu \overline x \in \Psi_{\<\<\J}, \\
  0&\textup{ otherwise,} 
 \end{cases}
 \mkern25mu\textup{and}\mkern25mu
\vG{\Psi}J(x,x') 
=
\begin{cases}
  J(x,x') &\textup{ if}\mkern8mu \overline{x'} \in \Psi, \\
  0&\textup{ otherwise.} 
 \end{cases}\\
\]
This checking is left to the reader, with the reminders that for any specialization $x''$ of $x'\<$, the $\cO_{\<\<X\<\<,\>x''}$-module structure on the stalk $J(x,x')_{x''}$ is induced by the natural homomorphism $\cO_{\<\<X\<\<,\>x''}\to\cO_{\<\<X\<\<,\>x}$, and that every element of $J_x$ is annihilated by a power of the maximal ideal of~$\cO_{\<\<X\<\<,\>x}\>$.
\end{proof} 
\vskip-3pt

Let $(X,\OX)$ be a ringed space, $\mathbf E$ an additive category and $\phi\colon \D(X)\to\mathbf E$ an additive functor. An $\OX$-complex~$F$ is \emph{(right-)$\phi$-acyclic} if the natural map 
$\phi F\to\R\phi F$ is a $\D(X)$-isomorphism.  (See \cite[p.\,50, Proposition 2.2.6]{Dercat}).

\begin{slem}\label{injective}
Let\/ $X$ be a scheme, $E^{\bullet}$ an\/ $\OX$-complex, $\Phi$ a finitary s.o.s.\ in\/~$X\<,$
and\/ $\I$ an $\OX$-base.  \va1

\textup{(i)} If\/ $E\in\Dqc(X)$\va{.6} and every $E^i$ is\/ $\vG{\Phi}$-acyclic  then the complex\/ $E^\bullet$ is\/ $\vG{\Phi}$-acyclic.\va1

\textup{(ii)} If\/ $X$ is locally noetherian and every $E^i$ is\/ $\vG{\I}$-acyclic, then \mbox{$E^\bullet$ is\/ $\vG{\I}$-acyclic.}
\end{slem}

\begin{proof}  
Using Remark~\ref{finitary restricts} and the fact that K-injectivity is preserved under restriction to open subsets---whence $\R\vG{\Phi}$ ``commutes" with such restriction---one finds that the assertions are local on~$X\<,$ so that $X$ may be assumed affine. 

In view of \eqref{Global Gam as lim} and since $\Phi$ is finitary, one can also assume that 
$\Phi=\Phi_{\!Z\<(I\>)}$ with $I$ generated by a finite sequence $\bt=(t_1,\dots,t_d)$ of global sections.
With~$\KK(t_i)$ the complex which vanishes everywhere except in degrees 0 and 1, where it is
\[
\cO_X^{(0)}\xto{\!\textup{natural}}\,\dirlm{}\big(\cO_X^{(0)}\xto{t_i} \cO_X^{(1)}\xto{t_i}\cO_X^{(2)}\xto{t_i}\cdots  \big)
\qquad\big(\cO_X^{(n)}\set\OX\ \forall\, n\ge0\big),\\[-3pt]
\]
and with $\KK(\bt)$ the bounded flat complex
$\KK(\bt)\set\otimes_{i=1}^d\KK(t_i)$, \cite[(3.2.3)]{AJL} gives an isomorphism
$\KK(\bt)\otimes E\iso \R\vG{\Phi}E$, which implies that the functor~$\vG{\Phi}$ 
is such that \cite[p.\,77, (a)]{Dercat} (dualized) applies, giving (i).\va1

As for (ii),  it's enough  that the maps induced by a K-injective resolution $E\to L$
 be isomorphisms $H^n\<\<\vG{\I} E\iso H^n\<\<\vG{\I} L\ (n\in\mathbb Z)$. For this,
\eqref{Gam as lim1} allows one to replace $\I$ by a quasi-coherent $\OX$-ideal~$I\<$ 
generated by a sequence $\bt=(t_1,\dots,t_d)$ of global sections. 
One has then an isomorphism $\R\vG{I}E\iso\KK(\bt)\otimes E$, 
(see proof of \mbox{\kf\cite[$(3.1.1)(2)'\Rightarrow(3.1.1)(2)$]{AJL}}), 
so \cite[p.\,77, (a)]{Dercat} (dualized) applies.
\end{proof}

\begin{subcosa}\label{K-flabby}
Again, $(X,\OX)$ is a ringed space. 

An $\OX$-complex~$E\in\Ab(X)$ is \emph{flabby} (or \emph{flasque}) if for~ every open $U\subset X$ the restriction map 
$\Gamma(X,E)\to \Gamma(U,E)$ is surjective; and \emph{quasi-flabby}  if the same holds for every \emph{quasi-compact}
open $U\subset X$.  
An \kf $\OX$-\kf complex $E$ is \emph{K-flabby}  (or \emph{K-flasque})
if for every s.o.s.~$\<\Phi$ \kf in~$X$ and for every open
$U\subset X\<,$ the natural  $\D(\Hr^0(U,\OU))$-map
$\Gamma^{}_{\<\!\Phi}(U,E)\to\R\Gamma^{}_{\<\!\Phi}(U,E)$ 
is an \emph{isomorphism}---see \cite[p.\,144, 5.19]{Sp}.
In other words, $E$ is K-flabby $\Leftrightarrow E$ is $\phi$\kf-acyclic for all functors~$\phi$ of the form 
$\Gamma^{}_{\<\!\Phi}(U,-)$. 

For example, if $E$  is K-injective then for any $\OX$-\kf complex~$C$, 
the \mbox{$\OX$-\kf complex}
$\sHom_\OX\<(C,E)$ is K-flabby \cite[p.\,142, 5.14 and p.\,141, 5.12]{Sp}.\vs1

If $E\to E'$ is a $\K(X)$-isomorphism then $E$ is K-flabby $\Leftrightarrow E'$ is K-flabby; 
and if $E\to E'$ is a quasi-isomorphism of K-flabby complexes then for all $\Phi$ and $U$ as above,
the induced map $\Gamma^{}_{\<\!\Phi}(U,E)\to\Gamma^{}_{\<\!\Phi}(U,E')$ 
is a quasi-isomorphism.\va1

If $E\to I$ is a K-injective resolution, then $E$ is K-flabby if and only if for all~$\Phi$, all~$U$ and all ~$n\in\mathbb Z$,  the induced map $\Hr^n\Gamma^{}_{\<\!\Phi}(U,E)\to \Hr^n\Gamma^{}_{\<\!\Phi}(U,I\>)$ is an~isomorphism. For any $\OX$-complex $C$ and $x\in X$, the stalk $(H^{n}\<\<\vG{\Phi}C\>)_x$ satisfies
\[
(H^{n}\<\<\vG{\Phi}C\>)_x=\Hr^{n}\>\dirlm{x\in U^{\mathstrut}}\<\Gamma^{}_{\<\!\Phi}(U,C\>)
=\dirlm{x\in U^{\mathstrut}}\<\Hr^n\Gamma^{}_{\<\!\Phi}(U,C\>).
\]
\vskip2pt\noindent
Hence if $E$ is K-flabby then the induced map $\vG{\Phi}E\to\vG{\Phi}I$ is a $\D(X)$-isomorphism, so that  K-flabby $\Rightarrow\vG{\Phi}$\kf-acyclic. 

If $f\colon W\to X$ is a map of ringed spaces, then any K-flabby
$\cO_W$-complex is $\fst$-acyclic \cite[p.\,147, 6.7(a) and p.\,141, 5.12]{Sp};
and the functor~ $\fst$ preserves K-flabbiness \cite[p.\,143, 5.15(b)]{Sp}. 
Upon replacing $E$ by a K-injective resolution,
it~follows then from \eqref{Gam and f} that  there is a natural functorial isomorphism
\begin{equation*}\label{RGamf}
\R\Gamma^{}_{\<\!\Phi_{\!f}}\<(W,E)\iso \R\Gamma^{}_{\<\!\Phi}(X,\R\fst E)
\qquad(E\in\D(W)).
\tag{\ref{K-flabby}.1}
\end{equation*}

Also, taking $f$ to be the natural map $(X,\OX)\to (X,\>\mathbb Z_X\<)$, one gets that 
any K-flabby $\OX$-complex is K-flabby as a complex of abelian sheaves.  
Hence for any integer~$n$ and $E\in\D(X)$, $\Hr^n_{\Phi}(E\>)$~depends, 
as an abelian group,  only on $X$ (not on~$\OX$)---and likewise for open $U\subset X$, whence  for the abelian sheaves $H^n_{\<\Phi}(E\>)$. \va2

\end{subcosa}

\begin{subcosa}\label{flabby}

An $\OX$-module $E$---as a complex $E^{\bullet}$ vanishing in all nonzero degrees---is~flabby if 
$\>\Hr^1_{\<Z}\mkern-1.5mu(X,E)=0$ for all closed $Z\subset X$ (see \cite[I, Corollaire 2.12]{SGA2}),
and only~if $\Hr^n_{\Phi}\<(X,E)=0$ for every s.o.s.\kern3pt$\Phi$ and $n>0$ 
(see \cite[p.\,174, 4.4.3(a)]{Go}). In~particular, any injective $\OX$-module is flabby.

The restriction of a flabby $\OX$-module~$E$ to any open 
$U\subset X$ is (clearly) a flabby $\mathcal O_U$-module; 
it follows that \emph{a~flabby $\OX$-module~$E\<$ is K-flabby.}

Conversely, \emph{if\/  $E$ is K-flabby} 
then~$\>\Hr^1_{\<Z}(X,E)\cong \Hr^1\Gamma^{}_{\<\<\!Z}(X,E^{\bullet})=0$,
and therefore $E$ \emph{is~flabby}. (Alternatively, see \cite[5.13(a)]{Sp}.)

Actually, \emph{any bounded-below quasi-flabby\/ $\OX$-complex
is K-flabby}. To prove this, use the~dual version of \cite[Proposition 2.7.2]{Dercat}, 
whose hypotheses hold for the class of 
flabby~$\OX$-modules by virtue~of the  second paragraph on page 147  and 
Th\'eor\`eme 3.1.2 + Corollaire on page 148 in~\cite{Go}.\va2 (Alternatively, see \cite[2.2(c) and~5.15(c)]{Sp}.)

Likewise, \emph{if\/ $X$ is quasi-noetherian then for every s.o.s.~$\Phi$ in\/~$X$ and 
quasi-compact open\/ $U\subset X,$ any bounded-below quasi-flabby\/ $\OX\<$-complex
is\/ $\Gamma^{}_{\<\!\Phi}(U,-)$-acyclic and\/ $\vG{\Phi}\<$-acyclic.} (Use
\cite[Proposition 4]{Kf} instead of \cite[Th\'eor\`eme 3.1.2]{Go}.)

\end{subcosa}

\goodbreak
\begin{slem}\label{Gam preserves flabby} 
 Let\/ $\<X\<$ be a topological space, $\Psi$  an s.o.s.\kern2pt in~$X,$\va1 $E$ an\/ \mbox{$\Ab(X)$-complex}.\va1 

\textup{(i)} Suppose that\/ $\Psi=\Phi_Y$ for some\/ $Y\subset X\<,$ or that every\/~$Z\in\Psi$ is noetherian. 
If\/~$E$ is flabby\va1 then so is\/ $\vG{\Psi}E\<$.\looseness=-1

\textup{(ii)} Suppose\/ $X$ quasi-compact and\/ $\Psi$ finitary. If\/ $E$ is quasi-flabby then so
is $\vG{\Psi}E$.\va1

\end{slem} 
  
\begin{proof} 
(i).~Let $U\subset X$ be open. By \eqref{intersect sos}, any $s\in\Gamma(U,\vG{\Psi}E)$ vanishes on~$U\setminus \<Z$ 
for some $Z\in\Psi$, hence extends to an 
$s'\in\Gamma(U\>\cup\>(X\<\setminus\<\<Z),E)$ that~vanishes on~
$X\setminus Z\>$. Since $E$ is flabby, therefore
$s'$ extends to an $s''\in\Gamma(X,E)$. This
$s''$ is an extension of~$s$ to~$\Gamma(X,\vG{\Psi}E)$.\va1

(ii).~Let $U\subset X$ be open and quasi-compact. By \eqref{intersect sos}, any $s\in\Gamma(U,\vG{\Psi}E)$ vanishes on~$U\setminus \<Z$ for some $Z\in\Psi$, and since $\Psi$ is finitary and $X$ quasi-compact, one may assume that $X\setminus Z$ quasi-compact. The section $s$
extends to a section 
$s'\in\Gamma(U\>\cup\>(X\<\setminus\<\<Z),E)$ that~vanishes on
$X\setminus Z\>$. Since $E$ is quasi-flabby and $U\>\cup\>(X\<\setminus\<\<Z)$ is quasi-compact,
therefore $s'$~extends to an $s''\in\Gamma(X,E)$. 
This~$s''$ is an extension of~$s$ to~$\Gamma(X,\vG{\Psi}E)$.\va1
\end{proof}

With $\phi\set\:$empty set, the (\emph{\kern-.4pt Krull}\kf) \emph{dimension} dim$\,X$ of a topological space $X\ne\phi$ is the supremum ($\le\infty$) of the set of those integers $n$ such that there
exists a strictly increasing sequence  $\phi\ne Z_0<Z_1<\cdots<Z_n$ of irreducible closed subsets of 
$X\<;$ and dim.$\,\phi\set-1$. 

\begin{slem}\label{L:flabby-fd} 
If\/ $X$ is a  finite-dimensional noetherian topological space, then any
flabby\/ $\Ab(X)$-complex is K-flabby.
\end{slem}

\begin{proof} Any open $U\subset X$ is noetherian; and
$\textup{dim}\,U\le \textup{dim}\,X$, since the $X$-closure~$\,\overline{\! Z}$ of 
an irreducible $Z$ closed in~$U$ is irreducible and such that 
$\,\overline{\! Z}\cap U=Z$.

For any s.o.s.\kern2pt $\Phi$ in $X\<,$ it holds then that 
\[
\Hr^{\>p}_\Phi\<(U,F)=0 \textup{ for all $F\in\Ab(U)$ and integers $p> \textup{dim}\,X,$}
\]
see \cite[Tag 02UZ]{Stacks}, whose proof works  with ``H" replaced by ``H$_\Phi.\<$"%
(Use \eqref{RGamf}, and \ref{RGam and colim1} below. Note too
that if $X$ is irreducible then any constant sheaf in~$\Ab(X)$ is flabby; 
moreover, if dim$\,X=0$ then
the only nonempty open subset of $X$ is $X$~itself, whence every $E\in\Ab(X)$ is flabby.)

Since every abelian sheaf embeds\va{.5} into  a flabby one \cite[p.\,147, 2nd paragraph]{Go}, it results as in 
the proof of (ii)$\Rightarrow$(iii)$\Rightarrow$(a) in  \cite[pp.\,76--77, (2.7.5)]{Dercat} (dualized) that 
any flabby $\Ab(X)$-complex  is $\Gamma^{}_{\<\!\Phi}(U,-)\>$-acyclic, thus K-flabby.\va4
\end{proof}

\centerline{* * * * *\va{-1}}

\begin{sprop}\label{derived intersect} 
Let\/ $X$ be a ringed space, $E$ an\/ $\OX$-complex, and each of\/ $\Phi$~ 
and\/ $\Psi$ an s.o.s.~in\/ $X\<.$ Suppose one of the following holds.\va1

\textup{(i)} $E\in\Dpl(X)$ and\/ $\Psi$ is as in \textup{\ref{Gam preserves flabby}(i).} 

\textup{(ii)} $X$ is quasi-noetherian, $E\in\Dpl(X),$ and\/ $\Psi$ is finitary.

\textup{(iii)} $X$ is noetherian and finite-dimensional.\va1

Then the natural map \(\kern-1pt arising from \eqref{intersect vG} is an isomorphism\va{-2}
\[
\mkern14mu\gamma^{}_{\Phi\<,\Psi}\colon\R\vG{\Phi\cap\>\Psi}E\iso
\R\vG{\Phi}\R\vG{\Psi}E.
\]

If, moreover, $E$~is cohomologically bounded-below
then the natural map \(\kern-1pt arising from \eqref{intersect sos}\) is an isomorphism\va{-2}
\[
\mkern14mu\overline\gamma^{}_{\Phi\<,\Psi}\colon\R \Gamma^{}_{\<\!\Phi\cap\>\Psi}(X,E)\iso
\R \Gamma^{}_{\<\!\Phi}(X,\R\vG{\Psi}E).
\]
\end{sprop}

\begin{proof} One can assume $E$ to be injective; and since
$\gamma^{}_{\Phi\<,\Psi}$ is an isomorphism if it is so locally, therefore, if $E\in\Dpl(X)$
then one can also assume $E\>$ bounded-below.
As~in~\ref{flabby}, bounded-below plus flabby implies K-flabby; so if (i) holds then
by~\ref{Gam preserves flabby}(i),
$\vG{\Psi}E$ is K-flabby, hence $\Gamma^{}_{\<\!\Phi}(X\<,-)$\kf- and $\vG{\Phi}$\kf-acyclic; 
if (ii) holds, argue similarly, replacing \ref{Gam preserves flabby}(i) with \ref{Gam preserves flabby}(ii);
and if (iii) holds, reach the same conclusion via \ref{L:flabby-fd}.
That~$\gamma^{}_{\Phi\<,\Psi}$ and $\overline\gamma^{}_{\Phi\<,\Psi}$ are isomorphisms follows.
\end{proof}
 
\begin{sprop}\label{derived intersect2}
Let\/ $X$ be a scheme, $E$ an\/ $\OX$-complex, and each of\/ $\Phi$ and\/~$\Psi$ 
an s.o.s.~in~$X\<$. \va1 

\textup{(i)} If\/ $E\in\Dqc(X)$ and both\/ $\Phi$ and $\Psi$ are finitary, then
the natural map \(\kern-1pt arising from \eqref{intersect vG}\) is an \emph{isomorphism}\va{-1}
\[
\gamma^{}_{\Phi\<,\Psi}\colon\R\vG{\Phi\cap\>\Psi}E\iso
\R\vG{\Phi}\R\vG{\Psi}E.
\]

\textup{(ii)} If\/ $X$ is locally noetherian and\/ $E$ is cohomologically bounded-below, then the natural map 
\(\kern1pt from \eqref{intersect sos}\) is an isomorphism
\[
\overline\gamma^{}_{\Phi\<,\Psi}\colon\R \Gamma^{}_{\<\!\Phi\cap\>\Psi}(X,E)\iso
\R \Gamma^{}_{\<\!\Phi}(X,\R\vG{\Psi}E).
\]
\end{sprop}

\begin{proof} 
The complex $E$\va{.5} can be assumed  K-injective, and furthermore, bounded-below if $E$ is cohomologically so. \va1

(i).~One checks that\va{.5} for $n\in\mathbb Z,$ the cohomology map 
$H^{n}\>\gamma^{}_{\Phi\<,\Psi}$ factors as
the following sequence of isomorphisms,\va{.5}
in which $V,\,W$ are such that $X\setminus V$ and~$X\setminus W$ are retrocompact in~$X\<$:\va{-3}
\begin{align*}
H^n\<\vG{\Phi\cap\>\Psi}E
&\iso
\dirlm{\lift.15,{\!\!\halfsize{$V$}\in{\lift1,\Phi,}},\!\!}\,
\dirlm{\lift.15,{\!\!\halfsize{$W$}\in{\lift1,\Psi,}},\!\!}\!  H^n\<\vG{V\cap W}E
&&\eqref{Global Gam as lim}\\[4pt]
&\iso
\dirlm{\lift.15,{\!\!\halfsize{$V$}\in{\lift1,\Phi,}},\!\!}\,
\dirlm{\lift.15,{\!\!\halfsize{$W$}\in{\lift1,\Psi,}},\!\!}\! H^n\R\vG{V}\vG{W}E
&&\textup{\cite[p.\,25, 3.2.5(ii)]{AJL}}\\[4pt]
&\iso
\dirlm{\lift.15,{\!\!\halfsize{$V$}\in{\lift1,\Phi,}},\!\!}\! 
 H^n\R\vG{V}\>\>\dirlm{\lift.15,{\halfsize{$W$}\in{\lift1,\Psi,}\!},\,}\!\vG{W}E
&&\textup{\ref{Dqc to itself}},\ \textup{\ref{RGam and colim2} (below)}\\[4pt]
&\iso
\dirlm{\lift.15,{\!\!\halfsize{$V$}\in{\lift1,\Phi,}},\!\!}\! 
 H^n\R\vG{V}\vG{\Psi}E\iso H^n\R\vG{\Phi}\vG{\Psi}E
 &&\eqref{Global Gam as lim},\ \eqref{HGam as lim.1}.
\end{align*}

\vskip4pt
\noindent Hence $\gamma^{}_{\Phi\<\<,\Psi}$ is an isomorphism.\va2

(ii).~By \ref{Gam preserves injective} the bounded-below $\OX$-complex~$\vG{\Psi}E$ is injective, hence K-injective, hence $\Gamma^{}_{\<\!\Phi}(X,-)$-acyclic, 
so that $\overline\gamma^{}_{\Phi\<,\Psi}$ is indeed an isomorphism.
\end{proof} 

\begin{sprop}\label{derived intersect4}
Let\/ $X$ be a locally noetherian scheme, $E$ an\/ $\OX$-complex and each of\/ $\I$ 
and\/ $\J$ an\/ $\OX\<$-base. The natural map \(\kern 1pt from \eqref{local intersect sos}\) is an isomorphism
\[
\gamma^{}_{\I\<,\J}\colon\R\vG{\I\cap\>\J}E\iso
\R\vG{\I}\R\vG{\J}E\>.
\]
If\/ $X$ is noetherian and finite-dimensional, then 
the natural map \(\kern 1pt from \eqref{intersect bases}\) is an isomorphism\va{-3}
\[
\overline\gamma^{}_{\I\<,\J}\colon\R\Gamma^{}_{\<\!\I\cap\>\J}(X\<,E)\iso
\R\Gamma^{{}}_{\<\!\I}(X\<, \R\vG{\J}E).
\]
\end{sprop}
\noindent
\emph{Proof.}
One can assume $E$ to be K-injective and injective. 

By~\ref{Gam preserves injective} and~\ref{injective}(ii), $\vG{\J}E$ is $\vG{\I}$\kf-acyclic, 
and so $\gamma^{}_{\I\<,\J}$ is an isomorphism.

Next, $\vG{\I}E$ is injective (\ref{Gam preserves injective}), hence flabby, so one has natural isomorphisms
\[
\R\Gamma^{}_{\<\!\I}(X\<,E)\iso \Gamma^{}_{\<\!\I}(X\<,E) 
\underset{\eqref{intersect bases}}\iso \Gamma(X\<,\vG{\I}E) 
\underset{\ref{L:flabby-fd}}\iso \R\Gamma(X\<,\vG{\I}E) 
\iso \R\Gamma(X\<,\R\vG{\I}E),\\[-1pt]
\]
via which, one checks, $\overline\gamma^{}_{\I\<,\J}$ factors as the sequence of natural isomorphisms\va{-1}
\[
\R\Gamma^{}_{\<\!\I\cap\J}(X\<,E) 
\!\iso\! \R\Gamma(X\<,\R\vG{\I\cap\J}E)
\!\iso\! \R\Gamma(X\<,\R\vG{\I}\R\vG{\J}E)
\!\iso\!\R\Gamma^{}_{\<\!\I}(X\<,\R\Gamma^{}_{\<\!\J}E).\va2\quad \square
\]

\centerline{* * * * *\va{-3}}

\begin{subcosa}\label{flabby res}
For a topological space~$X\<$, bounded-below complexes $E\in\Ab(X)$ have canonical (Godement) flabby 
resolutions $E\to G(E\>)$, with $G(E)$ bounded below and varying functorially with $E$ (see \cite[proof of 3.9.3.1]{Dercat}). 

If $X$ is quasi-noetherian, then the functor~$G$, with ``flabby" replaced by ``quasi-flabby," extends to unbounded $E\>$:
with $E^{{\sst\ge} -n}$ the complex obtained from $E$ by replacing $E^m$
with 0 for all  $m<-n$, and with $E^{{\sst\ge} -n}\to E^{{\sst\ge} -(n+1)}$ the obvious map, one has\va{-2}
$E=\dirlm{} \<E^{{\sst\ge} -n}\>$;
and since a filtered direct limit of flabby (hence quasi-flabby) sheaves is 
quasi-flabby \cite[p.\,641, Corollary 7\kf]{Kf}, one can set 
$
\smash{G(E\>)\set \dirlm{n\in\mathbb Z^{\mathstrut}} G(E^{{\sst\ge} -n}).}
$ 
\end{subcosa}

\vskip6pt
\begin{sprop} \label{RGam and colim1}
Let $\<X\!$ be a\va1 quasi-noetherian topological space, $\Phi$~a finitary s.o.s.~ in~$X\<,$ $A$ a small filtered category, and\/ $\mathscr M$\va2  a functor from~$A$ to the category 
of\/ $\Ab(X)$-complexes. 
If\/ $\dirlm{\lift.5,\sst A,} \mathscr M$ is bounded-below,\va{2} 
or if\/ $X$~is noetherian and of finite dimension,
then for every\/ $n\in\mathbb Z,$\va{1.5} the natural maps are isomorphisms\va{-1}
\begin{align*}
\dirlm{A} \!(H^n_{\<\Phi}\smcirc\mathscr M\>)&\iso 
H^n_{\<\Phi}\>\>\dirlm{A}\<\< \mathscr M,\\[3pt]
\dirlm{A} \!(\Hr^n_\Phi(X\<,-)\smcirc\mathscr M\>)&\iso 
\Hr^n_\Phi(\<X\<,\dirlm{A}\<\< \mathscr M\>).
\end{align*}
\vskip3pt
\noindent In particular, $\R\vG{\Phi}$ and\/ $\R\Gamma^{}_{\<\!\Phi}(X,-)$ commute with  small direct sums.
\end{sprop}

\begin{proof}  
With $G$ as in \ref{flabby res}, 
$\>\dirlm{A}\<\mathscr M\to\dirlm{A}\<(G\smcirc\mathscr M\>)$ is a quasi-isomorphism\va3 whose target is, by 
\cite[Corollary 7\kf]{Kf}, a flabby,\va2 hence as in \ref{flabby} or by~\ref{L:flabby-fd}, K-flabby, hence $\vG{\Phi}$-acyclic, complex. \va1

The first isomorphism is then the natural composite isomorphism\va{-1}
\begin{align*}
\dirlm{A}\!( H^n_{\<\Phi}\smcirc\mathscr M\>)
 \iso \dirlm{A}\! (H^n\mkern-1.5mu\smcirc\vG{\Phi}\smcirc G\smcirc\mathscr M\>)
&\iso H^n\>\dirlm{A}\<\< (\vG{\Phi}\smcirc G\smcirc\mathscr M\>)\\[3pt]
&\underset{\ref{Gam and lim}}{\iso}
H^n\<\vG{\Phi}\>\dirlm{A}\!(G\smcirc\mathscr M\>)
 \iso 
H^n_{\<\Phi}\>\>\dirlm{A}\<\< \mathscr M.
\end{align*}
\vskip1pt
\noindent The second is obtained similarly, via \ref{L:flabby-fd} and \ref{lims of Gams}.

As for direct sums, the standard argument associates to any set $I$ the ordered (by~inclusion) set $A$ of finite subsets of $I$, 
regards $A$ in the usual way as a filtered category, and uses commutativity of the additive functor 
$H^n_{\<\Phi}$ with finite direct sums to get, for any~family $(M_i)_{i\in I}$ of $\Ab(X)$-complexes, ~any $n\in\mathbb Z$, and~\mbox{$M_\alpha\set\oplus_{i\in\alpha}\>M_i\ (\alpha\in A)$,} natural isomorphisms:\va{-1}
\begin{align*}
H^n\big(\!\oplus_{i\in I}\R\vG{\Phi}M_i\big)
&\iso\oplus_{i\in I}H^n_{\<\Phi}M_i
\iso\dirlm{\alpha\in \<A^{\mathstrut}}H^n_{\<\Phi} M_\alpha\\[4pt]
&\iso H^n_{\<\Phi}(\>\dirlm{\alpha\in \<A^{\mathstrut}}M_\alpha)\iso 
H^n_{\<\Phi}\big(\!\oplus_{i\in I}M_i\big)
  =H^n\R\vG{\Phi}\big(\!\oplus_{i\in I}M_i\big).
\end{align*}
\vskip3pt
\noindent Thus the natural map is an isomorphism \va{-1}
\[
\oplus_{i\in I}\>\>\R\vG{\Phi}M_i\iso \R\vG{\Phi}(\oplus_{i\in I}M_i).
\]

Similar considerations hold with $\Gamma^{}_{\!\<\Phi}(X,-)$ in place of $\vG{\Phi}\>$.
\end{proof}

\begin{sprop} \label{RGam and colim2}
Let $X$ be a  scheme,  $\Phi$~a finitary\/  s.o.s.\kern2.4pt in\/~$X\<,$ $A$ a small filtered category, 
$\mathscr M$  a functor from~$A$ to the category of\/ $\OX$-complexes with
quasi-coherent homology,  and\/ $n\in\mathbb Z$.
The natural map is an isomorphism
\begin{equation*}
\dirlm{A} \!(H^n_{\<\Phi}\smcirc\mathscr M\>)\iso
H^n_{\<\Phi}\>\>\dirlm{A}\<\< \mathscr M.
\end{equation*}
In particular, $\R\vG{\Phi}$ commutes with small direct sums in\/ $\Dqc(X)$.
\end{sprop}

\begin{proof}
Using Remark~\ref{finitary restricts} and the fact that K-injectivity is preserved under restriction to open subsets---whence $\R\vG{\Phi}$ ``commutes" with such restriction---one finds that the first assertion is local on~$X\<,$ so that $X$ may be assumed affine.

From \eqref{HGam as lim.1} it follows that it's enough to treat the case $\Phi=\Phi_{\<\<Z}$, with $Z\subset X$ closed and such that $X\setminus Z$ is retrocompact in~$X$; therefore it may be assumed that 
$Z=\Supp(\OX/{\mathbf t}\OX\<)$ with~$\mathbf t$ a finite sequence in~$\Gamma(X,\OX\<)$.
Then the assertion is a simple consequence of the fact that for complexes with quasi-coherent homology, 
applying $\R\vG{\<Z}$ is the same as tensoring with 
the complex~$\KK({\mathbf t})$ (see~proof of \ref{injective}).

The argument for direct sums is as in the proof of \ref{RGam and colim1}.
\end{proof}

\begin{sprop}\label{RGam and colim3}
Let $X$ be a  locally noetherian scheme,  $\I$~an\/  $\OX\<$-base, $A$ a small filtered category,  
$n\in\mathbb Z,$ and\/
$\mathscr M$  a functor from~$A$ to the category of\/ $\OX$-complexes.\va1
The natural map is an isomorphism\va{-2}
\[
\dirlm{A} \!(H^n_{\<\I}\<\<\smcirc\mathscr M\>)\iso 
H^n_{\<\I}\>\dirlm{A} \mathscr M.
\]
In particular, $\R\vG{\I}$ commutes with small direct sums in\/ $\D(X)$.\va2
\end{sprop}

\begin{proof}
Imitate the proof of \ref{RGam and colim2}, 
replacing \cite[(3.2.3)]{AJL} in the proof of~\ref{injective} by \cite[(3.1.1)(2)]{AJL}.
(Alternatively, using \ref{2RGams} deduce the result from \ref{RGam and colim2}.)
\end{proof}

\noindent\emph{Remark.} More generally, \ref{RGam and colim1}--\ref{RGam and colim3} hold 
when $A$ is a \emph{pseudo-filtered} category \cite[p.\,216, Exercise 2]{Ma2}. 

\subsection{Coreflections}\label{tensorcompat}
\stepcounter{thm}
This section expands on \emph{coreflectiveness,} both abstractly and in the context of ringed spaces. In the following section there is a discussion of  \mbox{\emph{$\otimes$-compatible coreflectiveness}} in the context~of symmetric monoidal categories, leading\- to the subsequently important notion of \emph{idempotent pairs} in such categories.

\begin{sdef}\label{coreflection} Let $\D$ be a category, with identity functor~$\mathbf 1_{\D}$,
let $\Gamma\colon\D\to\D$ be a functor and\/ $\iota\colon\Gamma\to\mathbf 1_{\D}$ a functorial map. 
The pair $(\Gamma\<,\iota)$ is  a \emph{coreflection} \emph{of}~$\D$ 
(or \emph{coreflecting in}~$\D$, or \emph{colocalizing in}~$\D$)
if for all $E\in\D$,  the functorial maps\/
$\Gamma\big(\iota(E)\big)$ and\/ $\iota(\Gamma E)$ are equal
isomorphisms from $\Gamma\Gamma E$ to $\Gamma E$.

 The functor $\Gamma$ is a \emph{coreflector} if there exists an
 $\iota$ such that $(\Gamma,\iota)$ is a coreflection.\vs1
\end{sdef}

 \begin{slem}[well-known] \label{coconds} 
The pair\/ $(\Gamma\<,\iota)$ is coreflecting in\/ $\D \!\iff\!$ for~all\/ \mbox{$F,G\in\D$} the map induced by\/ $\iota(G)$ is an \emph{isomorphism}
\begin{equation*}\label{gamadj1}
\Hom_\D(\Gamma\<F\<,\>\Gamma G)\iso \Hom_\D(\Gamma\<F\<,\> G).
\tag{\ref{coconds}.1}
\end{equation*}
\end{slem}

\begin{proof}
For $\Rightarrow,$ one checks,
using $\Gamma\big(\iota(G)\big)=\iota(\Gamma G)$ and the functoriality of~$\iota$, that the 
natural composite map
$$
\Hom_\D\<(\Gamma\< F, G)\lra
\Hom_\D\<(\Gamma\Gamma\< F, \Gamma G)\!\iso\!
\Hom_\D\<(\Gamma\< F, \Gamma G)
$$  
is inverse to the map in \eqref{gamadj1}. 

For~$\Leftarrow,$  simple considerations applied to~ \eqref{gamadj1} with $\Gamma G$ in place of~$G$ show that $\iota(\Gamma G)$ is an isomorphism; and functoriality of~$\iota$ implies\va{-2}
\[
\iota(G)\smcirc\iota(\Gamma G)=\iota(G)\smcirc \Gamma(\iota(G)),\\[2pt]
\]
and so \eqref{gamadj1} with $F=\Gamma G$ gives that $\iota(\Gamma G)=\Gamma(\iota(G))$.
\end{proof}

\begin{exams} \label{coreflect exams}
(a) Let $(\Gamma,\iota)$ be coreflecting in $\D$, and let $\D'\subset \D$ be a full subcategory such that 
$\Gamma\D'\subset\D'$. Let $\Gamma'\colon\D'\to\D'$ be the restriction $\Gamma|_{\D'}$.
Then $\iota$~induces a functorial map $\Gamma'\to\b1_{\D'}$, and $(\Gamma'\<,\iota')$ is 
coreflecting in~$\D'\<$.

(b) Setting $\Psi=\Phi$ in \eqref{intersect vG}, one gets that for any s.o.s.~$\Phi$ on a ringed space $X\<,$ the functor $\vG{\Phi}$ and 
its inclusion into $\mathbf 1_{\<\A(X)}$ constitute a coreflection of $\A(X)$; and likewise, via \eqref{local intersect sos}, for any functor $\vG{\I}$ with $\I$ an $\OX$-base on a scheme~$X\<$.\va1

(c)  Let $X$ be a ringed space, $E$ an\/ $\OX$-complex, and $\Phi$ an s.o.s.~in $X\<.$ Propositions~\ref{derived intersect} and~\ref{derived intersect2} give that if one of the following conditions (i)--(v) holds, then there exists a natural isomorphism
\begin{equation*}\label{gamphiphi}
\gamma^{}_{\Phi\<,\Phi}\colon\R\vG{\Phi}E \iso \R\vG{\Phi}\R\vG{\Phi}E.
\tag{\ref{coreflect exams}.1}
\end{equation*}

(i) $E\in\Dpl\<(X)$, and $\Phi=\Phi_Y$ for some $Y\subset X$.

(ii) $E\in\Dpl\<(X)$,  and  every member of $\Phi$ is quasi-compact.

(iii) $X$ is quasi-noetherian, $E\in\Dpl\<(X)$, and $\Phi$ is finitary.

(iv) $X$ is noetherian and finite\kf-dimensional.

(v) $X$ is a scheme, $E\in\Dqc(X)$, and $\Phi$ is finitary.\va2

The next lemma implies that with 
$\iota^{}_\Phi\colon \R\vG{\Phi}\to {\mathbf 1}$ the natural map, 
both $\R\vG{\Phi}\iota^{}_\Phi$ and~ 
$\iota^{}_\Phi(\R\vG{\Phi})$ are  inverse to \eqref{gamphiphi},
so they are equal isomorphisms from $\R\vG{\Phi}\R\vG{\Phi}$ to $\R\vG{\Phi}$. 
Since $\R\vG{\Phi}\Dpl(X)\subset\Dpl(X)$ (locally verifiable, so one need only consider
bounded-below complexes\dots), and by \ref{Dqc to itself}, $\R\vG{\Phi}\Dqc(X)\subset\Dqc(X)$, therefore:
\va2

\emph{If}~(i), (ii) or\/~(iii) \emph{holds, then\/ $(\R\vG{\Phi},\iota^{}_\Phi)$ is coreflecting in~$\Dpl(X);$ 
if\/ \textup{(iv)} holds, then $(\R\vG{\Phi},\iota^{}_\Phi)$ is coreflecting in\/~$\D(X)$;
and if\/ \textup{(v)} holds, 
 $(\R\vG{\Phi},\iota^{}_\Phi)$ is coreflecting in~$\Dqc(X)$}. \va2

Similarly, using \ref{derived intersect4} one gets:\va2 

\emph{If\/ $X$ is a locally noetherian scheme 
and\/ $\I$ is an\/ $\OX\<$-base, then\/  $(\R\vG{\I},\iota_\I)$ is co\-reflecting in\/ $\D(X)$---and also, 
by\/ \textup{\ref{RGam respects Dqc},} in\/~$\Dqc(X)$.}

\begin{slem}\label{symm}
For systems of supports\/ $\Phi,\;\Psi$ in a topological space, and bases\/ $\I,\;\J$ over a scheme,  the subtriangles in the following natural diagrams commute.
\[
\def\1{$\R\vG{\Phi\cap\>\Psi}$}
\def\2{$\R\vG{\Psi}$}
\def\3{$\R\vG{\Phi}$}
\def\4{$\R\vG{\Phi}\R\vG{\Psi}$}
 \bpic[xscale=3.5,yscale=1.75]
  \node(11) at (1,-1){\1};
  \node(12) at (2,-1){\2};
  
  \node(21) at (1,-2){\3};
  \node(22) at (2,-2){\4};
 
   \draw[->] (11) -- (12) node[below=9, midway, scale=0.85] {\kern44pt\textup{\circled2}};
   \draw[<-] (21) -- (22) ;
   
   \draw[->] (11) -- (21) ;   
   \draw[<-] (12) -- (22) ;
   
    \draw[->] (11) -- (22) node[below=-1, midway, scale=0.85]
                          {\textup{\circled1} \kern17pt$\gamma^{}_{\Phi\<,\<\Psi}\mkern70mu$} ;
    
 \epic
 \qquad
 \def\1{$\R\vG{\I\cap\>\J}E$}
\def\2{$\R\vG{\J}E$}
\def\3{$\R\vG{\I}E$}
\def\4{$\R\vG{\I}\R\vG{\J}E$}
 \bpic[xscale=3.5,yscale=1.75]
  \node(11) at (1,-1){\1};
  \node(12) at (2,-1){\2};
  
  \node(21) at (1,-2){\3};
  \node(22) at (2,-2){\4};
 
   \draw[->] (11) -- (12) node[below=9, midway, scale=0.85] {\kern44pt\textup{\circled4}};
   \draw[<-] (21) -- (22) ;
   
   \draw[->] (11) -- (21) ;   
   \draw[<-] (12) -- (22) ;
   
    \draw[->] (11) -- (22) node[below=-1, midway, scale=0.85]
                          {\textup{\circled3}\kern17pt$\gamma^{}_{\I\<,\J}\mkern60mu$
                           } ;
 \epic
\]  
\end{slem}

\begin{proof}
For commutativity of \circled1 it's enough  (by the universal property of derived functors) to check after composing with the natural map $\vG{\Phi\cap\>\Psi}\to\R\vG{\Phi\cap\>\Psi}\>$, for which purpose it's enough to have commutativity of  the subdiagrams in the following natural expansion of \circled1, commutativities that result directly from definitions. \va{-4}

\[
\def\1{$\vG{\Phi\cap\>\Psi}$}
\def\2{$\ \vG{\Phi}$}
\def\3{$\vG{\Phi}\vG{\Psi}$}
\def\4{$\R\vG{\Phi\cap\>\Psi}$}
\def\5{$\R\vG{\Phi}$}
\def\6{$\R\vG{\Phi}\R\vG{\Psi}$}
\def\7{$\R\vG{\Phi}\vG{\Psi}$}
 \bpic[xscale=1.85,yscale=1.25]
  \node(11) at (2.1,-2.5){\1} ;
  \node(15) at (2.1,-3.5){\2} ;
  
  \node(22) at (2.5,-1){\4} ;
  \node(23) at (.5,-4.1){\5} ;
  
  \node(33) at (4.5,-4.1){\6} ;
  \node(43) at (2.9,-3.5){\7 } ;
 
  \node(53) at (2.9,-2.5){\3} ;
  
%
   \draw[->] (22) -- (23) ;
   \draw[->] (43) -- (33) ;
%
   \draw[->] (43) -- (33) ;
   \draw[->] (53) -- (43) ;
   \draw[->] (11) -- (15) ;   
   \draw[->] (11) -- (22) ;   
   \draw[->] (22) -- (33) ;
   \draw[->] (33) -- (23) ;
   \draw[double distance=2] (11) -- (53) ;
   \draw[->] (15) -- (23) ;
 \epic
\]  

 \vskip-1pt
 That \circled2, \circled3 and \circled4 commute is shown similarly. (Details left to the reader.)
\end{proof}

(d) Variants of the foregoing examples will emerge in the contexts of topological rings and of noetherian formal schemes (Propositions~\ref{Gamma' coreflects}, \ref{vGI tensor-coreflects} and \ref{Dqct to itself2}).
\end{exams}

\vskip7pt
\centerline{* * * * *}

\begin{subcosa}\label{essimage0}
The \emph{essential image} $\D_\Gamma$\va{.6} of a functor~$\Gamma\colon\D\to\D$ is 
the strictly full subcategory of~$\D$ spanned by the objects\va{1.5} 
~$\Gamma\< F\ (F\in\D)$\va{1.1}. The functor~$\Gamma$ factors as
$\D\xto{\lift1.15,\Gamma^0\<,}\D_\Gamma\overset{\lift.85,\>\>j,}\hookrightarrow\D$ where\va1 $j$ is the inclusion functor.

It is easy to see that if $(\Gamma\<,\iota)$ is coreflecting in\/ $\D$ then an object $E\in\D$ lies in~  
$\D_\Gamma$  if and only if $\iota^{}(E)$ is an isomorphism $\Gamma E\iso E$.
 \end{subcosa}

\begin{slem}\label{coref and adj} 
The pair $(\Gamma\<,\iota)$ is coreflecting in\/ $\D$ $\!\iff\!$ there is an adjunction 
$j\<\dashv\Gamma^0$ with counit $\iota$.  Thus\/ $\Gamma$ is a coreflector if and only if\/ 
$\Gamma^0$ is right-adjoint to~$j$ \(that is, if and only if\/ $\D_\Gamma$ is 
a \emph{coreflective subcategory} of\/ $\D$ 
 \cite[p.\,91, bottom]{Ma2}\).
\end{slem}

\begin{proof}
Let $E\in\D_\Gamma$, so that there is  a $\D$-isomorphism 
$\alpha\colon jE\iso\Gamma F\ (F\in\D)$. For~any $G\in\D$, the square in the following diagram clearly commutes. \va{-1}
\[
\def\1{$\Hom_\D(\Gamma F\<,\Gamma G)$}
\def\2{$\Hom_\D(jE\<,\Gamma G)$}
\def\3{$\Hom_\D(jE\<,\>j\Gamma^0 G)$}
\def\4{$\Hom_\D(\Gamma F\<, G)$}
\def\5{$\Hom_\D(jE\<, G)$}
\def\6{$\Hom_{\D_\Gamma}(E\<, \Gamma^0 G)$}
 \bpic[xscale=3.6,yscale=1.35]
  \node(11) at (1,-1){\1};
  \node(12) at (2,-1){\2};
  \node(13) at (3,-1){\3};
  
  \node(21) at (1,-2){\4};
  \node(22) at (2,-2){\5};
  \node(23) at (3,-2){\6};
  
   \draw[->] (11) -- (12) node[above, midway,scale=.75]{$\Iso$}
                                    node[below=.5, midway, scale=0.75] {$\via\alpha$} ;
   \draw[double distance=2] (12) -- (13) ;
   
    \draw[->] (21) -- (22) node[above, midway,scale=.75]{$\Iso$}
                                     node[below=.5, midway, scale=0.75] {$\via\alpha$} ;
   
   \draw[->] (11) -- (21)  node[left, midway,scale=.75]{via $\iota G$} 
                                     node[right=-1, midway]{$\lift .1,a,$} ; 
   \draw[->] (12) -- (22)  node[right, midway,scale=.75]{via $\iota G$}  
                                     node[left=-1, midway]{$\lift.5,b,$} ; 
   \draw[double distance=2] (13) -- (23) ; \epic
\]
\vskip-4pt\noindent
Hence $(\Gamma,\iota)$ is coreflecting $\Leftrightarrow\:a$ is an isomorphism (see \ref{coconds}) 
$\Leftrightarrow\:b$ is an~isomorphism $\Leftrightarrow\:b$ gives an adjunction $j\<\dashv\Gamma^0$ whose counit\va{.5} (the image under $b$ of 
the identity~map of $j\Gamma^0G=\Gamma G$) is $\iota(G)$.
\end{proof}

\begin{subcosa}\label{essimage}
To illustrate, let  $(X\<,\OX)$ be a ringed space, let $\Phi$ be an s.o.s.~in~$X\<$, and 
let $\A_\Phi(X)\subset\A(X)$ be the full subcategory spanned by the 
\emph{$\Phi$-torsion} $\OX$-modules, that~is, those~$M$ such that $\vG{\Phi}M=M$. 
Then $\A_\Phi(X)$ is the essential image of~$\vG{\Phi}$, since for any $\OX$-isomorphism $M\iso \vG{\Phi} N$, \eqref{intersect vG} shows that $M$ is $\Phi$-torsion. 

If $X$ is a scheme, then in the preceding paragraph one can replace $``\A\,$'' by $``{\Aqc}$'' (see \ref{Cor2Gammas}); and if $\I$ is an $\OX$-base, one can replace ``$\>\Phi$''  by ``$\>\>\I$.''
If $\Phi=\Phi_\I$ then $\A_\I(X)\subset \A_\Phi(X)$ and $\A_{\>\qc\I}(X)=\A_{\>\qc\Phi}(X)$
(see \ref{Gam and qc}, \ref{2Gammas}).\va1

\goodbreak
The next lemma gives conditions on the ringed space $X\<$ and the s.o.s.~$\Phi$ ensuring 
that $\A_\Phi(X)$  is a \emph{Serre subcategory} of~$\A(X)$  \cite[Tag 02MN\kf]{Stacks}, so that 
$\A_\Phi(X)$~\emph{is~plump in} $\A(X)$ (see section \ref{notation}).
Similarly, when $X$ is a scheme and $\I$ a finitary $\OX$-base, then $\A_\I(X)$ (resp.~$\A_{\>\qc\I}(X)$) is a Serre---hence plump---subcategory of~$\A(X)$ (resp.~$\Aqc(X)$).\looseness=-1

\begin{slem}\label{serre subcat} Let\/ $X$ be a ringed space, $M'\xto{\,f\,} M\xto{\,g\,} M''$  an exact sequence 
of\/ $\OX$-modules, $\Phi$ an s.o.s.~in\/ $X\<,$ and 
when\/ $X$ is a scheme, $\I$ a finitary\/ $\OX$-base.\va1

\textup{(i)} Suppose\/ $X$ has a base of quasi-compact open sets, and either that\/ $\Phi$ is finitary or 
that\/ $\Phi=\Phi_Y\ (Y\subset X)$. If\/ $M'$ and $M''$ are in\/ $\A_{\Phi}(X)$ then\/ $M\in\A_\Phi(X)$.

\textup{(ii)} When\/ $X$ is a scheme, if\/ $M'$ and $M''$ are in\/ $\A_\I(X)$ then\/ $M\in\A_\I(X)$. Hence
if\/ $M'$ and $M''$ are in\/ $\A_{\>\qc\I}(X)$ and\/ $M\in\Aqc(X)$ then\/ $M\in\A_{\>\qc\I}(X)$.

\end{slem}

\begin{proof}
(i). Fix an open $U\subset X$ and $m\in\Gamma(U,M)$. One needs that any $x\in U$ has an open neighborhood $V\subset U$ such that $\supp_V(m)\in\Phi|_V$. By assumption, \mbox{$x\in V\subset U$}
with $V$ quasi-compact and open, and such that $\supp_{V}(g(m))\subset Z'\cap V$ for some $Z'\in\Phi$. 

If $\Phi$ is finitary, then one can assume that
$X\setminus \<Z'$ is retrocompact in $X\<$, so that $V\setminus \<Z'$ is quasi-compact.
Over $V\setminus \<Z'$, $g(m)=0$, so $m\in\textup{im}\>(f)\<$, whence by \ref{soslocal},\looseness=-1 
\[
\supp_{V\setminus \<Z'}(m)= Z''\cap(V\setminus \<Z')\textup{ for some } Z''\in\Phi.
\] 
Thus $\supp_{V}(m)\subset (Z'\cup Z'')$, and so $\supp_V(m)\in\Phi|_V$.\va1

(ii). Fix an open $U\subset X$ and $m\in\Gamma(U,M)$. One needs, first, that each $x\in U$ has an open neighborhood $V\subset U$ over which $m$ is annihilated by some $I\in\I$. By~assumption, $x$ has an open neighborhood $V\subset U$ over which $g(m)$ is annihilated by some $I'\in\I$ with $I'|_{V}$ generated
by finitely many of its sections over $V\<$. Hence, over~$V\<$, $I'm\subset\textup{ im}\>f\<$,
so over some open neighborhood $V'\subset V\<$,  $I'm$ is annihilated by some $I''\in\I$. 
One can then take $V\set V'\<$, $I\set I''I'$. 

The last assertion follows at once.
\end{proof}

Upgrading to the derived level, let $\D_{\<\Phi}(X)\subset\D(X)$ be the
full subcategory spanned by the complexes whose homology\- modules
are all in~$\A_\Phi(X)$. 

Under the hypotheses of ~\ref{serre subcat}(i), $\A_\Phi(X)$ is plump in $\A(X)$, so
the exact homology sequence of a triangle entails that
$\D_{\<\Phi}(X)$ is a \emph{triangulated subcategory} of~$\D(X)$: if two
vertices of a $\D(X)$-triangle lie in $\D_{\<\Phi}(X)$ then so does the
third. 

Furthermore, if $\vG{\Phi}$ commutes with direct sums (see, e.g., \Pref{Gam and lim}(i)),
then $\D_{\<\Phi}(X)$ is a \emph{localizing subcategory} of
$\D(X)$, meaning here a triangulated subcategory closed under small
direct sums in~$\D(X)$.

Plumpness of $\A_\Phi(X)$ also implies that any complex in $\A_\Phi(X)$ is in $\D_{\<\Phi}(X)$.\va2

Similar statements hold, with $\I$ in place of $\Phi$, when $X$ is a scheme and $\I$ is
a finitary $\OX$-base. \va2
\end{subcosa}
\begin{sprop}\label{ringed space D_Phi}
 \textup{(i)} If\/ $(\R\vG{\Phi},\iota^{}_\Phi)$ is coreflecting in\/ $\D\set\D(X)$ or\/ $\Dpl(X)$ 
 \textup{(see~\ref{coreflect exams}(c)),} then
 $\D_{\<\Phi}(X)\cap \D$ is the essential image of\/~$\R\vG{\Phi}\colon\D\to\D$.
 
 \textup{(ii)} Similarly, if\/ $X$ is a locally noetherian scheme and\/ $\I$ an\/ $\OX$-base, then 
the essential image of\/~$\R\vG{\I}\colon\Dqc(X)\to\Dqc(X)$ 
is\/ $\D_{\I}(X)\cap\Dqc(X)$.
 
\end{sprop}

\begin{proof}
By~\ref{2RGams}, (ii) follows from (i). 
Application of the sentence preceding~\ref{coref and adj} to the coreflecting
pair $(\R\vG{\Phi},\iota^{}_\Phi)$ shows that
(i) results from the next lemma.
\end{proof}

\begin{slem}\label{gamadj}
For any s.o.s.~$\Phi$ in a ringed space\/ $X\<,$ an\/ $\OX\<$-\kf complex\/~$E\>$ 
lies in\/~$\D_{\<\Phi}(X)$ if and only if\/ $\iota^{}_\Phi(E\>)$
is an~isomorphism $\R\vG{\Phi} E\iso E\>$.
\end{slem}

\begin{proof}  Let $E\to I$ be a K-injective resolution. 
By \eqref{intersect vG}, $\vG{\Phi}I$ is a complex in~$\A_\Phi(X)$, so as noted above,
$\R\vG{\Phi}E\cong \vG{\Phi}I\in \D_{\<\Phi}(X)$, whence the essential image of $\R\vG{\Phi}$ is 
contained in~$\D_{\<\Phi}(X)$.\vs1

For the opposite inclusion it suffices to show that \emph{if\/ $E$---hence\/ $I\>$---is 
in\/ $\D_{\<\Phi}(X),$ then the natural map is an isomorphism $\R\vG{\Phi}E\iso E,$ that is, for every\/ $n\in\mathbb Z,$ the natural map\/ 
is an isomorphism} $H^n\<\vG{\Phi} I\iso H^n\<I$.

For any closed $Z\subset X\<,$ set $\vG{\<Z}\set\vG{\Phi_{\<\<Z}}$ and let $u^{}_{\<Z}\colon (X\setminus \<Z)\hookrightarrow X$ be the inclusion. Since $I$ is flabby, there is a natural exact sequence
\[
0\to\vG{\<Z}I\to I\to u^{}_{\<Z*}u_{\<Z}^*\mkern.5mu I\to 0
\]
whence an exact cohomology sequence
\begin{equation*}\label{cohseq}
\cdots\to H^n\<\<\vG{\<Z}I\to H^n\<I\to H^n u^{}_{\<Z*}u_{\<Z}^*\mkern.5mu I\to H^{n+1}\<\vG{\<Z}I
\to H^{n+1}I\to\cdots
\tag{\ref{gamadj}.1}
\end{equation*}
to which, by\eqref{HGam as lim.1}, application of the exact functor 
$\dirlm{\lift .15,{\!\!\halfsize{$Z$}\in\lift1,\Phi,},\!\!}\!$\va{-1} brings the  problem down to proving the next Lemma.
\end{proof}

\begin{slem} If\/ $J=(J^\bullet\<\<,d^{\>\bullet})\in\D_{\Phi}(X)$ is flabby,\va{-3} then\/ 
$\dirlm{\lift .15,{\!\!\halfsize{$Z$}\in\lift1,\Phi,},\!\!}\!H^nu^{}_{\<\<Z*}u_{\<\<Z}^* J=0$.
\end{slem}
\vskip1pt
\begin{proof} 
For $Z\in\Phi$, since $\A_{\>\Phi}\mkern-1.5mu(X)$ is plump in~$\A(X)$, 
therefore $\vG{\<Z}J\in \D_{\<\Phi}(X)$;
and  the exactness of \eqref{cohseq} with $J$ in place of $I\<$ shows that  
$u^{}_{\<Z*}u_{\<Z}^*\mkern.5mu J\in\D_{\<\Phi}(X)$.\vs1

Let $x$ be any point in $X\<,$ $V$ an open neighborhood of $x$, and 
$
h\in\Gamma(V, u^{}_{\<Z*}u_{\<Z}^*\mkern.5mu J^n)
$
such that $d^{\>n}h=0$. Since $u^{}_{\<Z*}u_{\<Z}^*\mkern.5mu J\in\D_{\<\Phi}(X)$,\vs{.5} 
$x$ has an open neighborhood $U\subset V$
where the element\vs{.3} 
$
\bar h\in\Gamma(U, H^n u^{}_{\<Z*}u_{\<Z}^*\mkern.5mu J\>)
$
given by $h$
is supported in a subset $Z'\cap U$ with $Z'\in\Phi$. Therefore, if $Z_1\set Z\cup Z'$ then the natural map
\[
\Gamma(U,H^n u^{}_{\<Z*}u_{\<Z}^*\mkern.5mu J\>)\to \Gamma(U,H^n u^{}_{\<Z_1*}u_{\<Z_1}^*J\>)
\] 
annihilates $\bar h$. Thus the stalk at $x$ of $\dirlm{\lift .15,{\!\!\halfsize{$Z$}\in\lift1,\Phi,},\!\!}\!H^n( u^{}_{\<Z*}u_{\<Z}^*\mkern.5mu J\>)$ vanishes.
\end{proof}
\vskip2pt

 The derived functor\/ $\R(\varGamma_{\<\!\Phi}^{\>0})$ is right-adjoint to the derived functor\-
 \[
 \boldsymbol j\set\R j\colon\D(\A_\Phi(X))\to\D(X),
 \] 
 see \cite[p.\,49, 5.2.2]{DFS} (in whose second line ``$\bj$ be the" should follow ``let"). And $\R\vG{\Phi}=\boldsymbol j\R(\varGamma_{\<\!\Phi}^{\>0})$. From \ref{gamadj} and \emph{loc.cit.\:}(2)$\>\Rightarrow$(1), one gets:

\begin{scor}\label{equivalence}
$\R(\varGamma_{\<\!\Phi}^{\>0})$ restricts to an\/ \emph{equivalence of categories} 
\[
\D_{\<\Phi}(X)\xto{\,\lift.5,\approx,\,}\D(\A_\Phi(X)),
\] 
with quasi-inverse given by $\boldsymbol j$.
\end{scor}

\centerline{* * * * *}

The \emph{support} $\Supp(E)$ of an $\OX$-complex $E$ is the set
of points at which $E$ is not exact, that is, the union of the
supports of all the homology sheaves of $E$.  

\begin{slem}\label{Supp in Y} 
For\/ $Y\subset X,\ \Phi_Y$ as in\/~\textup{\ref{supports},} and  $E\in\D(X),$ 
\[
\Supp(E)\subset Y\!\iff\! E\in\D_{\<\Phi_{\<Y}}\<\<(X).
\]
\end{slem}

\begin{proof}
This is a statement about the homology modules of $E$, so it suffices to note that 
for an $\OX$-module $M$, it follows directly from definitions that
\[
\Supp(M)\subset Y\!\iff\! M\in\A_{\>\Phi_Y}\<\<(X).\\[-16pt]
\]
\end{proof}

\begin{slem}\label{Supp in cup} For any s.o.s.~$\Phi$ in a ringed space\/ $X,$ and  $E\in\D(X),$
\[
\Supp(\R\vG{\Phi}E)\subset \bigcup_{Z\in\Phi}\mkern-4mu Z.
\]
\end{slem}

\begin{proof}
Since $E$ can be assumed to be K-injective, it suffices to note that since $\vG{\<Z}E$ vanishes outside $Z$,
therefore $\vG{\Phi} E= \dirlm{\lift .15,{\!\!\halfsize{$Z$}\in\lift1,\Phi,},\!\!}\!\vG{\<Z}E$ vanishes 
outside $\bigcup_{Z\in\Phi}\<Z$.
\end{proof}

\vskip8pt

\subsection{Idempotent pairs in symmetric monoidal categories}\label{idemclosed}
\stepcounter{thm} 
Part of the ``basic formal setup,'' a category-theoretic framework for duality, 
local and~global, to be built on in subsequent chapters, is the notion of \emph{idempotent pair} in a 
symmetric monoidal category $\D$---more precisely, in the slice category $\D/\cO$ with $\cO$ the unit object  (Definition~\ref{idempotent pair}).%
\footnote{\kf $\D/\cO$ has as objects the pairs $(C,\gamma)$ with $C$ an object of $\D$ and 
$\gamma\colon C\to\cO$ a $\D$-map, and as morphisms 
$\lambda\colon(B,\beta)\to(A,\>\alpha)$ the 
$\D$-morphisms \mbox{$\lambda_0\colon B\to A$} such that $\beta=\alpha\lambda_0$.
(Henceforth, absent potential for confusion we will not differentiate
notationally between $\lambda$ and $\lambda_0$.) } 
This notion is equivalent to that of \emph{$\>\otimes$-coreflection,} 
that is, coreflection $(\Gamma,\iota)$ with $\Gamma$  isomorphic to a functor 
$\Gamma_{\<\<\!A}(-)\set A\otimes-$ where $A$ is a fixed object and $\otimes$ is the monoidal product (Proposition~\ref{idpt3}).
This section and the following two review some basics about such pairs.

\begin{sdef}[\,{\cite[p.\,251\kern.4pt{ff}\kern1pt]{Ma2}}\,]\label{def monoidal}
A (\emph{symmetric}) \emph{monoidal category}
$$
\mathbf D=(\mathbf{D_0},\ot,\cO,\ma,\ml,\mr,\ms)
$$
consists of a category $\mathbf{D_0}$, a  ``product" functor
$\ot\colon\mathbf{D_0}\times\mathbf{D_0}\to\mathbf{D_0}$, an
object $\cO$ of~$\mathbf{D_0}$,
and functorial isomorphisms (for $A, B, C$  in $\mathbf {D_0}$)
\begin{align*}
\mkern-50mu\ma=\ma_{A,B,C}\colon (A\ot B)\ot C
   &\iso A\ot (B\ot C) 
	\tag{associativity}\\
\ml=\ml^{}_{\<A}\colon\cO\ot A\iso A\ \
   &\qquad \quad\; \mr=\mr^{}_{\!A}\colon A\ot\cO\iso A \tag{units} \\
\ms=\ms_{A,B}\colon A\ot B
   &\iso B\ot A \tag{symmetry}
\end{align*}

\noindent
such that
$\ms\smcirc\ms=\b1$ (identity map) and the following diagrams commute:\vs{8} 
\end{sdef}

\begin{gather*}\tag{\ref{def monoidal}.1}\label{mon1}\\[-38pt]
(A\ot\cO)\ot B
 \kern .7em\overset{\ma}
{\underset{\vbox to0pt{\vskip12pt
\hbox to 0pt{\hss$\scriptstyle\mr\>\ot\>\b1\hskip13ex\b1\>\ot\>\ml$\hss}\vss}}
 {\hbox to 5.4em{\rightarrowfill}}}\kern .7em
  {A\ot(\cO\ot B)}	 \\[-12pt]
\centerline{$\rotatebox{-45}{\hbox to 38pt{\rightarrowfill}}\mkern 40mu
\raisebox{4pt}{\rotatebox{-135}{\hbox to 38pt{\rightarrowfill}}}$}\\[-3pt]
\centerline{$\<\<A\ot \<B$}\\
\end{gather*}

\vs{-20}
\begin{equation*}\label{mon2}
\minCDarrowwidth.36in
\CD
((A\ot B)\ot C)\ot D
 @>\ma>>
  (A\ot B)\ot (C\ot D)
   @>\ma>>
    A\ot(B\ot(C\ot D))\\
@V\ma\>\>\ot\kf \b1VV
  @.
    @VV\b1\>\ot\>\>\ma V\\
(A\ot(B\ot C))\ot D
 @.
  \underset{\ma}{\Rarrow{14.7em}}
   @.
    A\ot((B\ot C)\ot D)
\endCD\tag*{\raisebox{4pt}{(\ref{def monoidal}.2)}}
\end{equation*}

\vs{10}
\begin{equation*}\label{mon3}
\minCDarrowwidth.36in
\CD
(A\ot B)\ot C
 @>\ma>>
  A\ot(B\ot C)
   @>\ms>>
    (B\ot C)\ot A\\
@V \ms\>\>\ot\kf 1VV
  @.
    @VV\ma V\\
(B\ot A)\ot C
 @>>\ma>
  B\ot(A\ot C)
   @>>\b1\>\ot\>\>\ms>
    B\ot(C\ot A)
\endCD\tag{\ref{def monoidal}.3}
\end{equation*}

\vs{20}
\begin{gather*}\tag{\ref{def monoidal}.4}\label{mon4}\\[-34pt]
\centerline{$
\underset{\vbox to 0pt{\vskip 10pt\hbox to 0pt{\hss\hskip9ex
$\scriptstyle\mr\ $\hss}\vss}}
{A\ot \cO}
 \kern .55em
\overset{\ms}{\hbox to 4.3em{\rightarrowfill}}\kern .6em
\underset{\vbox to 0pt{\vskip 9pt\hbox to
 0pt{\hss$\hskip-12pt\scriptstyle\ml$\hss}}}\cO\ot A
$} \\[-12.4pt]
\centerline{$\rotatebox{-45}{\hbox to 38pt{\rightarrowfill}}\mkern 10mu
\raisebox{4pt}{\rotatebox{-135}{\hbox to 38pt{\rightarrowfill}}}$}\\[-4pt]
\centerline{$A$}
\end{gather*}

Necessarily,  the following diagrams commute\vs{-9} too \cite[p.\,165, Exercise 1]{Ma2}.

\begin{gather*}\tag{\ref{def monoidal}.5}\label{mon5}
(A\ot B)\ot \cO
 \kern .7em\overset{\!\!\ma}
{\underset{\vbox to0pt{\vskip12pt\hbox to
 0pt{\hss$\mspace{18mu}\scriptstyle\mr\hskip13.5ex \b1\ot\>\mr$\hss}\vss}}
 {\hbox to 3.3em{\rightarrowfill}}}\kern .7em
  {A\ot(B\ot\cO)}
\qquad\quad
(\cO\ot A)\ot B
 \kern .7em\overset{\,\>\ma}
{\underset{\vbox to0pt{\vskip12pt
\hbox to 0pt{\hss$\mspace{-23mu}\scriptstyle\ml\>\ot\b1\hskip13ex\ml$\hss}\vss}}
 {\hbox to 3.3em{\rightarrowfill}}}\kern .7em
  {\cO\ot(A\ot B)}
\\[-12pt]	
\mkern-3mu\rotatebox{-35}{\hbox to 38pt{\rightarrowfill}}\mkern 40mu
\raisebox{4pt}{\rotatebox{-145}{\hbox to 38pt{\rightarrowfill}}}
\mspace{180mu}
\rotatebox{-35}{\hbox to 38pt{\rightarrowfill}}\mkern 40mu
\raisebox{4pt}{\rotatebox{-145}{\hbox to 38pt{\rightarrowfill}}}\\[-2pt]
A\ot B\ \mspace{285mu} A\ot B\ 
	 \\[-8pt]
\end{gather*}

\begin{exams}\label{exams-mon}
(a) Let $X$ be a ringed space.  
Derived tensor product makes $\D(X)$ into a monoidal category, with unit object
$\OX$ (\ref{derived tensor} below); and similarly for $\Dqc(X)$ (resp.~$\Dqct(X)$)
when $X$ is a scheme (resp.~noetherian formal scheme), see~\ref{monoid1}.\va1

\pagebreak[3]
(b) For a monoidal category $\D$, the slice category $\D/\cO$ has a monoidal structure 
with unit object~$(\cO,\b1_\cO)$, 
product $(A,\alpha)\>\,\underline{\!\otimes\!}\,\>(B,\beta)\set(A\otimes B,\mu\smcirc(\alpha\otimes\beta))$
where $\mu=\mr_\cO=\ml_\cO$ (see proof of~\ref{idpt1}), and  isomorphisms $\ma$, $\ml$, $\mr$ 
and~$\ms$ whose images under the functor $(A,\alpha)\mapsto A$ are the the corresponding isomorphisms in~$\D$. 
(Details are left to the reader.)
\end{exams}

Until otherwise indicated, $\D$ will be a fixed monoidal category. Sometimes, for simplicity, $A\ot\cO$ and
$\cO\ot A$ will be identified---harmlessly---with the object $A\in\D$.

\begin{sdef}\label{idempotent pair}
A \emph{$\D$-idempotent pair} $(A,\alpha$) is a $\D$-map
$\alpha\colon A\to\cO$ such that the  composite maps\va1
$
A\ot A \xto{\lift.9,\<\<\b1\ot\>\alpha\>\>,}A\ot\cO\xto{\lift.65,\ \mr,\ }A
$
and
$
A\ot A \xto{\<\alpha\>\ot \b1\>}\cO\ot A\xto{\ \ml\ }A
$
are \emph{equal isomorphisms}. An object $A\in\D$ is idempotent if such an $\alpha$ exists.

\end{sdef}

\begin{exams}\label{exams-idem}

(a) Let\/ $X$ be a locally noetherian scheme, and\/ $\I$ an\/~$\OX$-base.
The pair $(\R\vG{\I}\OX,\iota^{}_\I(\OX))$ is\/ $\Dqc(X)$- and $\D(X)$-idempotent (see \ref{RGamidem}).\va2

(b) Let\/ $X$ be a scheme, and\/ $\Phi$ a finitary s.o.s.\ in~$X\<$.
The pair $(\R\vG{\Phi}\OX,\iota^{}_\Phi(\OX))$ is $\Dqc(X)$- and $\D(X)$-idempotent (see \ref{RGamidem}).\va2

For additional such examples, involving topological rings, or noetherian formal\- schemes, see \Cref{example},
\Pref{A and Z} and  \Cref{idem and bases}. \va2

(c) Let $\D$ be a category with a terminal object $\cO$, and such that
any two objects $A,B\in\D$ have a product, denoted $A\ot B\<$.
With this $\cO$, $\ot$ (made into a functor), and obvious choices for $\ma$, $\ml$, $\mr$ and $\ms$,
one gets a monoidal category. One verifies, for any $A\in\D$ with  $\alpha\colon A\to\cO$ the unique map, that $(A,\alpha)$ is idempotent if and only~if $\alpha$~is a monomorphism (that is, for any $B\in\D$ there is at most one\va1
map $B\to A$).\looseness=-1

(d) In particular, let $(\D,\le)$ be a preordered set (set with a reflexive, transitive binary relation $\le$), 
considered as a category in the usual way: the objects are the elements of $\D$,
there is a unique map $A\to B$ if $A\le B$, and otherwise no such map
at all. Assume that $\D$ has a largest object~ $\cO$, 
and that any two objects $A,B\in\D$ have a greatest lower bound (=\:product), 
denoted $A\ot B$. With the obviously unique $\ma$, $\ml$, $\mr$ and $\ms$,
$\D$ is a monoidal category in which for \emph{any} object~$A$ with
$\alpha\colon A\to\cO$ the~unique map, $(A,\alpha)$ is idempotent.

Note that $B\le A\!\iff\! \alpha\ot \b1\colon A\ot B\to\cO\ot B\cong B$ \emph{is an isomorphism}.

Also, $B\cong A\!\iff\! B\le A \textup{\emph{ and }}A\le B$.\vs1
 
Such categories will be called \emph{preordered} monoidal categories. They can be viewed as small monoidal categories in which for any objects $A$ and $B$, there exists at most one map $A\to B$, and exactly one if $B=\cO$ or if $B=A\otimes A$.
\end{exams}

\begin{srem}\label{idpt and D/O}
The full subcategory $\mathbf I_\D$ of $\>\>\D\</\<\cO$ spanned by the $\D$-idempotent pairs 
is \emph{strictly} full: use the fact that 
if\/ $(A,\alpha)$ is idempotent and\/ $\lambda\colon B\iso A$ is a $\D$-isomorphism, 
then $(B, \alpha\lambda)$ is idempotent. 
Moreover, $\mathbf I_\D$ is a \emph{preordered monoidal} subcategory of~$\D/\cO$, see
\ref{idpt1}, \ref{idpt tensor} and \ref{one map} below. 

Note that $(A,\alpha)$ is $\D$-idempotent $\Leftrightarrow((A,\alpha), \,\alpha)$ is $(\D/\cO)$-idempotent.
\end{srem}

\enlargethispage*{-10pt}
\begin{slem}\label{idpt1}
The pair\/ $(\cO,\b1_\cO)$ is\/ $\D$-idempotent.
\end{slem}

\begin{proof}
The assertion means that the unit isomorphisms
$\ml=\ml_\cO\colon \cO\ot\cO\iso \cO$ and
$\mr=\mr_\cO\colon \cO\ot\cO\iso \cO$ are the same, or,
by \eqref{mon4}, that 
\emph{the automorphism\/ 
$\ms=\ms_{\>\cO\<,\cO}\colon\cO\ot\cO\iso\cO\ot\cO$ is the identity map.}\va1

\pagebreak[3]
By \eqref{mon3} (with $(A,B,C)$ replaced by $(\cO,A,\cO)$), the border of the 
following diagram of isomorphisms
commutes, for any $A\in \D$:
\vskip-2pt
\hspace{21pt}
\begin{minipage}[t][142pt][t]{4in}
$$
\begin{CD}
\underset{\raisebox{3pt}{\UnderElement{\ma}{\downarrow}{34pt}{}}}
{(\cO\ot A)\ot\cO}@.\mspace{50mu}@.
    \overset{\ms\>\ot \b1_{}}{\Rarrow{5em}}\mspace{50mu} @. 
\underset{\raisebox{3pt}{\,\UnderElement{}{\|}{34pt}{}}}    
{(A\ot\cO)\ot \cO}
@>\ma>\under{-6.2}{\displaystyle\square}>
\underset{\raisebox{3pt}{\UnderElement{}{\downarrow}{34pt}{\!\<\b1\ot\>\ms}}}
{A\ot(\cO\ot\cO)}
               \\[-7pt]
\mspace{50mu}\hbox to
0pt{$\mspace{20mu}\raisebox{14pt}{\rotatebox{-45}{\hbox to
29pt{\rightarrowfill}}}\raisebox{8pt}{$\mspace{-17mu} 
  \scriptstyle\ml^{}_{\<A}\ot\b1$}$\hss}
@.
\hbox to 0pt{$\mspace{66mu}\raisebox{18.3pt}{\rotatebox{-135}{\hbox
to 29pt{\rightarrowfill}}}\raisebox{8pt}{$\mspace{-54mu} 
  \scriptstyle\mr^{}_{\<\<A}\ot \b1$}$\hss}
                \\
\end{CD}
$$
$$
\mspace{-173mu}
\begin{CD}
\ A\ot\cO\\
\end{CD}
$$
$$
\begin{CD}
\mspace{50mu}\hbox to
0pt{$\mspace{16mu}\raisebox{17pt}{\rotatebox{45}{\hbox to
29pt{\rightarrowfill}}}\raisebox{23pt}{$\mspace{-18mu} 
  \scriptstyle\ml^{}_{\<\<A\ot\cO}$}$\hss}
@.
\hbox to 0pt{$\mspace{57mu}\raisebox{20.3pt}{\rotatebox{135}{\hbox to
29pt{\rightarrowfill}}}\raisebox{23pt}{$\mspace{-50mu} 
\scriptstyle\mr^{}_{\<\<A\ot\cO}$}$\hss}\\[-26pt]
\mspace{8mu}\cO\ot(A\ot\cO)@.\mspace{48mu}@. 
\underset{\under{1.2}{\ms}}{\Rarrow{5em}}\mspace{50mu} 
@. \,(A\ot\cO)\ot \cO\,
 @>>\under{1.2}\ma> A\ot(\cO\ot\cO)
\end{CD}
$$
\end{minipage}

\vspace{-32pt}\noindent
The subtriangle on the left commutes by \eqref{mon5} (second diagram).
The ones  at the top and bottom commute by \eqref{mon4}. Furthermore, 
$\mr^{}_{\!A}\ot\b1=\mr^{}_{\!A\ot \cO}\>$,
as shown by the next diagram,
which commutes because $\mr$ is functorial:
\[
\begin{CD}
(A\ot\cO)\ot\cO @>\mr^{}_{\<\<A}\ot\b1>> A\ot\cO \\
@V\mr^{}_{\<\<A\ot\cO} VV @V\simeq V\mr^{}_{\!A} V \\
A\ot\cO @>\Iso>\under{1.2}{\mr^{}_{\!A}} > A
\end{CD}
\] 
It follows that
subrectangle $\square$ commutes, whence $\b1\ot\ms$ is the identity map,
whence so is $\ms$, as one sees by taking $A=\cO$ and applying
$\ml_{\cO\ot\cO}$.
\end{proof}

\begin{srem}\label{s=1} 
\emph{For any idempotent pair\/ $(A,\alpha)$,\va{.5} the symmetry
automorphism\/  $\ms^{}_{\!A,A}\colon A\ot A\iso A\ot A$ 
 is the identity map.}\va{.5} Indeed, \eqref{mon4} shows that
the following diagram---whose rows compose to 
the same isomorphism---commutes:
$$
\begin{CD}
A\ot A @>\b1\>\ot\>\alpha>> A\ot\cO @>\mr^{}_{\!A}>> A\\
@V\ms^{}_{\!A\<,A} VV @VV\ms^{}_{\!A\<,\cO}V @|\\
A\ot A @>>\alpha\>\ot\>\b1>\cO\ot A @>>\under{1.3}\ml_{\<\<A}> A
\end{CD}
$$

\emph{More generally,} by \ref{idpt and D/O} and \ref{idpt tensor} below, 
\mbox{$(A,\alpha)\>\,\underline{\!\otimes\!}\,\>(A,\alpha)$} is $(\D/\cO)$-idempotent, and so by
\ref{one map}, its only $(\D/\cO)$-endomorphism is the identity map.
\end{srem}

\begin{slem}\label{f*idpt}
Let\/ $\xi\colon \D_2\to\D_1$ be a functor between monoidal categories having respective units $\cO_2,\,\cO_1$ and product functors $\ot_2,\,\ot_1$.   Let\/ $(B,\beta)$ be\/ $\D_2$-idempotent.  

Suppose there exists 
a\/ $\D_1\<$-map\/ $u\colon \xi\cO_2\to\cO_1$ and a bifunctorial\/ $\D_2$-isomorphism\/ 
\[
v(E\<,F\>)\colon \xi E\ot_1\<\xi F\iso \xi(E\ot_2 F)\qquad (E,F\in\D_2)
\]
such that subdiagrams\/ \textup{\circled1} and\/ \textup{\circled4} of the following natural diagram commute.\va3
\[
\def\1{$\xi\<(B\ot_2\<B)$}
\def\2{$\xi\<(\cO_2\ot_2\<B)$}
\def\3{$\cO_1\ot_1\xi B$} 
\def\4{$\xi B\ot_1\xi B$}
\def\5{$\xi\cO_2\ot_1\xi B$} 
\def\6{$\xi B$} 
\def\7{$\xi\<(B\ot_2\<\cO_2)$}
\def\8{$\xi B\ot_1\<\xi\cO_2$}
\def\9{$\xi B\ot_1\cO_1$}
 \bpic[xscale=3, yscale=1.25]

   \node(11) at (1.9,-3){\1};
   \node(12) at (2.9,-2){\2}; 
   \node(13) at (4,-3){\6};

   \node(31) at (1,-3){\4};  
   \node(32) at (1.8,-1){\5};
   \node(33) at (3.1,-1){\3};
   
   \node(42) at (2.9,-4){\7};
    
   \node(52) at (1.8,-5){\8};
   \node(53) at (3.1,-5){\9};

    \draw[->] (11)--(12) node[above=-3, midway,scale=.75]{$\xi\<(\beta\ot_2\<\b1)\mkern95mu$}
                                   node[below, midway,scale=.75]{$\mkern10mu\simeq$} ;  
    \draw[->] (12)--(13) node[above=-.5, midway, scale=.75]{$\rotatebox{-22}{$\Iso$}$} ; 
    
    \draw[->] (31)--(32) node[above=1,midway,scale=.75]{$\xi\beta\ot_1\<\b1\mkern75mu$} ;
    \draw[->] (32)--(33) node[above=1, midway, scale=.75]{$u\ot_1\<\b1$} ;
    \draw[->] (52)--(53) node[below=1,midway,scale=.75]{$\b1\ot_1\<u$} ;
 
    \draw[->] (31)--(11) node[above, midway, scale=.75]{$\Iso$} 
                                   node[below,midway,scale=.75]{$v$};  
  
   
   \draw[<-] (13)--(33) node[above=1,midway,scale=.75]{$\mkern20mu\simeq$} ; 
   
   \draw[<-] (42)--(52) node[above=-.5, midway, scale=.75]{\rotatebox{22}{$\Iso$}} 
                                   node[below,midway,scale=.75]{$\mkern12mu v$} ; 
   
   \draw[->] (11)--(42) node[below=-3, midway,scale=.75]{$\xi\<(\<\b1\ot_2\<\beta)\mkern95mu$}
                                  node[above, midway,scale=.75]{$\mkern10mu\simeq$} ;  
   \draw[->] (42)--(13) node[above=-.5, midway, scale=.75]{$\rotatebox{20}{$\Iso$}$} ; 
   \draw[<-] (12)--(32) node[above=-.5, midway, scale=.75]{\rotatebox{-22}{$\Iso$}} 
                                   node[below=-.5,midway,scale=.75]{$v\mkern 15mu$} ; 
   \draw[->] (31)--(52) node[below=-2,midway,scale=.75]{$\b1\ot_1\<\<\xi\beta\mkern75mu$} ;  
   \draw[->] (53)--(13) node[below=1,midway,scale=.75]{$\mkern20mu\simeq$} ;
   
  \node at (3.05,-1.55)[scale=.85]{\textup{\circled1}} ;
  \node at (1.9,-2)[scale=.85]{\textup{\circled2}} ;
  \node at (1.9,-4)[scale=.85]{\textup{\circled3}} ;
  \node at (3.05,-4.55)[scale=.85]{\textup{\circled4}} ;
  \node at (2.9,-3)[scale=.85]{\textup{\circled5}} ;
  \epic
\]
\vskip-2pt\noindent
Then $(\xi B, u\smcirc\xi \beta)$ is\/ $\D_1$-idempotent.
\end{slem}

\begin{proof}
The commutativity of subdiagrams \circled2 and \circled3 is given by the functoriality of~$v\>$; and that of~
 \circled 5 holds by the idempotence of $(B,\beta)$.
These commutativities, plus those of \circled1 and \circled4, imply that the border of the diagram commutes, and consists entirely of isomorphisms, whence the conclusion.\va{-2}
\end{proof}

\begin{srem}\label{f*idpt2}
The hypotheses in \ref{f*idpt} are satisfied if, $f\colon X_1\to X_2$ being a map of ringed spaces, 
$\xi$ is $\bL f^*\colon\D(X_2)\to\D(X_1)$ and $u,\:v$ are the natural isomorphisms:
commutativity of subdiagram \circled1 follows by the \emph{duality principle} \cite[p.\,106]{Dercat} from that of the first diagram in \cite[p.\,103, (3.4.2.2)]{Dercat}, and that of \circled4 is shown similarly.

For another instance, see \ref{Gam preserves idem} below.
\end{srem}

\subsection{Idempotent pairs and $\boldsymbol\otimes$-coreflections}\label{idem and tensor}
\stepcounter{thm} 
The main results in this section are \ref{idpt3} and~\ref{coreflections}, whose corollary, \ref{RGamidem}, motivates much of the subsequent approach to duality theory.\va2

Fix a symmetric monoidal category $\D=(\mathbf{D_0},\ot,\cO,\ma,\ml,\mr,\ms)$.\va3

Sending an object $A\in \D$ to the natural functor $\Gamma_{\<\<\!A}\colon\D\to\D$ taking $F$ to~ $A\ot F$ 
gives an \emph{equivalence} from the category $\D$
to the category of  
\mbox{\emph{$\ot$-endofunctors of}~$\D$,} that is, 
those functors\/ $\Gamma\colon\D\to\D$ such that there 
exists a functorial isomorphism 
$\Gamma\cO\ot  F\iso \Gamma F.$ 
There is a quasi-inverse equivalence taking\/ $\Gamma$ to\/ $\Gamma\cO$.\looseness=-1

These quasi-inverse equivalences lift to quasi-inverse equivalences between $\D/\cO$ and the category 
$\mathbf E^\ot$ of pairs $(\Gamma,\iota)$ with
$\iota\colon \Gamma\to {\bf 1_\D}$ a map of endofunctors of $\D$
such that there exists a functorial isomorphism 
\begin{equation}\label{PsiE}
\psi(F)\colon\Gamma\cO\ot F\iso \Gamma F\qquad(F\in\D)\\[-1pt]
\end{equation}
making the following diagram commute:\va2
\begin{equation}\label{tensor iso}
\begin{CD}
 \Gamma\cO\ot F@>\Iso>\under{1.2}{\psi(F)}>\Gamma F\\[-3pt]
 @V\iota(\cO)\>\ot\>\b1VV@VV\iota(F)V \\
\cO\ot F @>\Iso>\under{1.25}{\ml^{}_{\<\<F}}>  F
 \end{CD}
\end{equation}
\vskip2pt

The lifted quasi-inverse equivalences act (objectwise) as follows: 
\steq
\begin{equation*}\label{equiv}
\begin{aligned}
(A,\alpha)&\mapsto \big(\Gamma_{\<\<\!A}, \>\iota_\alpha\colon \Gamma_{\<\<\!A}E=A\ot E 
\xto{\!\alpha\otimes\b1\>} \cO\ot E\xto{\ml_{\<E}} E\big)\qquad (E\in \D);\\[3pt]
(\Gamma,\>\iota)&\mapsto \big(\Gamma\cO,\>\iota(\cO)\big).\vs1
\end{aligned}\tag{\theequation}
\end{equation*}
(Note that if $\psi(F)$ is the natural functorial isomorphism $(A\ot \cO)\ot F\iso A\ot F$, then 
with $\Gamma\set\Gamma_{\<\<\!A}$ and $\iota\set\iota_\alpha$, the diagram \eqref{tensor iso}
commutes.)\va1

As $\ml^{}_\cO=\mr^{}_\cO$ (see \ref{idpt1}), it  is straightforward to verify that \eqref{equiv} does give rise naturally to quasi-inverse functors, that is, there are functorial isomorphisms\va{-2}
$$
(\Gamma^{}_{\<\<\Gamma\cO},\>\iota_{\iota(\cO)})\iso (\Gamma\<,\iota),
\qquad \big(\>\Gamma_{\<\<\!A}\cO,\>\iota_\alpha(\cO)\big)\iso (A,\alpha).
$$

The monoidal structure on 
$\mathbf E^\ot$ corresponding under this lifted equivalence to the one on~$\D/\cO$ mentioned in 
\ref{exams-mon}(b) has product $\bar\ot$ such that
\[
(\Gamma,\iota)\,\bar\ot\, (\Gamma',\iota')=(\Gamma\smcirc\Gamma'\!,\>\> \iota\smcirc\Gamma(\iota')).
\]

\begin{sprop}\label{Psi for coref} 
For any $(\Gamma,\iota)\in\mathbf E^\ot$  and\/ $E,\>F\in\D$ there are  isomorphisms\/ 
$$
\Gamma\<E\ot F\underset{\psi(E,F)}\iso\Gamma(E\ot F\>)
\underset{\psi'(E,F)}\osi E\ot \Gamma\<F
$$  
making the following diagram commute$\>:$
\[
\def\1{$\Gamma E\ot F$}
\def\2{$\Gamma(E\ot F\>)$}
\def\3{$E \ot \Gamma F$}
\def\4{$E\ot F$}
 \bpic[xscale=3.5, yscale=2]

   \node(11) at (1,-1){\1};
   \node(12) at (2,-1){\2};
   \node(13) at (3,-1){\3}; 
   
   \node(22) at (2,-2){\4};
     
   \draw[->] (11)--(12) node[above, midway,scale=.75]{$\psi(E,F)$} ;  
    
   \draw[->] (13)--(12) node[above, midway,scale=.75]{$\psi'(E,F)$} ;
   
   \draw[->] (12)--(22) node[left, midway, scale=.75]{$\iota(E\ot F)^{}\<\<$} ;  

   \draw[->] (11)--(22) node[below=-5, midway, scale=.75]{$\iota(E)\<\ot\< \b1_{\<F}\mkern110mu$} ;  
   \draw[->] (13)--(22) node[below=-5, midway, scale=.75]{$\mkern110mu\b1_{\<E}\<\ot\< \iota(F)$} ;
 
  \node at (1.65,-1.3)[scale=.85]{\textup{\circled1}} ;
   \node at (2.35,-1.3)[scale=.85]{\textup{\circled1}$'$} ;
  \epic
\]

If, moreover, $(\Gamma,\iota)$ is a coreflection of\/ $\D,$ then\/ $\psi(E,F)$ and\/ $\psi'(E,F)$ are unique.\vs{1.5}
 \end{sprop}
 
 \begin{proof} 
The easily-checked (via \eqref{mon5} and \eqref{tensor iso}) commutativity of the following natural diagram shows that the composite isomorphism
\[
\psi(E,F\>)\colon\Gamma\<E\ot F\underset{\under{1}{\eqref{PsiE}}}\iso
\,(\Gamma\cO\ot E)\ot F \underset{\under{1}{\ma}\>\>}\iso
\Gamma\cO\ot(E\ot F\>)\>\underset{\under{1}{\eqref{PsiE}}}\iso
\,\Gamma(E\ot F\>)
\]
\vskip1pt\noindent
makes subdiagram \circled1 commute.\va1 
It follows that the natural composite isomorphism
\[
\psi'(E,F\>)\colon E\ot\Gamma F\iso \Gamma F\ot E\xto{\psi(F\<,E)}\Gamma(F\ot E)\iso \Gamma(E\ot F\>) 
\]
makes \circled1$'$ commute.
\goodbreak
\[
\def\1{$(\Gamma\cO\ot E)\ot F$}
\def\2{$\Gamma\cO\ot(E\ot F\>)$}
\def\3{$(\cO\ot E) \ot F$}
\def\4{$\cO^{\mathstrut}\ot(E\ot F\>)$}
\def\5{$\Gamma\<E\ot F$}
\def\6{$E\ot F$} 
\def\7{$\Gamma(E\ot F\>)$} 
 \bpic[xscale=4, yscale=1.35]

   \node(11) at (1,-1){\1};
   \node(13) at (3,-1){\2}; 
   
   \node(21) at (1.5,-2){\3}; 
   \node(22) at (2.5,-2){\4};
   
   \node(31) at (1,-3){\5};  
   \node(32) at (2,-3){\6};
   \node(33) at (3,-3){\7};

    \draw[->] (11)--(13) node[above, midway,scale=.75]{${\ma}$} ;  
    
    \draw[->] (21)--(22) node[above, midway,scale=.75]{${\ma}$} ;

    \draw[->] (31)--(32) node[below=1pt,midway,scale=.75]{$\b1_E\<\ot\iota(F\>)$} ;
    \draw[->] (33)--(32) node[below=1pt,midway,scale=.75]{$\iota(E\ot F\>)$} ;

   \draw[->] (31)--(11) node[left=1pt, midway, scale=.75]{$\psi^{\<-\<1}(E)\ot\b1_{\<F}$} ;  
  
   \draw[->] (13)--(33) node[right=1pt, midway, scale=.75]{$\psi(E\ot F\>)$} ;  
   
   \draw[->] (11)--(21) ;  
   \draw[->] (21)--(32) ;
   \draw[->] (22)--(32) ; 
   \draw[->] (13)--(22) ;  
  \epic
\]

When $(\Gamma,\iota)$ is a coreflection, the unicity of $\psi$ and $\psi'$ follow from  \ref{coconds}
(with both $F$ and~$G$ replaced by $E\ot F$).\vs{1.5}
\end{proof}

\begin{srem} \label{PsiF and PsiOF}
One checks that a map $\psi(F)$ makes \eqref{tensor iso} commute if and only if  
$\psi(F)=\Gamma\>\>\ml^{}_{\<\<F}\smcirc\psi(\cO\<,\>F\>)$
for some $\psi(\cO,F)$ as in \ref{Psi for coref}. Hence when such a $\psi(\cO,F)$ is unique then so is
such a $\psi(F)$.
\end{srem}

\centerline{* * * * *}

\begin{sdef}\label{ot type}
A pair $(\Gamma,\iota)$ with\/ $\Gamma\colon\D\to\D$ a functor and\/ $\iota\colon\Gamma\to{\bf1}$ a map of functors is  a \emph{$\ot$-coreflection} \emph{of}~$\>\D$ 
(or $\ot$-\emph{coreflecting in} $\D$) if it is a coreflection of~$\>\D$ that lies in $\mathbf E^\ot\<$.
 The functor $\Gamma$ is a \emph{$\ot$-coreflector} if there exists an
 $\iota$ such that $(\Gamma,\iota)$~is a $\ot$-coreflection.
\end{sdef}

\begin{sprop}\label{idpt3} 
The above equivalence between\/ $\D/\cO$ and\/ $\mathbf E^\ot$ \textup{(see \eqref{equiv})}
induces an equivalence between the category\/~$\bf I_\D$  of\/ $\D$-idempotent pairs and the category of\/
\mbox{$\ot$-coreflections} of\/ $\D$.
\end{sprop}

\begin{proof}
To be shown is that $(A,\>\alpha)$ is idempotent if and only if
$(\Gamma,\iota)\set(\Gamma_{\<\<\!A},\iota^{}_{\alpha})$ is a $\ot$-coreflection.

Suppose first that $(A,\>\alpha)$ is idempotent. As before, $\Gamma\set\Gamma_{\<\<\!A}$ is 
a $\ot$-endofunctor of $\D$.
That $(\Gamma,\iota)$ is coreflecting means that for
any $E\in\D$, the following
diagram commutes and moreover, the maps $\bf1\ot(\alpha\ot\bf1)$ and $(\alpha\ot\b1)\ot\b1$ are isomorphisms:
$$
\begin{CD}
A\ot(A\ot E) @>\alpha\>\ot\>(\b1\ot\b1) >>\cO\ot(A\ot E)   \\
@V\bf1\ot(\alpha\ot\bf1) VV @VV\>\ml^{}_{A\ot E} V \\
A\ot (\cO\ot E) @>>\under{1.2}{\b1\ot\>\ml^{}_{\<\<E}}>A\ot E\vs3
\end{CD}
$$
Using \eqref{mon1}  and the functoriality of $\ma$, one expands this diagram as 
$$
\begin{CD}
A\ot(A\ot E)  @>\ma^{-\<1}>\under{-2}{\displaystyle\square_1}>(A\ot A) \ot E
    @>(\alpha\>\ot\>\b1)\>\ot\b1>\under{-2}{\displaystyle\square_2}>(\cO\ot A)\ot E   
    @>\ma>\under{-2}{\displaystyle\square_3}>\cO\ot(A\ot E)   \\
@V\bf1\ot(\alpha\>\ot\bf1) VV @V(\b1\>\ot\>\alpha)\! V\!\ot\>\b1V
@V\under{1}{\ml^{}_{\!A}}\ot\!V\!\<\b1V @VV\ml^{}_{\<A\ot E} V \\
A\ot (\cO\ot E)@>>\under1\ma^{-\<1}>(A\ot \cO)\ot E 
    @>>\under{1.2}{\mr^{}_{\<\<A}\>\ot\>\>\b1}>A\ot E@=A\ot E\vs3
\end{CD}
$$
The top row consists entirely of isomorphisms, so $(\alpha\ot\b1)\ot\b1$ is an isomorphism. 
The commutativity of square~$\square_1$ holds because
$\ma$ is functorial,  and since $1\ot\alpha$ is an isomorphism, therefore so is $\bf1\ot(\alpha\ot\bf1)$.
The~commutativity of $\square_2$ holds by idempotence of $(A,\>\alpha)$,
and of $\square_3$ by \eqref{mon5}. So $(\Gamma,\iota)$ is indeed 
$\ot$-coreflecting.\vs1

Suppose, conversely, that  $(\Gamma,\iota)$ is a $\ot$-coreflection. What's needed is that
the maps $p\set\ml^{}_{\>\Gamma\<\cO}\smcirc(\iota(\cO)\>\ot\>\b1)$ and 
$q\set\mr^{}_{\Gamma\cO}\smcirc(\b1\>\ot\>\iota(\cO))$ from 
$\Gamma\cO\ot\Gamma\cO$ to $\Gamma\cO$ are equal.

As in the proof of~\ref{idpt1}, $\ml^{}_{\cO}=\mr^{}_{\<\<\cO}$, and so
commutativity of the subdiagrams of the following diagram is clear, whence 
$\iota(\cO)\smcirc p=\iota(\cO)\smcirc q$. As there is an isomorphism 
$\psi(\Gamma\cO)\colon \Gamma\cO\ot\Gamma\cO\iso\Gamma\Gamma\cO$ (see \eqref{PsiE}), 
\ref{coconds} implies that, indeed, $p=q$.
\end{proof}
\vskip-10pt
\[
\def\1{$\Gamma\cO\ot \Gamma\cO$}
\def\2{$\Gamma\cO\ot\cO$}
\def\3{$\Gamma\cO$}
\def\4{$\cO$}
\def\5{$\cO\ot \Gamma\cO$}
\def\6{$\cO\ot \cO$} 
 \bpic[xscale=1.25, yscale=.9]

   \node(11) at (0,-1){\1};
   \node(15) at (4,-1){\2}; 
   
   \node(23) at (4,-4){\3}; 
    \node(24) at (3.04, -3.3){\4};
   
   \node(41) at (0,-4){\5};  
   \node(45) at (2,-2.5){\6};

    \draw[->] (11)--(15) node[above=1, midway,scale=.75]{$\b1\ot\iota(\cO)$} ;  
    
    \draw[->] (23)--(24) node[above=-2.92, midway,scale=.75]{$\mkern38mu\iota(\cO)$} ;

    \draw[->] (41)--(45) node[above=-8, midway,scale=.75]{\kern-14pt\rotatebox{26}{$\b1\ot\iota(\cO)$}} ;

   \draw[->] (11)--(41) node[left=1, midway, scale=.75]{$\iota(\cO)\ot\b1$} ;  
  
    \draw[->] (15)--(45) node[above=-5, midway, scale=.75]{\rotatebox{26}{$\mkern-12mu\iota(\cO)\ot\b1$}} ;

   \draw[->] (15)--(23) node[right=1, midway, scale=.75]{$\mr^{}_{\Gamma\cO}\mkern40mu$} ;  
   \draw[->] (41)--(23) node[below=1, midway, scale=.75]{$\mkern55mu\ml^{}_{\Gamma\cO}$} ; 
   \draw[->] (45)--(24) node[above=-2.5, midway, scale=.75]{\kern34pt$\ml^{}_{\<\cO}\<\<=\<\mr^{}_{\<\<\cO}$} ;  
  \epic
\]

\begin{scor}\label{idem<->otcoreflect}
The natural functors taking\/ $A$ to\/ $\Gamma_{\<\<\!A}$\ \(respectively, $\Gamma\<$ to\/ $\Gamma\cO)$ are
quasi-inverse equivalences between the category of idempotent
$\D$-objects and that of\/ \mbox{$\otimes$-coreflectors} of\/~$\D$.\hfill$\square$
\end{scor}

\begin{sprop}\label{idpt2} 
For an object\/ $(A,\>\alpha)$ in\/ $\D/\cO$, 
the following are equivalent$\>.$\vs2

{\rm(i)} $(A,\>\alpha)$ is a $\D$-idempotent pair.\vs2

{\rm(ii)} For all\/ $F,\:G\in\D$ the composite map\vs1
\[
j^{}_{F\<,\>G}\colon\Hom_\D(A\ot F,\> A\ot G\>)
\xto[\under{1.2}{\textup{via}\;\alpha}]{}
\Hom_\D(A\ot F,\> \cO\ot G\>)
\xto[\under{1.2}{\<\textup{via}\;\ml^{}_G\!}]{\Iso}\Hom_\D(A\ot F, \>G\>)\\[-3pt]
\]
is an  isomorphism.\vs2

{\rm(iii)} The maps\/ $j^{}_{\<A,A}$ and\/ $j^{}_{\<A,\cO}$ 
in\/ \textup{(ii)} are injective, and\/ $j^{}_{\cO\<,\>A}$ is surjective.
\end{sprop}

\begin{proof} (i)${}\Leftrightarrow{}$(ii). By Lemma~\ref{coconds}, (ii) says that $(\Gamma_{\<\<\!A},\iota_\alpha)$ (see \eqref{equiv}) is  coreflecting, which, by ~\ref{idpt3}, just means that $(A,\alpha)$ is idempotent. 

(ii)${}\Rightarrow{}$(iii). Trivial.\vs1

(iii)${}\Rightarrow{}$(i). Suppose  (iii)
holds. In the (obviously) commutative diagram
$$
\begin{CD}
A\ot A @>\alpha\>\ot \b1_{\<\<A}>>\cO\ot A 
  @>\ms^{}_{\cO\<,A}>>A\ot\cO \\
@V\b1_{\<\<A}\>\ot\>\alpha VV @VV\b1_\cO\>\ot\>\alpha V @VV\alpha \>\ot\b1_\cO V\\
A\ot \cO @>>\under{1.2}{\alpha\>\ot\b1_\cO}>\cO\ot\cO 
@>>\under1{\ms^{}_{\cO\<,\cO}}>\cO\ot\cO
\end{CD}
$$
the map $\ms^{}_{\cO\<,\>\cO}$ 
is the identity of $\cO\ot\cO$ (see proof of~\Lref{idpt1}),\vs{-1}
so by injectivity of~$j^{}_{A,\cO}\>$,  $\b1_{\<\<A}\ot\alpha$ factors as 
$A\ot A\xto{\!\alpha\>\ot \b1_{\<\<A}\>}
\cO\ot A\xto{\under{.9}{\ms^{}_{\cO\<,A}}}A\ot\cO$,  whence
$$
\mr^{}_{\<\<A}\smcirc(\b1_{\<\<A}\ot\alpha)=\mr^{}_{\<A}\smcirc\>\ms^{}_{\cO\<\<,\mkern.5muA}\smcirc (\alpha\>\ot \b1_{\<\<A})\underset{\eqref{mon4}}= \ml^{}_{\<\<A}\smcirc (\alpha\>\ot \b1_{\<\<A}).
$$
It will suffice, therefore, to show that $\alpha\>\ot \b1_{\<\<A}$ is an isomorphism.

Surjectivity of~$j^{}_{\cO\<,\>A}$ entails the existence of a map
$\chi\colon A\ot \cO\to A\ot A$ such that 
$$
(\alpha\ot\b1_{\<\<A})\smcirc\chi=
\ms^{}_{A,\cO}\colon A\ot \cO\iso \cO\ot A,
$$
whence  $(\alpha\ot\b1_{\<\<A})\smcirc\chi\smcirc\>\ms^{}_{\cO\<\<,\mkern.5muA}=\b1_{\cO\ot A}$. 
Moreover, 
$$
(\alpha\ot\b1_{\<\<A})\smcirc\chi\smcirc\ms^{}_{\cO\<\<,\mkern.5muA}\smcirc(\alpha\ot\b1_{\<\<A})
=\alpha\ot\b1_{\<\<A} =(\alpha\ot\b1_{\<\<A}) \smcirc \b1_{\<\<A\ot A}\>,
$$
and since $j^{}_{A,A}\>$ is injective, therefore 
$
\chi\smcirc\ms^{}_{\cO\<\<,\mkern.5muA}\smcirc(\alpha\ot\b1_{\<\<A})= \b1_{\<\<A\ot A}. 
$

Thus $\alpha\ot\b1_{\<\<A}$ is indeed an isomorphism, with inverse $\chi\smcirc\ms^{}_{\cO\<\<,\mkern.5muA}$\kf.
\end{proof}

\begin{sprop}\label{Gam preserves idem} 
Let\/ $(\Gamma,\iota)$ be a $\ot$-coreflection of\/ $\D$. 
If the pair\/ $(B,\beta)$ is\/ $\D$-idempotent then so is\/ $(\Gamma \<B,\iota(\cO)\<\smcirc\<\Gamma\beta)$. 
\end{sprop}

\begin{proof} 
One checks, via \ref{Psi for coref}, that the following natural diagrams commute.
\[\mkern-3mu
\def\1{$\Gamma\cO\ot \Gamma \<B$}
\def\2{$\Gamma(\cO\ot B)$}
\def\3{$\Gamma(\cO\ot\Gamma\<B)$}
\def\4{$\Gamma\Gamma\<B$}
\def\5{$\cO\ot \Gamma\<B$}
\def\6{$\Gamma\<B$} 
 \bpic[xscale=1.2, yscale=.83]

   \node(11) at (2,-1){\1};
   \node(13) at (2,-3){\3};
   \node(15) at (2,-5){\2};

   \node(24) at (3.65,-3){\4};
   
   \node(41) at (5,-1){\5};  
   \node(45) at (5,-5){\6};

   \draw[->] (11)--(41) node[above, midway, scale=.75]{$\iota(\cO)\ot\b1$} ;  
     
   \draw[->] (15)--(45) node[above=1, midway, scale=.9]{$\Iso$} ;  
   
   \draw[->] (11)--(13)  node[right=1, midway,scale=.75]{$\simeq$}
                                   node[left, midway,scale=.75]{$\psi(\cO,\Gamma B)$} ;  
   \draw[->] (13)--(15) node[left=1, midway,scale=.75]{$\Gamma(\b1_{\cO}\ot \iota(B))\!$} ;  
 
    \draw[->] (41)--(45) node[left=1, midway,scale=.75]{$\simeq$} ;
   
   \draw[->] (13)--(24) node[above=.5, midway, scale=.75]{$\Iso$} ;
   \draw[->] (13)--(41) node[below, midway, scale=.75]{$\mkern60mu\iota(\cO\ot\Gamma\<B)$} ; 
   \draw[->] (24)--(45) node[above=-1, midway, scale=.75]{$\mkern25mu\simeq$} 
                                   node[below=-3.5, midway, scale=.75]{$\Gamma(\iota B)\!=\<\<\iota(\Gamma \<B) 
                                            \mkern120mu$} ;  
  \epic
\mkern25mu
\def\1{$\Gamma \<B\ot \Gamma\cO$}
\def\2{$\Gamma(B\ot\cO)$}
\def\3{$\Gamma(\Gamma\<B\ot\cO)$}
\def\4{$\Gamma\Gamma\<B$}
\def\5{$ \Gamma\<B\ot\cO$}
\def\6{$\Gamma\<B$} 
 \bpic[xscale=1.2, yscale=.83]

   \node(11) at (2,-1){\1};
   \node(13) at (2,-3){\3};
   \node(15) at (2,-5){\2};

   \node(24) at (3.65,-3){\4};
   
   \node(41) at (5,-1){\5};  
   \node(45) at (5,-5){\6};

   \draw[->] (11)--(41) node[above, midway, scale=.75]{$\b1\ot\iota(\cO)$} ;  
   \draw[->] (15)--(45) node[above=1, midway, scale=.9]{$\Iso$} ;  
   
   \draw[->] (11)--(13)  node[right=1, midway,scale=.75]{$\simeq$}
                                   node[left, midway,scale=.75]{$\psi'(\Gamma B,\cO)$} ;  
   \draw[->] (13)--(15) node[left, midway, scale=.75]{$\Gamma(\iota(B)\ot \b1_\cO)\!$} ;  
 
    \draw[->] (41)--(45) node[left=1, midway,scale=.75]{$\simeq$} ;
   
   \draw[->] (13)--(24) node[above=.5, midway, scale=.75]{$\Iso$} ;
   \draw[->] (13)--(41) node[below, midway, scale=.75]{$\mkern60mu\iota(\Gamma\<B\ot\cO)$} ; 
   \draw[->] (24)--(45) node[above=-1, midway, scale=.75]{$\mkern25mu\simeq$} 
                                   node[below=-3.5, midway, scale=.75]{$\Gamma(\iota B)\!=\<\<\iota(\Gamma \<B)
                                            \mkern120mu$} ;  
  \epic
\]

Then by \ref{f*idpt}, with $\xi\set\Gamma$, $u\set\iota(\cO)$ and 
$v(E,F)\set\Gamma(\b1_{\<E}\ot \iota(F))\smcirc\psi(E,\Gamma F)$
(see~\ref{Psi for coref}), Proposition~\ref{Gam preserves idem} follows from the fact---to be shown---that
\[
\Gamma(\iota(B)\ot\b1_\cO)\smcirc\psi'(\Gamma B,\cO)=v(B,\cO)\set\Gamma(\b1_B\ot\iota(\cO))\smcirc\psi(B,\Gamma\cO),
\] 
that is, the border of the following natural diagram, with $\zeta=\zeta(B)$ the composite isomorphism
\[
B\ot\Gamma\<\cO\xto{\!\psi'(B,\cO)}\Gamma(B\ot\cO)\xto{\;\Gamma(\mr^{}_{\!B}\<)\;}\Gamma\<B
\xto{\ \ml_{\Gamma\<B}^{-\<1}\ }\cO\ot \Gamma\<B,
\]
 commutes.\va5
\begin{small}
\[\kern-1pt
\def\1{$\Gamma\<B\<\ot\< \Gamma \cO$}
\def\2{$(\Gamma\cO\<\ot\< B)\<\ot\< \Gamma \cO$}
\def\3{$\Gamma\cO\<\ot\< (B\<\ot\< \Gamma \cO)$}
\def\4{$\Gamma(B\<\ot\<\Gamma\cO)$}
\def\5{$\Gamma(B\<\ot\<\cO)$}
\def\6{$\Gamma\cO\<\ot\< \Gamma\<B$}
\def\7{$(\Gamma\cO\<\ot\< \cO)\<\ot\< \Gamma\<B$}
\def\8{$\Gamma\cO\<\ot\< (\cO\<\ot\< \Gamma\<B)$}
\def\9{$\Gamma(\cO\<\ot\<\Gamma\<B)$}
\def\0{$\Gamma(\Gamma\<B\<\ot\<\cO)$}
 \bpic[xscale=2.75, yscale=1.75]

   \node(11) at (1,-1){\1};
   \node(12) at (2,-1){\2};
   \node(13) at (3.1,-1){\3};
   \node(14) at (4.12,-1){\4};
   \node(15) at (5.05,-1){\5}; 

   \node(21) at (1,-2){\6};
   \node(22) at (2,-2){\7};
   \node(23) at (3.1,-2){\8};
   \node(24) at (4.12,-2){\9};
   \node(25) at (5.05,-2){\0}; 

   \draw[->] (11)--(12) ;  
   \draw[->] (12)--(13) ;  
   \draw[->] (13)--(14) ;  
   \draw[->] (14)--(15) ;
  
   \draw[->] (21)--(22) node[above=1, midway, scale=.75]{$\via$} 
                                  node[below=1, midway, scale=.75]{$\psi(\cO)^{-\<1}$} ;  
   \draw[->] (22)--(23) ;  
   \draw[->] (23)--(24) ;  
   \draw[->] (24)--(25) node[above=1, midway, scale=.75]{$\via$} 
                                   node[below=1, midway, scale=.75]{$\psi'(\cO,B)$} ;
    
   \draw[->] (11)--(21)  ;
   \draw[->] (13)--(23)  node[left=1, midway, scale=.75]{$\via\> \zeta\!$} ;
   \draw[->] (14)--(24)  node[left=1, midway, scale=.75]{$\via\> \zeta\!$} ;                                
   \draw[->] (25)--(15)  ;  
 
  \node at (2,-1.52)[scale=.85]{\circled1} ;
  \node at (3.61,-1.52)[scale=.85]{\circled2} ;
  \node at (4.56,-1.52)[scale=.85]{\circled3} ;
 \epic
\]
\end{small}

\vskip-7pt
Subdiagram \circled2 clearly commutes. \va1

Subdiagram \circled1 expands as follows, with 
$\psi'(B)\set \Gamma(\mr^{}_B)\smcirc\psi'(B,\cO)$:

 \[
\def\1{$\Gamma\<B\<\ot\< \Gamma \cO$}
\def\2{$(\Gamma\cO\<\ot\< B)\<\ot\< \Gamma \cO$}
\def\3{$\Gamma\cO\<\ot\< (B\<\ot\< \Gamma \cO)$}
\def\6{$\Gamma\cO\<\ot\< \Gamma\<B$}
\def\7{$(\Gamma\cO\<\ot\< \cO)\<\ot\< \Gamma\<B$}
\def\8{$\Gamma\cO\<\ot\< (\cO\<\ot\< \Gamma\<B)$}
 \bpic[xscale=3.5, yscale=2.2]

   \node(11) at (1,-1){\1};
   \node(12) at (1.96,-1){\2};
   \node(13) at (3.1,-1){\3};
 
   \node(21) at (1,-2){\6};
   \node(22) at (1.96,-2){\7};
   \node(23) at (3.1,-2){\8};

   \draw[->] (11)--(12) node[above, midway, scale=.75]{$$} ;  
   \draw[->] (12)--(13) node[above=1, midway, scale=.9]{$$} ;

   \draw[->] (21)--(22) node[above, midway, scale=.75]{$\via$}
                                  node[below, midway, scale=.75]{$\psi(\cO)^{-\<1}$} ;  
   \draw[->] (22)--(23) node[above=1, midway, scale=.9]{$$} ;  
    
   \draw[->] (11)--(21)  node[left=1, midway, scale=.75]{$$} ;
   \draw[->] (13)--(23)  node[left=1, midway, scale=.75]{$\via\> \zeta\!$} ;
 
   \draw[->] (13)--(21) node[above=-5, midway, scale=.75]{\rotatebox{17}{$\via\>\psi'\<(B)$}};
   
  \node at (1.5,-1.4)[scale=.85]{\circled4} ;
  \node at (2.5,-1.7)[scale=.85]{\circled5} ;
  
 \epic
\]

Since $\ml^{}_\cO=\mr^{}_\cO$ (see proof of \ref{idpt1}), it follows from \ref{PsiF and PsiOF} that 
$\psi(\cO)=\mr^{}_{\Gamma\cO}$, so by \eqref{mon4}, the bottom row
composes to the map $\b1_{\Gamma\cO}\ot \ml^{-\<1}_{\Gamma \<B}\>$. The commutativity of \circled5 results then from the definition of $\zeta$. 

\goodbreak
Subdiagram \circled4 expands naturally as
\[
\def\1{$\Gamma\<B\ot \Gamma\cO$}
\def\2{$(\Gamma\cO\<\ot B)\ot \Gamma\cO$}
\def\3{$\Gamma\cO\<\ot(B\ot \Gamma\cO)$}
\def\4{$\Gamma\cO\<\ot\Gamma\<B$}
\def\5{$\Gamma\cO\<\ot(\Gamma\cO\<\ot B)$}
\def\6{$\Gamma\cO\<\ot\Gamma(B\ot\cO)$} 
\def\7{$\Gamma\cO\<\ot\Gamma(\cO\ot B)$}
 \bpic[xscale=3.3, yscale=1.75]

   \node(11) at (1,-1){\1};
   \node(12) at (2,-1){\2};
   \node(13) at (3,-1){\3}; 
   
   \node(21) at (1,-2){\4};
   \node(23) at (3,-2){\5};
   
   \node(31) at (1,-3){\6};  
   \node(33) at (3,-3){\7};

   \draw[->] (11)--(12) node[above, midway, scale=.75]{$$} ;  
   \draw[->] (12)--(13) node[above, midway, scale=.75]{$$} ; 

   \draw[<-] (21)--(23) node[above, midway,scale=.75]{$\via\> \psi (B)$} ;  
 
   \draw[<-] (31)--(33) node[above, midway, scale=.75]{$$} ;  
 
   \draw[->] (11)--(21)  ;  
   \draw[<-] (21)--(31) ;
   
   \draw[->] (13)--(23)  node[right=1, midway,scale=.75]{$$}
                                   node[left, midway,scale=.75]{$$} ;  
   \draw[->] (23)--(33) node[right=1, midway,scale=.75]{$\via\> \psi(\cO\<,\<B)$} ;  
   
   \draw[->] (12)--(23) node[above=.5, midway, scale=.75]{$$} ;
   \draw[->] (33)--(1.2,-2.17) node[below, midway, scale=.75]{$$} ; 
   
  \node at (1.95,-1.45)[scale=.85]{\circled4$_1$} ;
  \node at (2.75,-1.45)[scale=.85]{\circled4$_2$} ; 
  \node at (2.45,-2.4)[scale=.85]{\circled4$_3$} ;
  \node at (1.45,-2.66)[scale=.85]{\circled4$_4$} ;
 \epic
\]

The commutativity of subdiagram \circled4$_1$ is obvious. That of \circled4$_3$ is given by \ref{PsiF and PsiOF}, and of \circled4$_4$ by \eqref{mon4}. 

Finally, \eqref{mon3} gives that for any $A\in\D$, the border of the following 
natural diagram of isomorphisms commutes:
\[
\def\1{$(A\<\ot\<A)\ot B$}
\def\2{$A\<\ot(A\ot B)$}
\def\3{$(A\ot B)\ot\< A$}
\def\4{$A\<\ot(B\ot A)$}
 \bpic[xscale=3.2, yscale=1.9]

   \node(11) at (1,-1){\1} ;
   \node(12) at (2,-1){\2} ;
   \node(13) at (3,-1){\3} ;
 
   \node(21) at (1,-2){\1} ;
   \node(22) at (2,-2){\2} ;
   \node(23) at (3,-2){\4} ;

   \draw[->] (11)--(12) ;  
   \draw[->] (12)--(13) ;

   \draw[->] (21)--(22) ;  
   \draw[->] (22)--(23) ;  
    
   \draw[->] (11)--(21)  node[left=1, midway, scale=.75]{$\ms_{\<A,A}$} ;
   \draw[double distance=2] (12)--(22) ;
   \draw[->] (13)--(23) ;
 
   \draw[->] (13)--(22) node[right, midway, scale=.75]{$\mkern22mu\ms_{\<A\ot B\<,\>A}$};
   
  \node at (1.5,-1.5)[scale=.85]{\circled4$_5$} ;
  \node at (2.25,-1.45)[scale=.85]{\circled4$_6$} ;
  
 \epic
\]
Clearly \circled4$_6$ commutes. If $A=\Gamma\cO$---which, by \ref{idpt3}, is idempotent---then $\ms_{\<A,A}$ is the identity (see \ref{s=1}), so \circled4$_5$ commutes, whence so does the unlabeled subdiagram, which is just \circled4$_2$.

Thus \circled4, and hence \circled1, commutes.\va2

Subdiagram \circled3, with the leading ``$\>\>\Gamma\>\>$'' omitted,  expands naturally as
\[
\def\1{$B\ot\Gamma\cO$}
\def\2{$B\ot\cO$}
\def\3{$\Gamma B$}
\def\4{$ B$}
\def\5{$\cO\ot \Gamma B$}
\def\6{$\Gamma B\ot\cO$}
\def\7{$\Gamma(B\ot\cO)$}
 \bpic[xscale=2.15, yscale=1.25]

   \node(11) at (1,-1){\1};
   \node(14) at (4,-1){\2};
 
  \node(21) at (2,-2){\7} ;
  \node(22) at (2,-3){\3};
  \node(23) at (3,-3){\4};
   
  \node(31) at (1,-4){\5};
  \node(34) at (4,-4){\6};

   \draw[->] (11)--(14) node[above, midway, scale=.75]{$$} ;  
   
   \draw[->] (22)--(23) node[above=1, midway, scale=.9]{$$} ;

   \draw[->] (31)--(34) node[below=1, midway, scale=.75]{$\psi'\<(\cO\<,\<B)$};  
    
   \draw[->] (11)--(31)  node[left=1, midway, scale=.75]{$\zeta$} ;
   \draw[->] (34)--(14)  node[right=1, midway, scale=.75]{$$} ;
 
   \draw[->] (11)--(21) node[above=-4, midway, scale=.75]{$\mkern90mu \psi'\<(B\<,\<\cO)$} ;
   \draw[->] (21)--(22) ;
   \draw[->] (21)--(14) ;
   \draw[->] (23)--(14) ;
   \draw[->] (22)--(31) ;
   \draw[->] (22)--(34) ; 
   
  \node at (1.45,-2.55)[scale=.85]{\circled3$_1$} ;
  \node at (2.4,-1.45)[scale=.85]{\circled3$_2$} ;
  \node at (2.4,-3.65)[scale=.85]{\circled3$_3$} ;
 \epic
\]
\vskip-2pt
Subdiagram \circled3$_1$ commutes by the definition of $\zeta$, and \circled3$_2$, \circled3$_3$ by \ref{Psi for coref}.
The commutativity of the unlabeled subdiagrams is obvious.
Thus \circled3 commutes.

This completes the proof of \ref{Gam preserves idem}.
\end{proof}

Via \ref{idpt3}, two alternate formulations of \ref{Gam preserves idem} are: 
\begin{sprop}\label{idpt tensor}
\textup{(i)} Let $\mu\colon\cO\ot\cO\iso\cO$ be the map\/ $\ml^{}_\cO=\mr^{}_\cO$
\(\kern-.4pt see ~\textup{\ref{idpt1})}. If\/ $(A,\alpha)$ and\/ $(B,\beta)$ are\/ $\D$-idempotent pairs, then 
so is\/ $\big(A\ot B, \>\mu\smcirc(\alpha\ot\beta)\big)$.\va1

\textup{(ii)} If\/ $(\Gamma_{\!1},\iota_1)$ and\/ $(\Gamma_{\!2},\iota_2)$ are\/ $\ot$-coreflections of\/ $\D$ then so is\/
$(\Gamma_{\!2}\smcirc\Gamma_{\!1},\> \iota_2\smcirc\Gamma_{\!2}(\iota_1))$.\va2
\rightline{$\square$}
\end{sprop}

\begin{subcosa}\label{derived tensor}
\emph{To illustrate,} let $(X,\OX)$ be a ringed space. The category $\D(X)$ carries a well-known monoidal structure with $\otimes\set\Otimes\>$, $\cO\set\OX$,\va{-5} and $(\ma,\ml,\mr,\ms)$ the standard isomorphisms. (It suffices to check the axioms on the full subcategory spanned by the K-flat complexes.)%
\footnote{
An $\OX$-\kf complex $P$ is \emph{K-flat} if for every $\OX$-quasi-isomorphism 
$Q_1\to Q_2$ the resulting map $P\otimes Q_1\to P\otimes Q_2$ is also a quasi-isomorphism;
or equivalently, if for  every exact $\OX$-\kf complex~$Q$,
the complex~$P\otimes Q$ is also exact.  Every $\OX$-\kf complex $Q$ admits a  
\emph{K-flat resolution,} i.e., there exists a quasi-isomorphism $P\to
Q$ 
with $P$~K-flat \cite[p.\,139, 5.6]{Sp}. If $P$~is~K-flat then for any $\OX$-\kf complex~$Q$, 
the natural maps, with \smash{$\Otimes$} denoting \emph{left-derived tensor product,}\va{.6} are isomorphisms
$
\smash{P\Otimes Q\iso P\otimes Q, \quad Q\Otimes P\iso Q\otimes P,}
$
see \cite[p.\,147, 6.5]{Sp}, \cite[\S2.5]{Dercat}.
}
If $X$ is a scheme,\va{.5} then $\Dqc(X)\subset\D(X)$ contains $\OX$ and is closed under 
$\Otimes{}$,\va{-3} (see \cite[2.5.8.1]{Dercat}); thus it is a \mbox{\emph{monoidal subcategory}} of\/~$\D(X)$. 
\end{subcosa}

\vskip-2pt
Recall that if $\Phi$ is a finitary s.o.s.\ in~a scheme $X$, then $\R\vG{\Phi}\>\Dqc(X)\subset\Dqc(X)$ (Proposition~\ref{Dqc to itself}). 
\vskip-1pt

\enlargethispage*{3pt}
\begin{sprop}\label{coreflections}
\textup{(i)} Let\/ $X$ be a locally noetherian scheme, and\/ $\I$ an\/~$\OX$-base.
The pair $(\R\vG{\I},\iota^{}_\I\>)$ is a\/ $\ot$-coreflection of\/~$\Dqc(X)$ and of\/ $\D(X)$.\va1
 
\textup{(ii)} Let\/ $X$ be a scheme, and\/ $\Phi$ a finitary s.o.s.\ in~$X$.
The pair $(\R\vG{\Phi},\iota^{}_\Phi)$ is a\/ $\ot$-coreflection of\/~$\Dqc(X)$.\va1

\end{sprop}

\begin{proof}
(i).~As in \ref{coreflect exams}(c), $(\R\vG{\I},\iota^{}_\I\>)$ is a coreflection of~$\Dqc(X)$ and of $\D(X)$. 

Moreover, there is a natural functorial map
\begin{equation*}\label{tensor iso2}
\smash{\psi^{}_\I(E\<,F\>)\colon\R\vG{\I}E\OT{\<X}\<F\to \R\vG{\I}(E\<\<\OT{\<X} \<F)} \qquad(E,\>F\in\D(X)),
\tag{\ref{coreflections}.1}
\end{equation*}
defined as follows.

Note first that for any $\OX$-complexes $E,\>F$, one has
$\vG{\I}E\otimes\sX F=\vG{\I}(\vG{\I}E\otimes\sX F)$: by \ref{Gam and lim}, it's enough to show this
when $E$ and $F$ are $\OX$-modules, a simple task left to the reader. Hence
if $E$ is K-injective and $F$ is K-flat, and  $E\otimes\sX \<\<F\to G$ is a K-injective resolution, then
the image of the natural composite map 
\[
\vG{\I}E\otimes\sX \<\<F\to E\otimes\sX \<\<F\to G
\]
lies in~$\vG{\I}G$. Via standard considerations (e.g., \cite[p.\,69, 2.6.5]{Dercat}),  
the map $\psi(E,F)$ for \emph{arbitrary} $E,\>F\in\D(X)$ results.

From this description of $\psi^{}_\I(E\<,F\>)$ one gets a commutative diagram
\begin{equation*}
\begin{CD}
\R\vG{\I}\OX\<\<\Otimes_{\mkern-7mu X}\>F@>{\psi^{}_{\<\I}\<(\<\OX\<\<,\>F\>)}>>\R\vG{\I}F\\[-3pt]
 @V\iota^{}_{\<\I}\<(\OX)\>\ot\>\b1VV@VV\iota^{}_{\<\I}\<(F)V \\
\OX\<\<\Otimes_{\mkern-7mu X}\>F @>\Iso>\under{1.25}{\ml^{}_{\<\<F}}>  F
 \end{CD}
\end{equation*}
It remains to be shown that  $\psi^{}_\I(\OX\<,F\>)$ is an isomorphism (see \eqref{tensor iso}).\va1

Actually, $\psi^{}_\I(E\<,F\>)$ is an isomorphism for all $E$. This assertion is local, so assume  $X=\Spec(R)$  ($R$ a noetherian ring).
If $\I=\I_{\mkern-2.5mu J}$ for some quasi-coherent $\OX$-ideal~$J$, then  by \cite[(3.1.2)]{AJL}, 
$\psi^{}_\I(E\<,F\>)$ is indeed an isomorphism. 
Thus for arbitrary $\I$,  the natural  composite maps
\[
\vG{I}E\otimes\sX F\to \vG{I}(E\otimes\sX F)\to \vG{I}G\qquad(I\in\I)
\]
are all quasi-isomorphisms, and one can apply\va5 $\dirlm{\lift.25,{\lift.9,\halfsize{$I$},\in{\lift1,\I\>,}},\,}$ \kern-2pt to get a quasi-isomorphism 
$
\vG{\I}E\otimes\sX F\to \vG{\I}G, 
$
whose $\D(X)$-image $\psi^{}_\I(E\<,\>F\>)$ is an isomorphism, as desired.\va2

\smallskip
(ii). Proceed as in the proof of (i), with $\Phi$ in place of~$\I$ and 
\cite[p.\,25, (3.2.5)(i)]{AJL} in place of \cite[p.\,20, (3.1.2)]{AJL}. \va1

\begin{small}
Alternatively, assuming---as one may---that $X$ is affine, 
check, using \Pref{RGam and colim2}, that the $E\in\Dqc(X)$
for which $\psi^{}_\Phi(E)$ is an isomorphism span
a localizing subcategory $\D_{\sst \<\<\otimes}\subset\Dqc(X)$. 
Since $\OX\in \D_{\sst \<\<\otimes}$,  
\cite[p.\,222, Lemma 3.2]{Nmn} gives $\D_{\sst \<\<\otimes}=\Dqc(X)$. 
\end{small}
\end{proof}

From \ref{coreflections} and \ref{idpt3} one gets:
\begin{scor}\label{RGamidem}
Let\/ $X$ be a locally noetherian scheme, and\/ $\I$ an\/~$\OX$-base.
The pair\/ $(\R\vG{\I}\OX,\iota^{}_\I(\OX))$ is\/ $\Dqc(X)$-idempotent and\/ $\D(X)$-idempotent.\va1
 
More generally \(see\ \textup{\ref{2RGams}),} if\/ $\Phi$ is a finitary s.o.s.\ in~a scheme\/ $X,$
then the~pair\/ $(\R\vG{\Phi}\OX,\iota^{}_\Phi(\OX))$ is\/ $\Dqc(X)$-idempotent, hence\/ $\D(X)$-idempotent.\va1%
\hfill{$\square$}
\end{scor}

\subsection{Morphisms of idempotent pairs}\label{idemcat}
\stepcounter{thm} 
Notation remains as in \Sref{idemclosed}.

\begin{sprop}\label{one map}
 Let\/ $(A,\>\alpha)$ and\/ $(B,\beta)$ be\/ $\D$-idempotent pairs. 

There is \emph{at~most one}
morphism $\lambda\colon(B,\beta)\to(A,\alpha)$.
Such a\/ $\lambda$ exists if and only~if\/ 
$\ml_{\<B}\smcirc (\alpha\ot\<\b1_{\<\<B})\colon A\ot B\to  B$ is an isomorphism.
\end{sprop}

\begin{proof}  We'll need:

\begin{slem} \label{epi->iso}
Let $(C,\gamma)$ be a\/ $\D$-idempotent pair, and\/ $(B,\beta)\in\D/\cO$.
Suppose  the\/ $\D$-maps\/
$B\xto{q\>\>}C\xto{p\>\>}B$ satisfy $\beta p=\gamma$ and 
$p\>q=\b1_{\<\<B}$. Then\/ $p$ is an isomorphism.
\end{slem}
\begin{proof} 
The composition $\b1_{C\<}\otimes \gamma\colon
C\otimes C\xto{\b1_{C}\otimes\> p\>}C\otimes B
\xto{\b1_C\otimes \>\beta} C\otimes \cO$ is, by  \ref{idempotent pair}, an
isomorphism. So
$\>\b1_C\>\otimes\< p$ has both a left inverse 
and a right inverse,
and thus must~be an isomorphism.
 
In the following commutative diagram, the isomorphisms $j^{}_{\bullet,\bullet}$ are as in \ref{idpt2} 
(with $(A,\alpha)$ replaced by $(C,\gamma)$):
$$
\begin{CD}
\Hom(C,\> C\otimes C) @>\Iso>\under{1.4}{\<\textup{via\ }\b1_C\otimes\>p\,}> \Hom(C,\> C\otimes B)\\[-2pt]
@V j^{}_{\cO\<\<,\<C} V\simeq V @V\simeq V j^{}_{\cO\<\<,\<B}  V\\
\Hom(C,\>\cO\otimes C) @>>\under{1.4}{\<\textup{via\ }\b1_{\cO}\otimes \>p\,}>\Hom(C,\>\cO\otimes B)
\end{CD}
$$
Hence  $\phi\mapsto p\phi$ is an isomorphism from $\Hom(C,C)$ to $\Hom(C,B)$; 
and  since 
$$
pqp=\b1_{\<\<B}\mkern.5mu p=p=p\>\b1_C,
$$
therefore  $q\>p=\b1_C$,  so $p$~is indeed an isomorphism.
\end{proof}

Assuming that $\lambda$ exists,  and having in mind Remark~\ref{s=1} and Proposition~\ref{idpt tensor}, one finds that \Lref{epi->iso}, with
$(C,\gamma)\set(A\otimes B, \>\ml_\cO\smcirc(\alpha\otimes\beta))$, applies to 
$$
B\xto{\ml_{\<B}^{-\<1}}\cO\ot B\xto{\<\<(\beta\>\otimes\b1_{\<\<B}\<\<)^{-1}\!}
 B\otimes B\xto{\<\lambda\>\otimes\b1_{\<\<B} \> } A\otimes B
\xto{\<\alpha\>\otimes \b1_{\<\<B} \>} \cO\otimes B\xto{\;\ml_{\<B}\; }B,
$$
giving that  $\alpha\otimes\<\b1_{\<\<B}$ is an isomorphism, whence so is $\ml_{\<B}\smcirc (\alpha\ot\<\b1_{\<\<B})$;\va1
and conversely,  the composites
\[
B\xleftarrow{\;\ml_B\;}
\cO\ot B\xleftarrow{\>\alpha\>\otimes \b1_{\<\<B} \<}A\ot B\qquad\textup{and}\qquad
A\ot B\xto{\b1_{\<\<A}\ot\beta} A\ot\cO\xto{\;\mr^{}_{\!A}\;}A
\]
are $\D/\cO$-morphisms, so 
if  $\alpha\otimes\<\b1_{\<\<B}$ is an isomorphism then $\lambda$~exists.

Uniqueness of $\lambda$ results from the
following isomorphisms (the first and third induced by $\ml_{\<B}\smcirc(\alpha\otimes\<\b1_B)\colon A\ot B\iso B$), whose composition takes $\lambda$ to $\alpha\lambda=\beta\>$:
\[
\mkern-25mu
\Hom(B,\>A)
\!\iso\!
 \Hom(A\otimes B,\>A\otimes\cO)
\smash{\!\underset{j^{}_{\<\<B\<,\<\cO}}\iso\!}
\Hom(A\otimes B,\>\cO) 
\!\iso\!
\,\Hom(B,\>\cO).
\]
\vskip-16.4pt
\end{proof}

\vskip-12pt

\begin{srem}\label{I universal}
Recall from Remark~\ref{idpt and D/O} that
the $\D$-idempotent pairs span a strictly full subcategory $\mathbf I_\D$ of the slice category $\D/\cO$.

It follows from \Pref{one map} that $\mathbf I_{\D}$ \emph{is a preordered monoidal category} (see 
\Eref{exams-idem}(d)).
Indeed, $(\cO,\text{identity})$ is clearly a largest object; and,
maps of idempotent pairs $(C,\gamma)\to
(A,\alpha)$ and $(C,\gamma)\to(B,\beta)$ give rise naturally to
a composite $\D/\cO$-map 
$
(C,\gamma)\iso(C\ot C, \mu\smcirc(\gamma\ot\gamma))\to
(A\ot B, \mu\smcirc(\alpha\ot\beta)),
$ 
whence $(A\ot B,\mu\smcirc(\alpha\ot\beta))$ is, 
via the maps $\mr\smcirc(\b1_A\ot\beta)$ and $\ml\smcirc(\alpha\ot\b1_B)$,
 a greatest lower bound for $(A,\alpha)$ and $(B,\beta)$. So~the unit object and the product functor 
 in~$\mathbf I_\D$ are the same as those in $\D/\cO$, and the associated functorial maps
$\ma$, $\ml$, $\mr$ and $\ms$ are necessarily the same as those inherited from~$\D/\cO$. 
\end{srem}

\begin{srem}\label{auto idpt}
If $(A,\>\beta)$ and $(A,\alpha)$ are idempotent pairs then
there is a unique $\lambda\colon A\to A$
such that $\beta=\alpha\lambda$. This $\lambda$ is an
\emph{automorphism}, with inverse the unique $\lambda'\colon A\to A$
such that $\alpha=\beta\lambda'$. 
(By~\ref{one map}, $\alpha\lambda\lambda'=\alpha1^{}_{\<\<A}\implies 
\lambda\lambda'=1^{}_{\<\<A}\>$;
and similarly, $\lambda'\lambda=1^{}_{\<\<A}$.) 
Explicitly, using ~\ref{s=1} and the functoriality of $\ms$, one finds that 
\[
\alpha\ot\beta=\beta\ot\alpha\colon A\ot A\to \cO\ot\cO,
\] 
whence the following diagram commutes:
$$
\def\1{\cO\ot A}
\def\2{A\ot A}
\def\3{A\ot \cO}
\def\5{\cO\ot \cO}
\def\6{A}
\def\7{\cO}
\CD
\1 @<\alpha\>\ot\>1<< \2 @>{1\>\ot\>\beta}>> \3@=
\underset{\UnderElement{}{\downarrow}{7.6ex}{\!\<\mr}}\3\\[-4pt]
@V\ml VV @V 1\>\ot\>\alpha VV @VV\alpha\ot\>1 V\\
\6 @<<\under1\mr< \3 @>>\beta\>\ot\>1> \5\\[-2pt]
@| @. @VV\mr V\\
\6 @. \underset{\under1{\beta}}{\Rarrow{10em}} @.\7 @<<\under1\alpha< \,\6,
\endCD
$$
so that $\lambda$ is the composite isomorphism
$$
A\underset{\<\ml^{-\<1}}\iso \cO\ot A\,\underset{(\alpha\ot\b1)^{-\<1}}{\iso}\,A\ot A\,
\underset{\under{.9}{\b1\ot \beta}\ }{\iso}\,A\ot\cO\underset{\under{.6}\mr}\iso A.
$$

Conversely, it follows, e.g., from (i)$\,\Leftrightarrow\,$(ii) 
in~\ref{idem<->otcoreflect}, that if\/ $(A,\alpha)$ is idempotent and\/ $\lambda\colon B\iso A$ is a $\D$-isomorphism 
then $(B, \alpha\lambda)$ is idempotent. 

Thus,  \emph{the automorphism group of
any idempotent\/~$A\in\D$ acts faithfully and transitively on the set of\/ $\alpha\colon A\to\cO$ 
such that\/ $(A,\alpha)$ is an idempotent pair.}
\end{srem}

\begin{srem}\label{isos}
\emph{Idempotent pairs\/ $(B,\beta)$ and\/ $(A,\alpha)$ are isomorphic
$\Leftrightarrow$ there exists a\/ $\D$-isomorphism} $\lambda\colon B\iso A$. The implication
$\Rightarrow$ is trivial. Conversely, if such a $\lambda$ exists then 
$(B,\beta)$ and $(B,\alpha\lambda)$ are both idempotent 
whence, as in \ref{auto idpt}, there is an automorphism $\kappa\colon B\to B$ such that
$\beta=\alpha\lambda\kappa\>$; so $\lambda\kappa\colon (B,\beta)\to
(A,\alpha)$ is an isomorphism of idempotent pairs.
\end{srem}

\begin{sdef}\label{defcurly}
For idempotent $B$ and $A$,  $B\preccurlyeq A$ 
means there exist $\beta$ and~$\alpha$  and
a map---unique, by~\ref{one map}---of 
idempotent pairs $(B,\beta)\to (A,\>\alpha)$, 
a condition which is independent of the choice of~$\beta$ and $\alpha$.
\end{sdef}

By \Rref{isos},  $B$ is $\D$-isomorphic to $A$ $\!\iff\!$ 
$B\preccurlyeq A$ and $A\preccurlyeq B$. 

Of course $A\preccurlyeq A$, and
$C\preccurlyeq B$ together with $B\preccurlyeq A$ implies
$C\preccurlyeq A$. So we have a preordering on the idempotent $\D$-objects, such that $\cO$ is a largest object and, as in Remark~\ref{I universal}, $A\ot B$ is a greatest lower bound for $A$ and~$B$.\va{10}

\centerline{* * * * *\va{-3}}

\begin{sdef}\label{D_A}
For $A\in\D$, the category $\D_{\!A}\set\D_{\Gamma_{\lift.7,\<\<\!\sst A,}}\subset \D$ is the essential image of the
functor $\Gamma_{\<\<\!A}(-)\set A\otimes-$. 
\end{sdef}

\begin{slem}\label{D_idem} If\/ $\alpha\colon A\to \cO$ is a\/ $\D$-morphism such that\/ $\alpha\>\ot \b1\colon A\ot A\to \cO\ot A$ is an iso\-morphism, then\/ $E\in\D_{\!A}$ if and only if\va{-1} 
$$
\iota^{}_{\alpha}(E)\set\ml_E\smcirc (\alpha\ot\b1)\colon A\ot\< E\to E\\[-2pt]
$$ 
is an isomorphism.
\end{slem}

\begin{proof}
``If\kf\kf" is trivial, and ``only if\kf\kf" follows from the commutativity of the following diagram, with $F\in\D$
such that $A\ot F\cong E$ (see \eqref{mon5}):\va{-3}
\[
\def\1{$A\ot (A\ot F\>)$}
\def\2{$(A\ot A)\ot F\>)$}
\def\3{$\cO\ot (A\ot F\>)$}
\def\4{$(\cO\ot A)\ot F$}
\def\5{$A\ot F$}
 \bpic[xscale=3.75, yscale=1.3]

   \node(11) at (1,-1){\1} ;
   \node(12) at (1,-2){\2} ; 
   \node(13) at (3,-1){\5} ;
   
   \node(21) at (2,-1){\3} ; 
   \node(22) at (2,-2){\4} ;
   \node(23) at (3,-2){\5} ;
   
   \draw[->] (11)--(21) node[above=.5, midway, scale=.75]{$\alpha\<\ot\<\<(\b1\<\<\ot\<\b1)$} ;
   \draw[->] (21)--(13) node[above=.5, midway, scale=.75]{$\ml_{A\ot F}$} ;
   
   \draw[->] (12)--(22) node[below=1, midway, scale=.75]{$(\alpha\<\ot\<\<\b1)\<\<\ot\<\b1$}
                                  node[above=1, midway, scale=.75]{$\Iso$} ;  
   \draw[->] (22)--(23) node[below=1, midway, scale=.75]{$\ml_A\ot \b1$}
                                   node[above=1, midway, scale=.75]{$\Iso$} ;

    \draw[->] (11)--(12) node[left=1, midway,scale=.75]{$\ma$} 
                                   node[right,midway,scale=.75]{$\simeq$}  ;  
    
    \draw[->] (21)--(22) node[right=1, midway,scale=.75]{$\ma$} 
                                   node[left,midway,scale=.75]{$\simeq$}  ; 
      
    \draw[double distance=2] (13)--(23) ;                             
                                   
 \epic
\]
\vskip-27pt
\end{proof}
\vskip5pt
If\/ $\alpha\colon A\to\cO$ is as in \ref{D_idem} then\/
for $E\in\D_{\!A}$, one has  the functorial isomorphism $\ml_\alpha(E)\set \iota_\alpha(E)$;  and likewise,
the  functorial isomorphism \va{-1}
\[
\mr_\alpha(E)\set \mr^{}_{\!E}\smcirc(\b1\ot\alpha)\colon E\ot A\iso E.
\]

Clearly, $A\cong A\ot\cO\in\D_{\!A}$. For any $E\cong A\ot G\in\D_{\!A}$ and $F\in\D$, 
one has $E\ot F\in\D_{\!A}$ and $F\ot E\in\D_{\!A}$. 
In particular, $\D_{\!A}$ is closed under $\ot$.

\begin{slem}\label{monoidal D_A}
Let\/ $\alpha\colon A\to\cO$ be as in~\textup{\ref{D_idem}}, and\/ $\D_{\<*}\subset \D_{\!A}$ a full subcategory such that\/ $A\in\D_{\<*}$ and such that 
if\/ $E\<,F\in\D_{\<*}$ then\/  $E\ot F\in\D_{\<*}$. 
Then\/ $(\ot,A,\ma,\ml_\alpha,\mr_\alpha, \ms)$ is  a monoidal structure on~$\D_{\<*}\<$.
\end{slem}
  
\begin{proof}
For any $F$ and $B$ in $\D$, one has the diagram\va{-2}
\[
\def\1{$(F\ot A)\ot B$}
\def\2{$F\ot (A\ot B\>)$}
\def\3{$(F\ot \cO)\ot B$}
\def\4{$F\ot (\cO\ot B\>)$}
\def\5{$(A\ot F)\ot B$}
\def\6{$(\cO\ot F)\ot B$} 
\def\7{$F\ot B$} 
 \bpic[xscale=3.5, yscale=1.3]

   \node(11) at (1,-1){\1};
   \node(13) at (3,-1){\2}; 
   
   \node(22) at (2,-2){\3}; 
   \node(23) at (3,-2){\4};
    
   \node(31) at (1,-3){\5};  
   \node(32) at (2,-3){\6};
   \node(33) at (3,-3){\7};

    \draw[->] (11)--(13) node[above=1, midway,scale=.75]{$\ma$} ;  
    
    \draw[->] (22)--(23) node[above=1, midway,scale=.75]{$\ma$} ;

    \draw[->] (31)--(32) node[below=1, midway,scale=.75]{$(\alpha\<\ot\<\<\b1)\<\<\ot\<\<\b1$} ;
    \draw[->] (32)--(33) node[below=1, midway,scale=.75]{$\ml\<\ot\<\<\b1$} ;

   \draw[->] (11)--(31) node[left=1, midway, scale=.75]{$\ms\<\ot\<\<\b1$} ;  
  
   \draw[->] (22)--(32) node[left=1, midway, scale=.75]{$\ms\<\ot\<\<\b1$} ;  
   
   \draw[->] (13)--(23) node[right=1, midway, scale=.75]{$\b1\<\<\ot\<(\alpha\<\ot\<\<\b1)$} ;  
   \draw[->] (23)--(33) node[right=1, midway, scale=.75]{$\b1\<\<\ot\<\ml$} ;
      
   \draw[->] (11)--(22) node[above=-2.5, midway, scale=.75]{$\mkern105mu(\b1\<\ot\<\<\alpha)\<\ot\<\b1$} ;  
   
   \draw[->] (22)--(33) node[above=-2.5, midway, scale=.75]{$\mkern55mu\mr\<\ot\<\<\b1$} ;  

  \node at (2.5,-1.5)[scale=.9]{\circled1} ;
  \node at (1.3,-2.23)[scale=.9]{\circled2} ;
  \node at (2.2,-2.6)[scale=.9]{\circled4} ;
  \node at (2.83,-2.37)[scale=.9]{\circled3} ;

 \epic
\]
\vskip-4pt\noindent
The subdiagrams commute: \circled1 and \circled2 clearly, \circled3  by \eqref{mon1}, and  \circled4 by \eqref{mon4}. Therefore  \circled3 plus ~\circled1
give that \eqref{mon1} with~$(A,\cO,\mr,\ml)$ replaced by $(F,A,\mr^{}_\alpha,\ml^{}_\alpha)$ commutes;
and with $B\set\cO$, \circled4 plus~\circled2 give that \eqref{mon4} with~
$(A,\cO,\mr,\ml\>)$ replaced by~$(F,A,\mr^{}_\alpha,\ml^{}_\alpha)$ commutes.
The rest is obvious.
\end{proof}

\vskip-3pt
In \ref{monoidal D_A}, $\ml_\alpha$ and $\mr_\alpha$ depend on $\alpha$. However, 
for $\D/\cO$-isomorphic objects $\alpha\colon A\to \cO$ and $\alpha'\colon A'\to \cO$ as in \ref{D_idem}, 
it holds that $\D_{\!A}=\D_{\<\<A'}$, 
and  the monoidal structures on $\D_{\<*}$ induced by $\alpha$ and $\alpha'$ are \emph{equivalent,} where equivalence of
monoidal structures $(\ot,\cO,\ma,\ml,\mr,\ms)$ and $(\ot,\cO',\ma,\ml',\mr',\ms)$ 
means, with $\lambda\colon \cO'\to \cO$ the isomorphism 
$\mr^{\>\prime}_{\!\cO}\smcirc(\ml_{\cO'}^{-\<1})=\ml'_{\cO}\smcirc(\mr_{\!\cO'}^{-\<1})$ (see \eqref{mon4}),
that for all $E\in\mathbf{D_0}$ one has\va{-1}
\[
\ml'_{\<E}=\ml_{\<E}\smcirc(\lambda\ot \b1_{\<E})\colon \cO'\ot E\to E\quad\textup{and}
\quad \mr^{\>\prime}_{\!E}=\mr^{}_{\!E}\smcirc(\b1_{\<E}\ot \lambda)\colon E\ot \cO'\to E;
\]
in other words, the identity functor of $\D_{\<*}$ along with the identity map of 
\mbox{$E\otimes F$} $(E\<,F\in\D_{\<*})$ and the isomorphism $\lambda$ form an isomorphism of monoidal categories. (Details are left to the reader.)

\begin{sprop}\label{preccurly}
{\rm(i)} For\/ $\D$-idempotents\/ $B$ and\/ $A,$\va{-1} 
$$
B\preccurlyeq A \!\iff\! B\in\D_{\!A}\!\iff\! \D_{\<B}\subset\D_{\!A}.\\[-1pt]
$$
In particular, $\D_{\<B}=\D_{\!A}\!\iff\! B\cong A$.\va1

{\rm(ii)} Let \/ $(A,\>\alpha)$ be a $\D$-idempotent pair, and let\/ $\D_{\!A}$ have the monoidal structure given in\/ \textup{\ref{monoidal D_A}}. 
The map\/ $\Theta_{\<\<\alpha}$
that sends $(B,\lambda)\in\D/A$ to\/ $(B,\alpha\lambda)\in\D/\cO$~restricts to a bijection from the set of\/
$\D_{\!A}$-idempotent pairs
to the set of\/ $\D$-idempotent pairs\/~$(B,\beta)$ such that\/ 
$B\preccurlyeq A$. 

Thus  the\/ $\D_{\!A}$-idempotents are just the\/ $\D$-idempotents\/ $B$ such that $B \preccurlyeq A$.
\end{sprop}

\begin{proof}
(i) Left to the reader. (See Proposition~\ref{one map}.)

(ii) Let $(B,\lambda)$ be a $\D_{\!A}$-idempotent pair. 
Since $\ms_{\<B\<,\>B}\colon B\ot B\iso B\ot B$ is
the identity map (see \ref{s=1}), \eqref{mon4} ensures that the following
diagram commutes:
$$
\begin{CD}
B\ot B@>\Iso>\b1\ot\lambda >B\ot A @>\Iso>\b1\>\ot\>\alpha > B\ot \cO @>\Iso>\mr^{}_{\!B}> B\\[-3pt]
@|@V\lift.4,\ms^{}_{\<B\<\<,A\<\<}, V\simeq V @V\simeq V\lift.4,\ms^{}_{\<B\<\<,\cO}, V   @|\\
B\ot B@>\Iso>\lambda\ot\b1>A\ot B @>\Iso>\alpha\>\ot\> \b1> \cO\ot B@>\Iso>\ml^{}_{\<B}> B
\end{CD}
$$
\Pref{one map} gives that $\ml^{}_{\<B}\smcirc(\alpha\ot\b1)$ is an isomorphism, whence so is 
$\mr^{}_{\<B}\smcirc(\b1\ot\alpha)$, as are $\alpha\ot\<\b1$ and $\b1\ot\alpha$.
Hence $(B,\alpha \lambda)$ is $\D$-idempotent; and 
$\lambda$ is a map of $\D$-idempotent pairs $(B,\alpha\lambda)\to(A,\alpha)$, 
so that $B\preccurlyeq A$. Moreover, if $(B,\lambda')$ is a 
$\D_{\!A}$-idempotent pair such that $\alpha\lambda'=\alpha\lambda$,
then $\lambda$ and $\lambda'$ are maps from $(B,\alpha\lambda)$
to~$(A,\alpha)$, so by \Pref{one map},  $\lambda=\lambda'$. Thus $\Theta_{\<\<\alpha}$ acts injectively on $\D_{\!A}$-idempotent pairs.

Suppose $(B,\beta)$ is a $\D$-idempotent pair such that
$B\preccurlyeq A$, so that there exists a map of idempotent pairs $\lambda\colon(B,\beta)\to(A,\>\alpha)$. Then $\b1\ot\<\lambda\colon B\ot B\to B\ot A$ is an isomorphism, 
because its composition with the isomorphism
$1\ot \alpha \colon B\ot A\to B\ot \cO$ is the isomorphism $\b1\ot\beta$.
Similarly, $\lambda\ot\b1$ is an isomorphism. 

Again, $\ms_{\<B\<,\>B}$ is
the identity map, so the preceding commutative diagram gives 
$$
\mr_\alpha(B)\smcirc(\b1\ot\lambda)=\ml_\alpha(B)\smcirc(\lambda\otimes\b1),
$$ 
so that $(B,\lambda)$ is a $\D_{\!A}$-idempotent pair; and 
$
\Theta_{\<\<\alpha}(B,\lambda) = (B,\alpha\lambda)=(B,\beta).
$
Thus $\Theta_{\<\<\alpha}$~is surjective, as well as injective.

Verifying the last assertion is now straightforward.
\end{proof}

\begin{scor}\label{AotBidem}   Let\/ $(A,\alpha)$ and\/ $(B,\beta)$ be\/ $\D$-idempotent pairs.\va1

\textup{(i)} $\:\big(A\ot B, \mr_{\!A}\smcirc(\b1_{\<A}\ot \beta)\big)$ is\/ $\D_{\!A}$-idempotent.\va1

\textup{(ii)} $\<\big(A\ot B,\> \ml_B\<\smcirc\<(\alpha\ot \b1_{\<B})\big)$ is\/ $\D_{\<\<B}$-idempotent.

\end{scor}

\begin{proof}
By \ref{idpt tensor}, it holds that\va{-1}
\[
\big(A\ot B, \alpha\<\smcirc\<\mr_{\!A}\<\smcirc\<(\b1_{\<\<A}\ot\beta)\big)=
\big(A\ot B,\> \mr_\cO\<\smcirc\<(\alpha\ot \beta)\big)\\[-1pt]
\]
is $\D$-idempotent, and so (i)~results as in the latter part of the proof
of \ref{preccurly}(ii). The proof of (ii) is similar.
\end{proof}

\vskip-3pt
\centerline{* * * * *}

A \emph{closed category} is a  monoidal category
$\D$ (with product functor $\ot$) together with an \emph{internal hom} functor\va{-2}
$$
[-,-]\colon \D^{\textup{op}}\times \D\to \D\\[-2pt]
$$
and a trifunctorial isomorphism\va{-2}
\begin{equation}\label{closed}
\boh\colon\Hom_\D(E\ot F,G)\iso \Hom_\D(E,[F,G\>])
\qquad(E,F,G\in\D).\\[-2pt]
\end{equation}
(See, e.g., \cite[Definition (3.5.1)]{Dercat} and the references following
it.)

Elementary considerations show that the existence and functoriality of $\boh$ are equivalent to the existence for all $F$ and $G$ of 
an \emph{evaluation map,} functorial in $G$,
\begin{equation*}\label{closed'}\tag*{$\rm{(\theequation)}'$}
\textup{ev}=\textup{ev}^{}_{\<F\<,\>G}\colon [F,G\>]\ot F\to G
\end{equation*}
such that for all $E$, the map taking $\phi\colon E\to [F,G\>]$ to the map
$\textup{ev}\<\smcirc\< (\phi\>\ot\b1)\colon E\ot F\to G$ is an isomorphism
 $\Hom_\D(E,[F,G\>]\>)\iso \Hom_\D(E\ot F,G)$, and also such that for any map
 $F\to F'$ the following naturally induced  diagram commutes:
 \begin{equation*}\label{closed''}\tag*{$\rm{(\theequation)}''$}
\begin{CD}
[F'\<,G\>]\ot F @>>> [F\<,G\>]\ot F \\
@VVV @VV\textup{ev}V  \\
[F'\<,G\>]\ot F' @>>\textup{ev}> G
\end{CD}
 \end{equation*}

\smallskip
By definition, $\D$ has a unit object $\cO$ equipped with a 
functorial isomorphism
\begin{equation}\label{unit}
\mr^{}_{\<G}\colon G\ot\cO\iso G\qquad(G\in\D). 
\end{equation}
The natural composite isomorphism
$$
\Hom(F,G)\cong\Hom(F\ot\cO,G)\cong\Hom(F,[\cO,G\>])
\qquad(F,G\in\D)
$$
takes the identity map $\b1_G$ to a functorial isomorphism 
\begin{equation*}\label{unit'}
G\iso[\cO,G\>] \tag*{(\theequation)$'$}
\end{equation*}
corresponding under \ref{closed} to $\mr^{}_{\<G}\>$.\va2

For $(A,\alpha)\in\D/\cO$, one shows, via \ref{closed''} with $F\to
F'$ the map \mbox{$\alpha\colon A\to\cO$,} 
that the  map 
$G\ot A\to G\ot\cO\cong G$
induced by $\alpha$ factors as
\begin{equation}\label{gamlam}
G\ot\<A\,\underset{\textup{\ref{unit'}}}\iso \,[\cO,G\>]\ot\<A
\underset{\under{1.2}{\textup{via}\:\alpha}}\iso[A,G\>]\ot A
\xto[\under1{\textup{ev}}]{}G.
\end{equation}
\vskip-2pt\noindent
In particular,  the evaluation map is an isomorphism $[\cO,G\>]\ot \cO \iso G.$\va2

In \cite[pp.\,69--70]{DFS}, there are a number of formal relations 
which hold for any $\D$-coreflector $\BG$ that has a right adjoint
$\BL$---for example, if $(A,\alpha)$ is $\D$-idempotent,
the natural adjoint functors specified  objectwise by
$$
\BG G\set G\ot A,\qquad \BL G\set[A,G\>] \qquad(G\in\D).
$$
One such relation is the existence of an isomorphism
$\BG\iso\BG\BL$, which, for the preceding example, is just the
composition of the  first two maps
in \eqref{gamlam}.
\vspace{1.4pt}

\begin{srems}\label{internal and Hom}
(a) Internal hom is related to $\Hom_{\D}$ thus:
for $G\in\D$ set \begin{equation*}\label{Hzero}
\Hr^0 G\set\Hom_{\D}(\cO,G);\tag{\ref{internal and Hom}.1}
\end{equation*}
then there are natural isomorphisms
\begin{equation*}\label{H0[]}
\Hr^0[E,F\>]\iso \Hom_{\D}(\cO\ot E, F\>)\iso \Hom_{\D}(E, F\>)
\qquad(E,F\in\D).\tag{\ref{internal and Hom}.2}
\end{equation*}

(b) Let $\D$ be a closed category having an 
initial object~$A$, and 
$\alpha\colon A\to \cO$ the unique morphism. Then
$(A,\alpha)$ is idempotent. Indeed, \ref{closed} shows, for any
$F\in\D$, that $A\ot F$ is also an initial object, so that there are unique
maps $A\ot F\to A$ and $A\to A\ot F$, both isomorphisms. In particular, $\mr\smcirc(\b1\ot\alpha)$ and
$\ml\smcirc(\alpha\ot\b1)$ are equal isomorphisms from $A\ot A$ to $A$.

Clearly, $(A,\alpha)\preccurlyeq (B,\beta)$ for any idempotent $(B,\beta),$ i.e., $(A,\alpha)$ is initial in $\mathbf I_\D$.
\end{srems}

\begin{exams}\label{exams closed}
(a) For a ringed space $(X\<,\OX\<)$,\va{.5}  the derived category $\D(X)$ is
closed, with product {$\OT {\,S}\!$}\va1 (derived tensor product, see footnote in section~\ref{derived tensor}),\va{-4.5}
unit $\OX$, and $[E,F\>]\set\R\Homb_X(E,F\>)$. (For \eqref{closed} see e.g. \cite[p.\,147, 6.6]{Sp}, or in more detail,
\cite[\S2.6]{Dercat}.)
The maps $(\ma,\ml,\mr, \ms)$ are the obvious ones.

In particular, if $S$is  a ring (i.e., a ringed space~$(X,\OX)$ with $X$ a single point), then~$\D(S)$ is closed. \va2

(b) Suppose $\D$ is closed, let $(A,\alpha)$ be a $\D$-idempotent pair, and let
$\D_{\!A}\subset\D$ be the corresponding monoidal category (see \Dref{D_A} and \Lref{monoidal D_A}). The natural isomorphisms, with
$E,F,G\in\D_{\!A}$,
$$
\smash{
\Hom_\D(G\ot E, F\>)\underset{\eqref{closed}}\iso
\Hom_\D(G,[E,F\>])\underset{\eqref{idpt2}}\iso \Hom_\D(G,[E,F\>]\ot A)
}
$$
\vskip3pt\noindent
show that $\D_{\!A}$ is a closed category, whose internal hom is
$[E,F\>]_A\set[E,F\>]\ot A$.

\end{exams}

When $\D$ is closed one can expand on~\ref{idpt2}
and~\ref{coconds} in terms of $[-,-]$:

\begin{scor}\label{idpt4} For a closed category\/ $\D,$ and\/ $(A,\alpha)\in\D/\cO,$
\textup{\ref{idpt2}(ii)} holds if~and only if\kf\ 
for all\/ $F,G\in\D$ 
the following composite map is an isomorphism$\>:$
\begin{equation*}\label{gamiso''}
[A\ot\< F, \>A\ot\< G\>]\xto[\under1{\textup{via\;}\alpha}]{}[A\ot\< F, \>\cO\ot G\>]
\xto[\under1{\textup{via\;}\ml}\>]{\Iso}
[A\ot\< F, \>G\>].\tag{\ref{idpt4}.1}
\end{equation*}
Consequently \(see \eqref{idpt3}$),$ for any\/ 
$\ot$-coreflection $(\Gamma\<,\>\iota)$ 
the map induced by\/ $\iota(G)$ is an isomorphism
\begin{equation*}\label{gamiso*}
[\>\Gamma\<F\<,\>\Gamma G\>]\iso [\>\Gamma\<F\<,\> G\>].\tag{\ref{idpt4}.2}
\end{equation*}
\end{scor}

\begin{proof} That \eqref{gamiso''} is an isomorphism\va{.6} follows, upon application of the functor $\Hom(E,-)$\va{.6} 
with $E\in\D$ arbitrary, from the same for the map $j^{}_{\<E\ot A, \>F}$ in \ref{idpt2}(ii). 
The converse is given 
by application of the functor $\Hr^0$---see \eqref{H0[]}. 
\end{proof}

\subsection{Cohomology with supports: topological rings}\label{ALC}

\stepcounter{thm}

Prior considerations are rehearsed here in the context of topological rings.
This provides, among other things,  a formulation, encapsulated  in \ref{example} 
(appearing also in \cite[\S3.5]{Lectures}),
suited to subsequent developments, of some basic facts
about cohomology with supports. The underlying idea, which
will emerge fully only in the next section (see~\ref{top=idem})  is to establish a categorical
equivalence between ``decently topologized" noetherian rings  and noetherian rings $S$ furnished with idempotent $\D(S)$-pairs. 

This approach owes much to communications with Amnon Neeman.\vs2

\begin{subcosa}\label{toprings}
(Topologies on a commutative noetherian ring.)
A \emph{topological ring} $(S,\mU)$ is understood to be a \emph{noetherian}
ring $S$ with topology~$\mU$ such that addition and multiplication are continuous 
and  such that there is a basis $\mathfrak B$ of neighborhoods of~0 consisting
of ideals whose squares are open. (Any member of $\mathfrak B$ must itself be open, since an ideal $J$ that contains an open neighborhood $U$ of 0 also contains the open neighborhood $a+U$ of any $a\in J$.)  Such a topology on~$S$ will be called \emph{decent}.

For example, the \emph{preadic} $(S,\mU)$ are those having a $\mathfrak B$ consisting of all the powers of a single ideal \cite[p.\,172, (7.1.9)]{GD}. 

In a topological ring, any product $I_1I_2\dots I_n$ of open ideals  is open: induction reduces the proof to where $n=2$, and since $I_1\cap I_2$ contains some $J\in\mathfrak B$, therefore $I_1I_2$ contains the open ideal $U\set J^2\<$, and so, as above, $I_1I_2$ is open.

Since every open ideal contains a finite product of open prime ideals, therefore such products constitute a basis of neighborhoods of 0. 

Thus, for fixed $S$,  there is a bijection between decent  $\mU$ and sets $Y$ of prime ideals such that  for any prime ideals $p\subset p'\<$, $p\in Y\Rightarrow p'\in Y$, 
i.e., specialization-stable subsets of~$X\set$Spec$(S)$, or equivalently (see \S\ref{supports}),
between decent~$\mU$ and systems of supports (necessarily finitary) in $X\<$, or equivalently (see ~\ref{ex:supports}),  between decent~$\mU$ and $\OX$-bases. 

The specialization-stable subset of~$X$ corresponding to~$\mU$ consists of all $\mU$\kf-open prime ideals. 
The corresponding system of supports $\Phi_\mU$ consists of those closed subsets of~$X$ 
all of whose members are $\mU$\kf-open prime ideals, i.e., with~``$\;\tilde{}\;$" denoting sheafification, 
$
\Phi_\mU=\{\,Z(\tilde I\>)\mid I \textup{ is $\mU$\kf-open}\kf\}.
$
The corresponding $\OX$-base $\I_\mU$ is the set of sheafifications of 
$\mU$\kf-open $S$-ideals; and vice versa, for any $\OX$-base $\J$, the global section functor
takes the members of~$\>\J$ to the set of open ideals for a decent topology $\mU=\mU_{\>\J}$ such that
$\J=\I_\mU$.\va2

For a topological ring $(S,\>\mU)$, let $\iG\<\< =\iGp{\mU}$ be
the \emph{left-exact subfunctor of the identity functor on the category\/ $\A(S)$ of $S$-modules} such that for any $S$-module~$M\<$,
$$
\iG M=\{\,x\in M\mid\textup{for some open ideal } J, \>\>Jx=0\,\}.
$$ 
If $p$ is a prime $S$-ideal and $I_p$ is an injective hull of $S/\<p\>$---or of its fraction
field,  so that $I_p$ is an
$S_p$-module---then  $\iG{\<I_p}=0$ if $p$ is not open, and since every
element of
$I_p$ is annihilated by a power of~$p$,  $\iG{\<I_p}=I_p$ if $p$ is open. 
Thus $\iG$ determines the set of open primes, and hence determines the
topology $\mU\>$. \va{1.5}

With $\mfs M$ the sheafification of $M\<$, one has\va{-1}
\[\label{S to spec}
\iG M=\Gamma^{}_{\!\<\I_{\<\mU}}(X, \mfs M\>)=\Gamma^{}_{\!\<\Phi_{\<\mU}}(X, \mfs M\>)
=\Gamma(X,\vG{\Phi_{\<\mU}}\<\mfs M\>),\tag*{(\ref{toprings}.1)}
\]
\vskip-1pt\noindent
see~\ref{2Gammas} and \eqref{intersect sos}. Consequently, by~\ref{lims of Gams}, 
\emph{the~functor~$\iG\<\<$ commutes with small filtered colimits, hence with small direct sums.}

{\small More directly, if $x\in\dirlm{}M_\alpha$\va{-1}  is annihilated by an open ideal $J=(a_1,a_2,\dots,a_n)S$, then for some $\alpha$, $x$ is the natural image of 
an $x_\alpha\in M_\alpha$, and $a_ix_\alpha=0$ for all $i$, i.e.,~$Jx_\alpha=0$.\looseness=-1}
\vskip1pt
Moreover, $\iG$
\emph{preserves injectivity of\/ $S$-modules,}  since every injective $S$-module is
a direct sum of ones of the form $I_p$,  and any such direct sum is injective.\va2
 
\begin{small} 
In fact, $\mU\mapsto \iGp\mU\>$ is a bijection from decent
topologies on~$S$ t left-exact subfunctors~$\varGamma$ of the
identity functor on $\A(S)$ that commute with direct sums and preserve injectivity.
For, since $I_p$ is indecomposable, its injective submodule
$\varGamma I_p$ is  $I_p$ or 0\kf; and if~$p\subset p'$ then by
left-exactness, $\varGamma I_p\subset\varGamma I_{p'}$; so the set of
$p$ such that $\varGamma I_p=I_p$ is the set of open primes for a decent topology
$\mU\>$.  One checks then that $\varGamma=\iGp\mU\>$ by applying both functors
to representations of $S$-modules as kernels of maps between
injectives.\par
\end{small}
\end{subcosa}

During the rest of this section, $(S,\mU)$ will be a topological ring. By and large, the presented properties of\kf~$\iGp{}\set\iGp{\mU}$ and its derived functor $\R\iGp{}$  correspond, via sheafification, to previously discussed properties, over 
$\Spec(S)$, of $\vG{\I_{\<\mU}}$ and $\R\vG{\I_{\<\mU}}$. 

\begin{slem}\label{L:injective} Any injective\/ $S$-complex is $\Gamma'$-acyclic.
\end{slem}

\begin{proof}
The proof, via that of \mbox{\kf\cite[$(3.1.1)(2)'\Rightarrow(3.1.1)(2)$]{AJL}}), \emph{mutatis mutandis}, is like that of \ref{injective}(ii). 

For another proof---Koszul-free---see \cite[Lemma 3.5.1]{Lectures}.
\end{proof}

\begin{sprop}\label{oplus} 
Set\/ $\Hr^n_{\mU}\set \Hr^n\mkern-1.5mu\smcirc\R\iGp{\mU}$. Let\/ $A$ be a small filtered category, $\mathscr M$  a functor from~$A$ to the category of $S$-complexes, and\/ $n\in\mathbb Z$. Then the natural map is an isomorphism\va{-3}
\begin{equation*}
\dirlm{A} \<\<(\Hr^n_{\mU}\mkern-.5mu\smcirc\mathscr M)\iso 
\Hr^n_{\mU}\,\dirlm{A} \mathscr M\big.
\end{equation*}
\vskip2pt
\noindent
In particular, $\R\iGp{\mU}$ commutes with small direct sums in\/ $\D(R)$.
\end{sprop}   

\begin{proof}
As in the proof of \ref{RGam and colim2}, reduce to where $\mU$ has an open base consisting of powers of a single ideal
$\mathbf tS$, in which case the functor $\R\vG{\mU}$ is given by tensoring with the bounded flat complex 
$\Gamma(\Spec(S),\KK({\mathbf t}))$, rendering \ref{oplus} obvious.

Or, make use  of the existence of functorial K-injective resolutions \cite[Tag 079P]{Stacks}, commutativity 
of~$\Gamma'_\mU$ with small filtered colimits, preservation of quasi-isomorphisms by such colimits,  and  
($S$ being noetherian) injectivity of  filtered colimits of injective $S$-modules.
 \end{proof}

One has natural functorial maps\va{-2}
\[
\iota'_\mU\colon \R\iGp{\mU}\to {\mathbf1}
\]
\vskip-1pt\noindent
and, for decent topologies $\mU$ and $\mV$, with $\mU\>\cap\mV$ the topology whose open sets are those sets which are open for both $\mU$ and $\mV$ (the decent topology whose corresponding specialization-stable subset of~$X$ is the intersection of those of $\mU$ and $\mV$),\looseness=-1
\[
\gamma'_{\mU,\mV}\colon \R\iGp{\mU\>\>\cap\>\mV}\iso\R\iGp\mU\>\R\iGp{\mV}\>,
\]
an isomorphism because  $\iGp{\mU\>\cap\mV}=\iGp{\mU} \iGp{\mV}$, and, as above,
$\iGp{\mV}$ preserves injectivity, so Lemma~\ref{L:injective} can be applied.\va2

\begin{sprop} \label{Gamma' coreflects}
Let\/ $\mU$, $\mV$ be decent topologies on\/~$S$.
 The subtriangles in the following natural functorial diagram commute.\va{-1}
\[
\def\1{$\R\iGp{\mU\>\>\cap\>\mV}$}
\def\2{$\R\iGp\mV $}
\def\3{$\R\iGp\mU $}
\def\4{$\R\iGp\mU\R\iGp\mV $}
 \bpic[xscale=3.5,yscale=1.4]
  \node(11) at (1,-1){\1};
  \node(12) at (2,-1){\2};
  
  \node(21) at (1,-2){\3};
  \node(22) at (2,-2){\4};
 
   \draw[->] (11) -- (12) ;
   \draw[<-] (21) -- (22) ;
   
   \draw[->] (11) -- (21) ;   
   \draw[<-] (12) -- (22) ;
   
    \draw[->] (11) -- (22) node[below=-4.5, midway, scale=0.85]{$\gamma'_{\mU\<,\mV}\mkern40mu$}
                                     node[above=-3.5, midway, scale=0.85]{$\mkern40mu\simeq$} ;
 \epic
 \]
 \vskip-4pt\noindent
 In particular, $(\R\iGp\mU,\iota'_\mU)$ is coreflecting in\/ $\D(S)$.
\end{sprop}

\begin{proof}
Imitate the proof of \ref{symm}. (For the last assertion, set $\mV\set\mU\>$.)
\end{proof}

\begin{subcosa}\label{DmU}
Let $\A(S)$ be the category of small $S$-modules, and let
$\A_{\mU}(S)\subset\A(S)$ be the essential image of $\iG\set\iGp{\mU}$---the Serre subcategory 
(cf.~\ref{serre subcat}) whose objects are
the \mbox{\emph{$\mU$-torsion $S$-modules,}} that is, those $S$-modules~$M$ such that
$\iG M=M\<$, or equivalently, such that the localization $M_p$ vanishes for every non-open
prime \mbox{$S$-ideal} $p$.  
One can regard~$\iG\<$ as  being right-adjoint to the inclusion 
\mbox{$\A_{\mU}(S)\hookrightarrow\A(S)$.} 

At the derived level, let $\D_{\mU}(S)\subset\D(S)$ be the
full subcategory whose objects are those complexes $E$ whose homology\- modules
are all in~$\A_{\mU}(S)$, that is, whose localization $E_p$ is exact for
every non-open prime $S$-ideal $p$. Any complex in~$\A_\mU(S)$ is in $\D_{\mU}(S)$.
As in the remarks after \ref{serre subcat}, $\D_{\mathfrak U}(S)$ is a \emph{localizing subcategory} of~$\>\D(S)$.
\looseness=-1
\end{subcosa}

\begin{sprop}\label{idem}
An\/ $S$-complex~$E\>$ is in\/~$\D_{\mU}(S)$  if and only if the
natural map\/ $\iota'(E)\set\iota'_\mU(E\>)\colon\R\iG   E\to E\>$  is an isomorphism.
So\/ $\D_{\mathfrak U}(S)$ is the essential image of
the functor\/ $\R\iGp{}\colon\D(S)\to\D(S)$.
\end{sprop}

\begin{proof}  
Set $X\set\Spec(S)$. Let $\mfs$ be the sheafification functor, an equivalence of categories\- 
from $\A(S)$ to $\Aqc(X)$. 

One can assume that $E$ is injective. Since $S$ is noetherian,
the $\OX$-module~$\mfs E$ is injective, and the first assertion is given by the following logical equivalences:
\begin{align*}
E\in\D_{\mU}(S)
&\!\iff\!
\forall n\in\mathbb Z,\ \iGp{}\Hr^n\<\<E=\Hr^n\<\<E \\[3.3pt]
&\mkern-8mu\underset
{\lift1,\textup{\ref{S to spec}},}
{\!\iff\!}
\mkern-8mu\forall n\in\mathbb Z,\ \Gamma(X\<, \vG{\Phi_{\<\<\mU}}H^n\mfs E)=\Hr^n\<\<E\\
&\underset
{\lift.75,\textup{\ref{Cor2Gammas}},}
{\!\iff\!}
\forall n\in\mathbb Z,\ \vG{\Phi_{\<\<\mU}}H^n\mfs E=\mfs\Hr^n\<\<E\\[-1pt]
&\!\iff\!
\mfs E\in\D_{\Phi_{\<\<\mU}}\<(X)\\[1pt]
&\underset
{\lift.75,\textup{\ref{gamadj}},}
{\!\iff\!}
\vG{\Phi_{\<\<\mU}}\mfs E=\mfs E\\
&\underset
{\lift.75,\textup{\ref{Cor2Gammas}},}
{\!\iff\!}
\mfs\Gamma(X\<, \vG{\Phi_{\<\<\mU}}\mfs E)=\mfs E\\
&\mkern-8mu\underset
{\lift1,\textup{\ref{S to spec}},}
{\!\iff\!}
\mkern-8mu\mfs\iGp{}\<E=\mfs E\!\iff\! \iGp{}\<E=E.
\end{align*}
\vskip3pt
The last assertion results then from the last assertion in~\ref{Gamma' coreflects}.
\end{proof}
\begin{small}
Here is another argument for the first assertion in \ref{idem}.

If $\sigma_{\!E}^{}\colon E\to I_{\<E}$ is a K-injective resolution then
$\R\iGp{}\<E\cong\iGp{}\<I_{\<E}\in\A_\mU(S)\subset\D_{\mathfrak U}(S)$, and so
$\R\iGp{}\D(S)\subset\D_{\mathfrak U}(S)$. Thus if $\iota'(E\>)$ is an isomorphism then $E\in\D_\mU(S)$.

Conversely,  note via \ref{oplus} that  the $E\in\D_{\mathfrak U}(S)$ for
which $\iota'(E\>)$ is an isomorphism span a localizing
subcategory $\mathbf L\subset\D_{\mU}(S)$. 
Now 
\cite[p.\,526, Theorem 2.8]{chrom} says that any localizing subcategory 
$\mathbf L'\subset \D(S)$ is determined\- by the set of prime
ideals~$p$ such that $\mathbf L'\<$ contains the fraction field~$k(p)$ of~$S/\<p$.  
Since $k(p)$ is in $\D_{\mU}(S)\!\iff\! k(p)$~is
$\mathfrak U$-torsion $\!\iff\! p$ is open,   therefore
$\mathbf L=\D_{\mathfrak U}(S)$ if 
$\iota'(k(p))$ is an isomorphism for any open $p$, which indeed it is,
because
$k(p)$ admits a quasi-isomorphism into a bounded-below complex of
$\mathfrak U$-torsion $S$-injective modules (which~follows easily from the fact
that if an $\mathfrak U$-torsion module~$M$ is contained~in an injective
$S$-module~$J$ then $M$ is contained in the $\mathfrak U$-torsion injective
module~$\iGp{}\<\< J$).
\end{small}

\smallskip
Once again, set $\R\iGp{}\set\R\iGp{\mU}$ and $\iota'\set\iota'_\mU$.

\begin{scor}\label{coconds1}
For\/ $F\in\D_{\mU}(S)$ and\/ $G\in\D(S),$ 
$\iota'(G)\colon\R\iGp{} G\to G$ induces an isomorphism
$$
\Hom_{\D_{\mkern-1.5mu\mU}\mkern-1.5mu(S)}\<(F,\>\R\iG G)=
\Hom_{\D(S)}\<(F,\>\R\iG    G)
\iso\Hom_{\D(S)}\<(F,G).
$$
\end{scor}

\begin{proof}
In view of the last assertion in \ref{Gamma' coreflects}, this results from \ref{idem} and \ref{coconds}. 
\end{proof}

\begin{sprop}\label {esslimage}
The natural functor  is an equivalence of categories
\[
\smash{\D(\A_\mU(S))\xrightarrow[\under{3.3}{\,\approx\>\>}]{} \D_\mU(S),} 
\]
with quasi inverse\/ $\R\iGp{}|_{\D_{\mkern-1.5mu\mU\<(\<S)}}$.
\end{sprop}

\begin{proof}
Apply  \cite[p.\,49, 5.2.2]{DFS}
{\small (where the second ``let" should be  ``let $\bj$ be the").} 
\end{proof}

Let \smash{$\Otimes$}\vadjust{\kern.7pt} denote derived tensor product
in~$\D(S)$---defined via K-flat resolutions, see footnote in section~\ref{derived tensor}.

\begin{sprop}\label{Psi}
There is a unique bifunctorial\/ $\D(S)$-isomorphism
$$
\smash{\psi(E\<,F\>)\colon  \R\iG\< E \Otimes F \iso\R\iG (E\Otimes F\>)}
\qquad\bigl(E,\:F\in\D(S)\bigr)\\[-3pt]
$$
making the following diagram commute:
\[
\def\1{\raisebox{-13pt}{$\R\Gamma' E\OT{} F$}}
\def\2{\raisebox{-13pt}{$\R\Gamma'(E\OT{} F\>)$}}
\def\4{$E\OT{} F$}
 \bpic[xscale=3.1, yscale=1.6]

   \node(11) at (1,-1){\1};
   \node(12) at (2,-1){\2};
   
   \node(22) at (2,-2){\4};
     
   \draw[->] (11)--(12) node[above, midway,scale=.75]{$\psi(E,F)$} ;  
       
   \draw[->] (2,-1.2)--(22) node[below=-5, midway, scale=.75]{$\mkern95mu\iota'(E\OT{} F)$} ;  

   \draw[->] (1.2,-1.2)--(22) node[below=-5, midway, scale=.75]{$\iota'(E)\<\OT{}\< \b1_{\<F}\mkern110mu$} ;  
 
 \epic
\]

 \vskip-8pt 
Thus the coreflecting pair\/ $(\R\iGp{}\<\<, \iota')\ ($see \textup{\ref{Gamma' coreflects})} is $\ot$-coreflecting 
in\/ $\D(S)$.
\end{sprop}

\begin{proof} First, \smash{$\R\iG\< E\OT{}F \in\D_{\mU}(S)$}---just note that if  $E$ is K-injective, $F$ is~K-flat, and $p$ is a non-open prime $S$-ideal, then 
$(\iG\<E\ot_S \<F\>)_p\cong  (\iG\<E)_p\ot_{S_p} \<F_p=0$. Hence
the existence and uniqueness of the map $\psi(E\<,F\>)$ is given by \ref{coconds1}.

To show that $\psi(E,F\>)$ is an isomorphism one reduces, as in the proof of~\ref{coreflections}(i), to the preadic case (i.e., a basis of neighborhoods of 0 is given by the powers of a single ideal), and then applies \cite[(3.1.2)]{AJL}.

Alternatively, it's enough, by  \ref{Psi for coref}, to show that $\psi(S,F\>)$ is an isomorphism.
For variable $F\<$, $\psi(S,\>F\>)$ is compatible with triangles and direct sums; hence, and by \ref{oplus}, 
the $F$ for which $\psi(S,\>F\>)$ is an isomorphism span
a localizing subcategory $\mathbf F\subset\D(S)$. As $S\in \bf F$,   
\cite[p.\,222, Lemma 3.2]{Nmn} gives  $\bf F=\D(S)$. \va{-2}
\end{proof} 

\begin{scor}\label{example} 

The pair $(\R\iGp{}S, \iota'(S))$ is\/ $\D(S)$-idempotent, and
there is a unique  functorial isomorphism\va{-2}
$$
\smash{ \psi_{\<\<F}\colon \R\iGp{}\<S\OT{} \>\>F\iso \R\iGp{}\<F\qquad \big(F\in\D(S)\big)}\\[-2pt]
$$ 
making the following diagram commute:
\[
\CD\label{idptotimes}
\smash{\R\iGp{}\<S\OT{} F}@>\Iso>\lift1.2,\psi_{\!F},> \R\iGp{} \<F \\[-4pt]
@V \smash{\iota'_S\>\OOtimes\b1_{\<\<F\<}} VV @VV\under{1.4}{\iota'(F\>)} V\\[-2pt]
S\Otimes F@>\Iso>\lift1.2,\ml_F,>F\\[-2pt]
\endCD
\tag*{(\ref{example}.1)}
\]

\end{scor}

\begin{proof}
The first assertion follows from \ref{idpt3}, and the rest from \ref{Psi} with $E=S$. 
\end{proof}

\noindent\emph{Remark.} 
By \ref{D_A}, \ref{idem} and \ref{example}, $\D_\mU=\D_{\R\iGp{\mU}S}$.

\subsection{Idempotent pairs and topological rings}
\label{top+idem}
\stepcounter{thm}

In a later chapter, the discussion of Duality will involve, in particular,  the behavior of functors
vis-\`a-vis  compositions\va{-2} 
$$ (R,\mathfrak T\>)\xrightarrow{\varphi\,} (S,\mU)\xrightarrow{\psi\,}(T,\mV)
\xrightarrow{\chi\,}(U,\mathfrak W)\\[-1pt]
$$ 
of continuous topological-ring homomorphisms and vis-\`a-vis certain commutative ``base\kf-change" diagrams. 
For that discussion, the formal basics 
can be set up more efficiently, and more generally, in an expanded category
obtained by substituting idempotent pairs for topologies and
dropping noetherian hypotheses. This section  explicates the
expansion. 

\begin{subcosa}\label{equivalences}
For a noetherian ring $S$, a decent topology $\mU$ determines the isomorphism class of the idempotent pair\va{.6} $(A,\alpha)\set(\R\iGp{\mU}S, \iota'_{\mU}(S))$ 
(see~\ref{example}); and conversely, this $(A,\alpha)$ determines $\mU$, since an $S$-ideal $J\>$ is 
$\mU$-open\va{.4} if~and~only if  $S/J\in\D_{\mU\>}(S)$, that is,
by \ref{idem} and \ref{idptotimes}, if and only if  
\smash{$\alpha\OT S\b1\colon A\OT S S/\<J \to S\OT S S/\<J=S/\<J$} is an isomorphism.\par\vs{2}

\begin{small}
Alternatively, it holds that
an $S$-prime  ideal $p$ is $\mU$-open if and only if,  with $k(p)$  the fraction field of $S/p\>\>, $
\smash{$\alpha\OT S\b1\colon A\OT S k(p) \to S\OT S k(p)=k(p)$}\vs2 is an isomorphism. \par
\end{small}

So  the map\va{-3}
\begin{equation*}\label{top->idem}
\{\,\textup{decent topologies}\,\}\lra  \{\,\textup{isomorphism classes of $\D(S)$-idempotent pairs}\,\} 
\end{equation*} 
\vskip-3pt\noindent
that takes $\mU$ to the class of  $(\R\iGp{\mU}S, \iota'_{\mU}(S))$ has a left inverse.
In fact it is \emph{bijective}  \cite[p.\,65,  3.5.7\kf]{Lectures}, a result generalized to formal schemes below, 
in~\ref{A and Z}.
 It is also order-preserving: for decent topologies $\mU$, $\mV$ with
$\mU\subset\mV$ (as collections of open sets), and any $S$-module~$M,$  one has 
$\iGp{\mU}\<M\subset\iGp{\mV}M$, and hence
$\R\iGp{\mU}\<S\preccurlyeq \R\iGp{\mV}S$. More~generally,
from~\ref{Gamma' coreflects} and~\ref{Psi} one gets that for any decent topologies $\mU$, $\mV$,\va{-2}
\begin{equation*}\label{tensortop}
\smash{(\R\iGp{\mU}S\<\OT {\<S} \R\iGp{\mV}S,\> \iota'_\mU(S)\OT {\<S} \>\iota'_\mV(S))} 
\cong(\R\iGp{\mU\>\>\cap\>\>\mV}S, \>\>\iota'_{\mU\>\>\cap\>\>\mV}(S)).
\end{equation*}
\end{subcosa}
 
\begin{subcosa}\label{psi^*}
Next, a reformulation of continuity of maps of topological rings, in terms of
idempotent pairs. 

For a ring homomorphism $\psi\colon S\to T$, let  $\psi_*$ be the restriction-of-scalars functor from the category $\A(T)$ of $T$-modules to the category~$\A(S)$.
This functor is exact, so its derived functor, from $\D(T)$ to $\D(S)$, will also be denoted by ``$\,\psi_*$".

The extension-of-scalars functor
\mbox{$\,-\otimes_S T\,$} from $\A(S)$ to $\A(T)$, together with the counit map $\psi_*M\ot_ST\to M\  (M\in\A(T))$ given by scalar multiplication, is left-adjoint to $\psi_*$. Standard arguments
(cf.~e.g., \cite[\S(2.5.7)]{Dercat}) show that this functor has a left-derived functor $\psi^*\colon\D(S)\to\D(T)$, constructed objectwise by choosing for each $S$-complex~$E$ a K-flat resolution\vspace{.3pt} 
$\varsigma_{\mkern-1.5mu E}^{}\colon \<P_E\to E$, and setting\va2 $\psi^*\<\<E\set P_E\ot_S T$,  furnished with the $\D(T)$-map
$
\smash{\psi^*\<\<E=P_E\ot_S T
\xrightarrow[\under{3.9}{\varsigma_{\<\<E}^{}\>\otimes\>\>{\bf1}\,}]{}
E\ot_S T}.
$

There is a natural identification $\psi^*\<S=T$.

The functor $\psi^*$ is left-adjoint to $\psi_*$,
with counit map at any $T$-complex $F$ being the natural composite\va{-4}
\[
\psi^*\psi_*F=P_{\psi_*\<F}\ot_ST \xto{\varsigma\>\ot\b1 } \psi_*\<F\ot_S T\lra F,
\]
\vskip-2pt\noindent
cf.~\cite[3.1--3.2.2]{Dercat}, or see \cite[Tag 09T5]{Stacks}).

There is a unique bifunctorial isomorphism $\tau=\tau(E,E')\ (E,E'\in\D(S))$ such that the following otherwise natural diagram commutes.
\[
\CD
\psi^*(E\OT{S}E')\hbox to 0pt{$\mkern50mu\overset\Iso {\underset{\lift1,\tau,}{\Rarrow{6em}}}$}@. \psi^*\!E\OT{\>T}\psi^*\<\<E'\\[-3pt]
@VVV @VVV \\
(E\ot_{S}\<E')\<\ot_S\<T@>\Iso>> (E\ot_S\<T)\ot_{\>T} (E'\ot_S\<T)
\endCD
\]
This follows from \cite[(2.6.5)]{Dercat}, cf.~proof of \cite[(3.2.4(i)]{Dercat}.\va3

One checks that the following natural diagram commutes:
\begin{equation*}\label{^*unit}
\CD
\psi^*(E\OT{S}S)@>\Iso>\lift1.2,\tau,> \psi^*\<\<E\OT{\>T}\psi^*\!S\\[-2pt]
@V\simeq VV @| \\
\psi^*\<\<E@>\Iso>> \psi^*\<\<E\OT{\>T} T
\endCD
\tag*{(\ref{psi^*}.1)}
\end{equation*}

For an $S$-complex $E$ and a $T$-complex $F$, let 
$E\ot_\psi F$ be the $T$-complex 
\[
E\ot_\psi F\set (E\ot_ST)\ot_{\>T} F=E\ot_S F,\\[-3pt]
\]
and set \va{-3}
\[
\smash{E\>\OT\psi F\set \psi^*\<\<E\>\OT {\>T}F}.\\[2pt]
\]
As $P_E\otimes_ST $ is K-flat over $T$, there is a canonical $\D(T)$-map
\begin{equation*}
\smash{E\>\OT\psi F=(P_E\otimes_ST)\ot_{\>T} F\to E\otimes_\psi F, }
\end{equation*} 
making \smash{$\OT{\;\psi}$} \kern-2pta two\kf-variable\vspace{1.4pt} 
derived functor of $\otimes_\psi\>$.\vs1

In particular, since $S$ is K-flat as an $S$-complex vanishing in
all  nonzero degrees,  there is a canonical functorial
$\D(T)$-isomorphism
\begin{equation}\label{iso1} 
S\>\OT\psi F\iso   F\qquad\bigl(F\in\D(T)\bigr).
\tag*{(\ref{psi^*}.2)}
\end{equation}

There is a unique
bifunctorial \mbox{``projection"} isomorphism\va4
\begin{equation*}\label{isrho}
\smash[b]{
{\rho}\colon E^{\mathstrut}\>\Otimes_{\<\<\!S}\,\>\psi_* F\!\iso
\psi_*(\psi^*\<\<E\>\OT {\>T}F\>)=
\psi_*(E\>\Otimes_{\<\<\!\psi}\,F\>)}
\qquad\bigl(E\in\D(S),\;F\in\D(T)\bigr)
\tag*{(\ref{psi^*}.3)}
\end{equation*} 
whose composition\vs1 with the natural $\D(S)$-map
\smash{$\zeta\colon \psi_*(E\>\Otimes_{\<\<\!\psi}\,F\>)\to E\otimes_S \<F$} is the natural map\vs{.5}
\smash{$\beta\colon  E\>\Otimes_{\<\<\!S}\,\psi_*F\to E\otimes_S F$}---an isomorphism when $E$ is K-flat.\vs2
This~${\rho}$ can be identified with the natural $\D(S)$ isomorphism $P_E\ot_S F\iso (P_E\ot_S T)\ot_T F$.\vs2

One checks that  ${\rho}$ is an instance of the map $p_2$ in \cite[p.\,107, 3.4.6]{Dercat}.\va1

\end{subcosa}

Recall the definition of \kf $\D_{\!A}$ (see \ref{D_A}\kf), and the discussion of \kf$\D_\mU$ preceding \ref{idem}. 
Recall further Lemma~\ref{f*idpt}, which for $\xi\set\psi^*$ (as in \ref{psi^*}) and $u\set\b1$ gives that 
for any $\D(S)$-idempotent pair~$(A,\alpha)$,   the pair $(\psi^*\!A,\psi^*\<\alpha)$ is $\D(T)$-idempotent. 

\pagebreak[3]

\begin{sprop}\label{rho and psi1} 
Let\/ $\psi\colon S\to T$ be a ring homomorphism. Let\/ $(A,\alpha)$
be\/ $\D(S)$-idempotent and let $(B,\beta)$ be $\D(T)$-idempotent.
The following are equivalent.

\smallskip

{\rm (i)} $B\preccurlyeq \psi^*\<\<A,$ that is \textup{(\emph{see}~\ref{defcurly}),} there exists a\/ $\D(T)$-map\/ $\lambda:B\to\psi^*\<\<A,$ necessarily unique,  such
that\/ the following diagram commutes$\>:$
$$
\CD
B @>\lambda>> \psi^*\<\<A \\[-2pt]
@V\beta VV @VV\psi^*\alpha V \\[-1pt]
T@. \mkern-35mu\Equals{3em}\mkern35mu \psi^*\<S\\[3pt]
\endCD 
$$

{\rm (ii)} The map\/ \smash{$\alpha\<\Otimes_{\<\<\!\psi}\>\>\b1_{\<B} $} is an
isomorphism\/ 
$\smash{A\Otimes_{\<\<\!\psi}\,B \iso S\Otimes_{\!\psi}\,B}
\underset{\under{1.2}{\ref{iso1}}}=B$.\va3

{\rm (ii)$'$} The map\/ \smash{$\alpha\<\Otimes_{\<\<\!S}\>\>\b1_{\psi_{\<*}\<B} $} is an
isomorphism\/ 
$A\Otimes_{\<\<\!S}\,\psi_*B \iso
S\Otimes_{\<\<\!S}\,\psi_*B=\psi_*B$.\va1

{\rm (iii)} For\/ $E\in\D(S)$ the map\/ \smash{$\alpha\<\Otimes_{\<\<\!\psi}\>\b1_{\<E}$} induces a\/ $\D(T)$-isomorphism\/ 
$$
\smash{(A\OT S E)\OT \psi B \iso (S\OT S E)\OT \psi B =E\OT \psi B.}
$$

{\rm (iii)$'$} $\psi_*\D_{\<B}(T)\subset\D_{\!A}(S)$.\va4

If\/ $(S,\mU)$ and \/ $(T,\mV)$ are topological
rings,  and\/ $(A,\alpha)$
\(respectively $(B,\beta))$ is the idempotent pair\/ 
$(\R\iGp\mU S,\>\iota_\mU(S))$ \(respectively $(\R\iGp\mV T,\>\iota_\mV(T))),$
then each of the preceding conditions is equivalent to each of the following ones.

\smallskip 
{\rm (iv)} The map\/ $\psi$ is continuous.\vspace{1.4pt}

{\rm(v)} For\/ $G\in\D_\mV(T)$ the map\/ 
\smash{$\iota_\mU(S)\<\Otimes_{\<\<\!\psi}\>\b1_{G}$} is a\/ $\D(T)$-isomorphism 
$$
\R\iGp\mU S\>\Otimes_{\<\<\!\psi}\,G
\iso 
 S\>\Otimes_{\<\<\!\psi}\, G \underset{\ref{iso1}} = G.
$$

{\rm (v)$'$} $\psi_*\D_{\mV}(T)\subset\D_{\mU}(S)$.

\end{sprop}

\begin{proof}  
${\rm(i)}\<\Leftrightarrow\<\<{\rm(ii)}\<$. This results from ~\ref{one map}.\va2

${\rm(ii)}\Leftrightarrow{\rm(ii)}'$.
The map \smash{$\alpha\<\Otimes_{\<\<\!\psi}\>\>\b1_{\<B} $} is an isomorphism,\va1 that is, it  induces homology\- isomorphisms,
if and only if\va{.6} its image under the exact functor~$\psi_*$ does so; and by~\ref{isrho},  that image is 
(up to isomorphism) the map \smash{$\alpha\<\Otimes_{\<\<\!S}\>\>\b1_{\psi_{\<*}\<B} $} in~(ii)$'$. \vspace{1.4pt}

\vs2

${\rm(iii)}\Leftrightarrow{\rm(ii)}.$ Condition ${\rm (ii)}$ is the case $E=S$
of ${\rm (iii)}$. That ${\rm(ii)}\Rightarrow{\rm(iii)}$ results from 
the commutativity---elementary to check, e.g., by unwinding the relevant definitions and making use of \ref{^*unit})---of the natural diagram\va{-4}
\[
\def\1{\raisebox{-13pt}{$(A\OT S E)\OT \psi B$}}
\def\2{\raisebox{-13pt}{$(S\OT S E)\OT \psi B$}}
\def\3{\raisebox{-13pt}{$E\OT \psi B$}}
\def\4{\raisebox{-13pt}{$(E\OT S S)\OT \psi B$}}
\def\5{\raisebox{-13pt}{$E\OT \psi (S\OT \psi B)$}}
\def\6{\raisebox{-13pt}{$(E\OT S A)\OT \psi B$}}
\def\7{\raisebox{-13pt}{$E\OT \psi (A\OT \psi B)$}}
 \bpic[xscale=1.8,yscale=1.1]
  \node(11) at (1,-1){\1} ;
  \node(13) at (3,-1){\2} ;
  \node(15) at (5,-1){\3} ;
  
  \node(23) at (3,-2){\4} ;
  
  \node(34) at (4,-3){\5} ;
  
  \node(41) at (1,-4){\6 } ;
  \node(45) at (5,-4){\7} ;
  
   \draw[->] (11) -- (13) node[above=1,midway, scale=.75]{$\via\alpha$} ;
   \draw[double distance=2] (13)--(15) ;
   
   \draw[->] (41) -- (45) node[above=1,midway, scale=.85]{$\Iso$} ;
   
    \draw[->] (1,-1.25) -- (1,-3.73) node[left,midway, scale=.75]{$\simeq$} ;
    \draw[->] (3,-1.25) -- (3,-1.73) node[right,midway, scale=.75]{$\simeq$} ;
     \draw[->] (5,-3.73) -- (5,-1.27) node[right=1,midway, scale=.75]{(ii)} 
                                     node[left,midway, scale=.75]{$\simeq$} ;

   \draw[double distance=2] (3.6,-1.73) -- (4.68,-1.23) ;   
   \draw[->] (3.3,-2.35) -- (3.7,-2.73) ;
   \draw[->] (4.7,-3.73) -- (4.3,-3.3) node[below=-4,midway, scale=.75]{$\via\alpha\mkern50mu$} ;   
   \draw[->] (4.3,-2.73) -- (4.83,-1.28) ;
   \draw[->] (1.3,-3.73) -- (2.7,-2.35) node[below=-2,midway, scale=.75]{$\mkern40mu\via\alpha$} ;

 \epic
\] 
\vskip-5pt
${\rm(ii)}\Rightarrow{\rm(iii)}'\Rightarrow{}$(ii)$'$.  By~\ref{D_idem},  (iii)$'$ means that
for all \mbox{$G\in\D_{\<B}(T),$}\va{.6} 
the map \smash{$\alpha\OT{S} \b1\colon A\OT{S}\psi_*G\to S\OT{S}\psi_*G$}\va{1}
is an isomorphism, to prove which it suffices to consider those $G$ of the form \smash{$B\OT{\>T} F\ (F\in\D(T)$).}\va1  
For such~$G\<$, assuming (ii), apply the functor
\smash{$\psi_*\<\smcirc\>(\smash{-\OT{\>T}F})$} to the map \smash{$\alpha\<\Otimes_{\<\<\!\psi}\>\>\b1_{\<B}, $}
\va{1.4} and then use \ref{isrho} to get (iii)$'\<$. 
Conversely, (ii)$'$~is the case $G=B$ of (iii)$'$.\vspace{1.4pt}

(iii)$'{}\Rightarrow{}$(v)$'$. By the Remark after \ref{example},
$\D_{\mV}(T)=\D_{\<B}(T)$ and $\D_{\mU}(S)=\D_{\!A}(S)$.\vspace{1.4pt}

(iv)${}\Rightarrow{}$(v)$'$. For $G\in\D_\mV(T)$, $n\in\mathbb Z,$ and $a\in \Hr^nG$, the annihilator $\textup{ann}_T(a)$
is a $\mV$-open $T$-ideal,  so that if $\psi$ is continuous then $\textup{ann}_S(a)=\psi^{-1}\textup{ann}_T(a)$ is a $\mU$-open $S$-ideal. 
It follows that $\psi_*G\in\D_{\mU}(S)$.\vspace{1.4pt}

(v)${}'\Rightarrow{}$(iv). Let $q$ be a $\mV$-open prime $T$-ideal.
Every $x$ in the $T$-injective hull~$I_q$ of~$\>T/q\>$ is annihilated by a power of $q$, so
$\R\iGp\mV I_q\cong\iGp\mV I_q=I_q$.  Hence, \ref{idem} and~(v)$'$ give $\R\iGp\mU \psi_*I_q\cong \psi_*I_q$.\va{1.5}

Set $p\set \psi^{-1}q$. Then $\psi_*I_q$ has a natural $S_p$-module structure. So tensoring an
$S$-injective resolution of $\psi_*I_q$ by~$S_p$ produces an injective resolution~$J$ of $\psi_*I_q$ 
such that multiplication by any element in~$S\setminus p$ is an isomorphism of~$J$, so that if
$p$ is not $\mU$-open, then  
$0=\iGp\mU J\cong\R\iGp\mU \psi_*I_q\cong \psi_*I_q,$
which is absurd;  thus $p$~must be $\mU$-open. Since an ideal in a topological ring is open if and only if it contains 
an intersection of open prime ideals, it follows that $\psi^{-1}$ takes open $T$-ideals to open
$S$-ideals, whence $\psi$ is continuous. \vspace{1.4pt}

${\rm(v)}\Leftrightarrow{\rm(v)}'$. By~\ref{isrho}, application of $\psi_*$ to \smash{$\iota_\mU(S)\<\Otimes_{\<\<\!\psi}\>\b1_{G}$} produces the map  
$$
\iota_\mU(S)\<\Otimes_{\<\<\!S}\>\b1_{\psi_{\<*}\<G}\colon\R\iGp\mU S\>\Otimes_{\<\<\!S}\,\psi_*G
\lra 
 S\>\Otimes_{\<\<\!S}\, \psi_*G  = \psi_*G,
$$
\vskip-4pt\noindent
so that the map \smash{$\iota_\mU(S)\<\Otimes_{\<\<\!\psi}\>\b1_{\<G}$} is an isomorphism iff\va1 so is 
\smash{$\iota_\mU(S)\<\Otimes_{\<\<\!S}\>\b1_{\psi_{\<*}\<G}$}
(see the above proof that ${\rm(ii)}\Leftrightarrow{\rm(ii)}'\>$), that is,\va{1.5} by~\ref{idem} and \ref{example},  
iff $\psi_*G\in\D_{\mU}(S)$.

${\rm(v)}\Leftrightarrow{\rm(ii)}$.  Using \ref{example} \va1 one gets (v)~for \smash{$G\cong B\OT{\>T}F$}
from~(ii) by applying the functor $-\OT{\>T}F$.
Conversely, (ii) is the  case $G=B$ of (v).\va{-5}
\end{proof}

\begin{sscholium}\label{top=idem}
Consider  the category $\mathcal T$ of triples 
$(S,A\<,\alpha)$ with $S$ a commutative ring and 
$(A,\alpha)$ a $\D(S)$-idempo\-tent pair,  morphisms $(S,A\<,\alpha)\to (T,B,\beta)$  being ring
homomorphisms \mbox{$\psi\colon S\to T$} satisfying the equivalent conditions
(i), (ii), (ii)$'\<$, (iii) and~(iii)$'$ in \ref{rho and psi1}.  The functor that takes  $(S,\mU)$ to 
$(S,\R\iGp{\mU}S, \iota'_{\mU}(S))$, and homomorphisms to themselves, 
embeds the category\- of continuous homomorphisms of
topological rings \emph{fully faithfully} into~$\mathcal T\<$, with essential image  the full subcategory spanned 
by all $(S,A,\alpha)$ with $S$ noetherian (see \ref{equivalences}, \ref{rho and psi1}).
\end{sscholium}

\begin{subcosa}
Let $(S,\mU)$ be a topological ring, and $\psi\colon S\to T$ a homomorphism of noetherian rings. Let  $\mU\> T$
be the (decent) topology on\/ $T$ for which a basis of neighborhoods of\/~$0$ is the family of ideals\/ 
$\{JT\mid J\textup{ an $\mU$-open $S$-ideal}\>\}$. Then
\[
(\psi^*\R\iGp\mU S,\psi^*\<\iota'_{\mU}(S))\cong(\R\iGp{\mU\>T} T, \iota'_{\mU\>T}(T)).
\]

In view of the bijective order-preserving map from $T$-topologies to isomorphism classes of
$\D(T)$-idempotent pairs (see \ref{equivalences}),
this results from the equivalence of (i) and (iv) in \Pref{rho and psi1} and the fact
that $\mU\> T$ is the strongest among the $T$-topologies $\mV$  that make $\psi$
continuous.\va2

\begin{small}
Alternatively, if for any  $S$-ideal~$J\<$, $\mU_{\<\<J}^{}$ is the topology\va2
with the powers of $J\>$ as a basis of neighborhoods
of 0, then $\iGp\mU=\dirlm{}_{\lift-.02,\mkern-7mu J\>\textup{open},}\iGp{\mU_{\<\<J}^{}},$
allowing one to assume
$\mU=\mU_{\<\<J}^{}$, in which case one can use the representation of $\R\iGp{\mU}S$
by a Koszul complex\dots
(see proof of \ref{injective}). \par
\end{small}

\vskip3pt
It follows, under the assumptions preceding 
\ref{rho and psi1}(iv), that \emph{the~map\/~$\lambda$ in \textup{\ref{rho and psi1}(i)} 
is an isomorphism if and only if the topology\/  $\mV$ equals} $\mU T$. 
For,  \ref{rho and psi1}(i)$\Leftrightarrow$\ref{rho and psi1}(iv) shows that  the (continuous) 
identity map \mbox{$(T,\>\mU T)\to(T,\mV)$} has a continuous inverse if and only if 
$\R\iGp{\mV}T\preccurlyeq\R\iGp{\mU T}T\preccurlyeq \R\iGp{\mV}T$, that is, if and only if
$\lambda$ is an isomorphism (see line following \Dref{defcurly}).

In terms of prime ideals,  $\mV=\mU T$  signifies that
 a prime\/ $T\<\<$-ideal $p$ is $\mV$-open if and only if  
$\psi^{-1}(p)$ is $\mU$-open.

\end{subcosa}

\begin{subcosa}
Given morphisms $\varphi\colon(R,D,\delta)\to(S,A,\alpha)$ and
$\mu\colon(R,D,\delta)\to(U,B,\beta)$ (as in \ref{top=idem}), one checks, with 
$V\set S\otimes_R U\<$, and $\nu\colon S\to V\<$,  
$\xi\colon U\to V$ the canonical maps, that
$
\smash{(V,\>\nu^*\!A\>\Otimes_{\!\<\<V}\,\xi^*\<\<B,\> 
\nu^*\<\alpha\>\Otimes_{\!\<\<V}\,\xi^*\<\beta)}
$
is, together with\va{.5} 
$\nu$ and $\xi$, a fibered direct sum of $\varphi$ and $\mu$.
It follows (or can be shown directly) that if $S$, $U$ and $V$ are noetherian, and $(A,\alpha)$, $(B,\beta)$ correspond to the
$S$-topology~$\mU$ and the $T$-topology~$\mathfrak T$ respectively, then \
\smash{$(\nu^*\!\<A\>\Otimes_{\!\<\<V}\,\xi^*\<B,\> 
\nu^*\<\<\alpha\>\Otimes_{\!\<\<V}\,\xi^*\<\<\beta)$} corresponds to the
tensor-product topology $\mU V\cap \mathfrak T V$ on~$V\<$.
\end{subcosa}

\bigskip
\subsection{Cohomology with supports:\;formal schemes\va6}
\label{formal schemes0}

\stepcounter{thm}

In this section,\va{-6} $X$ will~be a \emph{noetherian formal scheme} \cite[p.\,407, (10.4.2)]{GD},
equipped with a specialization-stable subset~$Z$---or equivalently, with an s.o.s, see Section~\ref{supports}.
 An $\OX$-ideal will be called \emph{open} if it contains an ideal of definition of~$X\<$. A noetherian formal scheme~$X$ has an ideal of definition all of whose powers are ideals of definition, whence any power 
of an open $\OX$-ideal is open.   

The main results extend those in the preceding two sections, where $X$ is just an ordinary noetherian affine scheme.
For noetherian formal schemes, some basics  
on cohomology with supports are gone over in \ref{formal base}--\ref{RGam and colim4}\kf;  
the close relation (given in~ \cite{AJS2}) between 
specialization-stable subsets of~$X$ and  
idempotent pairs in the derived torsion category is reviewed in \ref{cor:GZ+tensor}\kf--\ref{x in Supp}\kf; 
and the interaction between derived torsion functors with 
maps of formal schemes\va1 is addressed in~\ref{liftidem}--\ref{rho and psi2}.\looseness=-1
  
The foundations of the theory of formal schemes, as presented in \cite[\S10]{GD}, are largely taken for granted. The notation and terminology to be used here can be chased down via the index in
\cite[p.\,125]{DFS}. Full~justification of statements to~be made requires, 
as indicated by references, numerous results which can be found in chapters 1--3
of~\cite{Dercat} (an exposition of standard material about unbounded
derived categories and the derived direct- and inverse\kf-image functors
associated to maps of ringed spaces) and in
\cite{DFS} (a~study of duality on formal schemes).  

\vskip4pt
\begin{subcosa}\label{formal base}

An $\OX$-base $\I$ is as in \ref{basedef}, with the constraint that members of\kf~$\I$ be~
\emph{coherent} and \emph{open}.  

Since $\OX$ is coherent \cite[p.\,428, (10.10.2.7)]{GD}), 
therefore $\OX\in\I$. 

An $\OX$-ideal belongs to such an~$\I$ if and only if so does its radical,
so if  $J$ is an ideal of definition (necessarily coherent, see \cite[p.\,429, (10.10.2.9)]{GD}) 
and $\>\>\overline{\!X\<}$~is the noetherian  scheme $(X,\OX\</J\>),$ and if
$\pi\colon\OX\to\cO_{\>\overline{\!X\<}}=\OX\</J$ is the canonical surjection, then there is a natural bijection\va{-1} 
\begin{equation*}\label{defvg0}
\I\mapsto \overline\I\set\{\,\cO_{\>\overline{\!X\<}\mkern2mu}\textup{-ideals } I\mid \pi^{-\<1}I\in\I\,\}
\tag*{(\ref{formal base}.1)}
\end{equation*}
\vskip-1pt\noindent
from the set of $\OX$-bases onto the set of $\cO_{\>\overline{\!X\<}}\,(\>=\!\OX\</J\>)$-bases.
Proposition~\ref{ex:supports} holds for open coherent $I$, giving an inclusion-preserving bijection 
from $\OX$-bases to specialization-stable subsets of~$X$ (see section~\ref{supports}). 
\end{subcosa}

\begin{subcosa}
Let $\A\set\A(X)$ be the abelian category of $\cO\sX$-modules. For any $\OX$-base $\I$,
one has the left-exact subfunctor\/ $\vG{\I}\colon\A\to\A$ of the identity functor, see ~\eqref{defvg}.

If $\I$ and $\J$ are $\OX$-bases then $\vG{\I}\vG{\J}=\vG{\I\cap\>\J}\>$, and so\va{.4} the functor~$\vG{\I}$ is idempotent, see~\eqref{local intersect sos}. Hence\va{.2} the essential image~$\A^{}_{\>\I}$ of~$\vG{\I}$ is~the full subcategory of $\A$ spanned by the $\OX$-modules $M$ such that $\vG{\I}M=M$. \va1
 
 \pagebreak[3]
 If $\J\subsetneq \I$ then $\A^{}_{\>\J}\subsetneq \A^{}_{\>\I}$. For, \va{-1}
 \[
M\in\A^{}_{\>\J}\mkern-6mu\implies\mkern-6mu \{M=\vG{\J}M\}\mkern-6mu  \implies\mkern-6mu  
\{\vG{\I}M\<=\<\vG{\I}\vG{\J}M\<=\<\vG{\I\cap\>\J}M\<=\<\vG{\J}M\<=\<M\}\mkern-6mu \implies\mkern-6mu M\<\in\<\A^{}_{\>\I}\>.
 \]
Furthermore,  if $I\in\I$  then\va{3}
 
\leftline{$\OX/I\in\A^{}_{\J}$\va7}
 
\leftline{$\mkern-5mu\!\iff\!\OX/I=\vG{\J}(\OX/I\>)=\dirlm{J\in\J^{\mathstrut}}\sHom_{\OX}(\OX/J,\OX/I\>)
=\dirlm{J\in\J^{\mathstrut}} (I\colon\! J)/I\va5
$}

\leftline{$\mkern-7mu\!\iff\!$ for some $J$, $1\in \Gamma(X,I\colon\! J)$ (because $\Gamma(X,-)$ respects $\dirlm{}\!)\!\iff\! I\in\J,$\va3} 

\noindent so\va1 if $I\notin\J$ then $\OX/I\in \A^{}_{\>\I}\!\setminus\<\A^{}_{\>\J}$. 
Thus $\J\ne\I\implies \A^{}_{\>\J}\ne\A^{}_{\>\I}$.\va1
 
 Also, for any $\I$, $\J$, it holds that $\vG{\I}\A^{}_{\>\J}=\A^{}_{\>\I\cap\>\J}=\A^{}_{\>\I}\cap\A^{}_{\>\J}$. For, as above, 
$M\in\A^{}_{\>\J}\Rightarrow\{\vG{\I}M=\vG{\I\cap\>\J}M\}$, so that
$\vG{\I}\A^{}_{\>\J}\subset\A^{}_{\I\cap\>\J}\subset \A^{}_{\>\I}\cap\A^{}_{\>\J};$\va{.5} 
and if $M\in  \A^{}_{\>\I}\cap\A^{}_{\>\J}$ then\va1 $M=\vG{\I}M\in\vG{\I}\A^{}_{\>\J}$.%
\footnote{More generally, for  endofunctors $\Gamma^{}_{\!\<\<1}$, $\Gamma^{}_{\!\<2}$ of a category $\D$,  
$\Gamma^{}_{\!\<\<1}\D_{\Gamma^{}_{\!\<2}}=\D_{\Gamma^{}_{\!\<\<1}\<\<\smcirc\Gamma^{}_{\!\<2}}=\D_{\Gamma^{}_{\!\<\<1}}\!\cap\D_{\Gamma^{}_{\!\<2}}.$}

Reasoning as in the proof of~\ref{serre subcat}(ii), one sees that $\A^{}_{\>\I}$  is a Serre---hence plump---subcategory of~$\A$. Also, as in the proof of \ref{Gam and lim}(ii) one sees that $\vG{\I}$ preserves small filtered colimits, so that $\A^{}_{\>\I}$ is closed under such colimits. 
\end{subcosa}

\begin{subcosa}\label{qct}
Let $\Aqc\subset\A$ (respectively~$\Avc\subset\A$) be\va{.5} the full subcategory spanned by the quasi-coherent
$\cO\sX$-modules (respectively~the $\cO_{\<\<X}$-modules\va1 which are small filtered\- colimits
of coherent ones---or equivalently, by \cite[p.\,33, 3.1.7]{DFS}, unions of coherent submodules). With $\I'$ the $\OX$-base comprising \emph{all}\va1 open coherent  $\OX$-ideals, and 
$\Aqct\set\Aqc\cap\A_{\>\I'}$, one has\va{.5} 
$\Aqct \subset\Avc\subset\Aqc$  \cite[p.\,32, 3.1.5 and p.\,48, 5.1.4]{DFS}. 
If $X$ is affine, then $\Avc=\Aqc$ \cite[p.\,32, 3.1.4]{DFS}. These are all \emph{plump} subcategories of $\A$, see \cite[p.\,34, 3.2.2 and p.\,48, 5.1.3]{DFS}.
It~is~clear that $\Avc$ is closed under small filtered colimits; and so is~$\Aqct$ \cite[p.\,48, 5.1.3]{DFS}.
As in the proof of \emph{loc.\,cit.,} 5.1.4, \emph{mutatis mutandis,}  one finds that for any $\OX$-base~$\I$, 
$\vG{\I}\Aqct\subset\Aqct$. Also,  $\vG{\I}\Avc\subset\Avc\>\>$: for if $(M_\alpha)_{\alpha\in A}$ is a directed system of coherent $\OX$-modules, then for all $\alpha\in\A$ and $I\in\I$, $\sHom_{\OX}(A/I,M_\alpha)$ is coherent \cite[p.\,33, 3.1.6(d)]{DFS}, and so\va{-1}
\[
\vG{\I}\>\>\dirlm{\alpha}\<\<M_\alpha
\cong\dirlm{\alpha}\<\<\vG{\I}M_\alpha
=\dirlm{\alpha}\dirlm{I\in\I^{\mathstrut}}\sHom_{\OX}(A/\<I,M_\alpha) \in\Avc\>.\\[2pt]
\]

If $X$ is an ordinary scheme then $\Aqct=\Avc=\Aqc$.
\end{subcosa}

Let $\Dqct\subset\Dvc\subset\Dqc$
be the full subcategories of~$\D$ spanned
by the complexes whose homology modules are all in~
$\Aqct$ (resp.~in~$\Avc$, resp.~in $\Aqc$). If $X$ is affine then $\Dvc=\Dqc$. If $X$ is an ordinary scheme then $\Dqct=\Dvc=\Dqc$.

Since~$\A_{\dots}$ is plump in $\A$, therefore  $\D_{\dots}$ is a
triangulated subcategory\- of~$\D$.\va{9
}

\centerline{* * * * *}

\vskip1pt

Let $Z\subset X$ be specialization-stable, and let $\Phi_{\<\<Z}$ be the set consisting of all subsets of $Z$ that are closed in~$X\<$. As noted in section~\ref{supports}, the map sending any such~$Z$ to $\Phi_{\<\<Z}$ is a bijection from the set of specialization-stable subsets of~$X$ to the set of systems of supports in~$X\<$.

Let \mbox{$\vG{\<Z}\set\vG{\Phi_{\!Z}}\colon\A\to\A$} be the functor of sections supported
in $Z$:\va{.6} for all $M\in\A$ and open $U\subset X\<,$\va{-1} 
$$
\bigl(\varGamma^{}_{\<\<\!Z}\> M\bigr)(U)\set 
\{\,\xi\in M(U)\mid \xi_x=0 \textup{ for all }x\in U\setminus Z\,\}.
$$
\vskip-1pt\noindent
The pair with components $\vG{\<Z}$ and its inclusion map into the identity functor is co\-reflecting in $\A$.

\begin{sprop}\label{Dqct to itself}
With preceding notation, it holds that\/ $\R\varGamma^{}_{\<\<\!Z}\Dqct\subset\Dqct$.
\end{sprop}

\begin{proof}
As $\Aqct$ is closed under \smash{$\dirlm{}\!\<$}, one can reduce via 
\eqref{HGam as lim.1}  to where $Z$ is closed in $X\<$, and then refer to 
the first sentence in \cite[\S2.1]{AJS2}).
\end{proof}

\begin{slem}\label{formal Gam preserves injective}
Let\/ $\I,$ $\J$ be\/ $\OX$-bases and\/ $E$ an  injective\/  $\OX$-complex.
Then\/ $\vG{\J}E$~is $\vG{\I}$-acyclic.
\end{slem}

\begin{proof}
To be proved is that the natural map $\vG{\I}\vG{\J}E\to\R\vG{\I}\vG{\J}E$ is an isomorphism, a local property, so that one can restrict to where  $X=\Spf(S)$ for some adic noetherian ring~$S$ with ideal of definition, say, $L$. 

The completion map $\kappa\colon X= \Spf(S)\to\Spec(S)=:X_0$ corresponds  to the \mbox{identity} map of $S$, considered as a continuous map of topological rings, with discrete source and $L$\kf-adically topologized target 
(see \cite[p.\,403, (10.2.1)]{GD}).
Topologically, $\kappa$~is the closed immersion $\Spec(S/\<L\>)\hookrightarrow \Spec(S)$.  

Accordingly, 
$\Phi_\I\set\{\,Z(I)\mid I\in\I\,\}$ can be regarded as an s.o.s.~in~$X_0$, with corresponding 
$\cO_{\<\<X_0}$-base~$\I_0$. Similarly, $\J$ determines an $\cO_{\<\<X_0}$-base~$\J_0$. \va1

To proceed, we'll need:

\begin{slem} \label{I and J}
With preceding notation,
\[
\I=\{\,I_0\OX\mid I_0\in\I_0\,\}.
\]
\end{slem}

\begin{proof}
For any $I_0\in\I_0$, $Z(I_0)\subset \Spec(S/L)$, so $\sqrt{I_0}$ contains the sheafification~$\widetilde L$. 
The locally-ringed-space map $\kappa$ being \emph{flat} 
\cite[p.\,185, (7.6.13), p.\,187, (7.6.18) and p.\,403, (10.1.5)]{GD}, one has\va{-1} 
\[
\kappa^*\<\<\sqrt{I_0}\cong \sqrt{I_0}\OX\supset\widetilde L\OX\cong\kappa^*\<\widetilde L\>.
\]
\vskip-1pt\noindent
Since $\widetilde L\OX\cong\kappa^*\<\widetilde L$ is an ideal of definition of $X$ (see \cite[p.\,420, (10.8.5), p.\,421, (10.8.8)(ii) 
and p.\,427, (10.10.1), second paragraph]{GD}), therefore the $\OX$-ideal $I_0\OX$ is~\emph{open}. Also, $I_0\OX\cong\kappa^*\<\<I_0$ is \emph{coherent} (see \cite[p.\,115, (5.3.14)]{GD}).  
Moreover,
$Z(I_0\OX)=\kappa^{-\<1}Z(I_0)\in\Phi_\I$, so $I_0\OX\in\I$.  
Thus~$\{\,I_0\OX\mid I_0\in\I_0\,\}\subset\I$.\va1

Let $I\in\I$, and let $I_0\subset\cO_{\<\<X_{\<0}\<}$ be the sheafication of 
$\Gamma(X,I)\subset\Gamma(X,\OX)=S$
(see \cite[p.\,402, (10.1.3)]{GD}).   Then $I_0\OX\cong \kappa^*\<\<I_0\cong I$, 
see \cite[p.\,420, (10.8.5), p.\,421, (10.8.8)(ii) with $i=\kappa$ and $\mathscr F=I_0$, and 
p.\,429, (10.10.2.9)(ii) with $M=\Gamma(X,I)$]{GD}.    Since $\kappa^*\<\<I_0$ is generated by its global sections, therefore so is $I$, whence $I=I_0\OX$.\va1

Since $\widetilde L\OX$ is an ideal of definition of $X\<$,
therefore $I\supset \widetilde L^n\OX$ for some $n>0$, so 
$\Gamma(X,I)\supset \Gamma(X,\widetilde L^n\OX)\supset L^n$, whence
$Z(I_0)\subset \Spec(S/L)$.
And since 
\[
    Z(I)=Z(I_0\OX)=\kappa^{-\<1}Z(I_0)\in\Phi_\I, 
\]
therefore $I_0\in\I_0$. Thus $\>\I\subset\{\,I_0\OX\mid I_0\in\I_0\,\}$.                                                                                                                                                                                                                                                                                                                                                                                                                                                                                                                                                                                                                                                                                                                                                                                                                                                                                                                                                                                                                                                                                                                                                                                                                                                                                                                                                                                                                                                                                                                                                                                                                                                                                                                                                                                                                                                                                                                                                                                                                                                                                                                                                                                                                                                                                                                                                                                                                                                                                                                                                                                                                                                                                                                                                                                                                                                                                                                                                                                                                                                                                                                                                                                                                                                                                                                                                                                                                                                                                                                                                                                                                                                                                                                                                                                                                                                                                                                                                                                                                                                                                                                                                                                                                                                                                                                                                                                                                                                                                                                                                                                                                                                                                                                                                                                                                                                                                                                                                                                                                                                                                                                                                                                                                                                                                                                                                                                                                                                                                                                                                                                                                                                                                                                                                                                                                                                                                                                                                                                                                                                                                                                                                                                                                                                                         
\end{proof}

Using the
commutativity of $\>\dirlm{}\!\<\<$\va{-1.2} with global sections over noetherian spaces 
\cite[p.\,641, Prop.\,6]{Kf},\va{1} plus \Lref{I and J}, plus the natural isomorphism $\kappa^*\!I_0\iso I_0\OX$\va{.3} 
one gets,\va{.5} for any \mbox{$\OX$-complex~$F\<$,}
 the natural composite isomorphism
 \[
 \CD
\begin{aligned}\label{zeta}
\xi_{\kappa\<,\I,F}\colon\kappa_*\vG{\I}F
&=\kappa_*\dirlm{\lift.01,{{\halfsize{$I$}\in\>}\lift.9,\halfsize{$\I$},},\,}\< \<
    \sHom_{\OX}\<\<(\OX/I,\>F\>)\\[1pt]
&\!\iso\dirlm{\lift.01,{{\halfsize{$I$}\in\>}\lift.9,\halfsize{$\I$},},\,}\< \<
    \kappa_*\sHom_{\OX}\<\<(\OX/I,\>F\>)\\[1pt]
&\!\iso\dirlm{\lift.01,{{\halfsize{$I_0$}\in\>}\lift.9,\halfsize{$\I_0$},},\,}\<\<
    \kappa_*\sHom_{\OX}\<\<(\kappa^*\<(\cO_{\<\<X_{\<0}}/\<\<I_0\>),\>F\>)\\[1pt]
&\!\iso\dirlm{\lift.01,{{\halfsize{$I_0$}\in\>}\lift.9,\halfsize{$\I_0$},},\,}\<\<
    \sHom_{\OX}\<\<(\cO_{\<\<X_{\<0}}/\<\<I_0,\>\kappa_*F\>)=\vG{\I_0}\<\kappa_*F.
\end{aligned}
\endCD
\tag*{(\ref{formal Gam preserves injective}.1)}
\]
\vskip6pt\noindent
Replacing $F$ by a K-injective resolution, one derives a functorial isomorphism
$
\bar \xi_{\kappa\<,\I,F}\colon\kappa_*\R\vG{\I}F\iso \R\vG{\I_0}\<\kappa_*F.
$
\va1

The map $\kappa$ being flat, the left adjoint~$\kappa^*$ of $\kappa_*$ is exact, and therefore $\kappa_*E$ is an injective $\cO_{\<\<X_0}$-complex, as is 
$\vG{\J_0}\<\kappa_*E$ (see \ref{Gam preserves injective}).
The obvious commutativity of the natural diagram\va{-2}
\[
\def\1{$\kappa_*\vG{\I}\vG{\J}E$}
\def\2{$\kappa_*\R\vG{\I}\vG{\J}E$}
\def\3{$\vG{\I_0}\<\kappa_*\vG{\J}E$}
\def\4{$\R\vG{\I_0}\<\kappa_*\vG{\J}E$}
\def\5{$\vG{\I_0}\<\vG{\J_0}\<\kappa_*E$}
\def\6{$\R\vG{\I_0}\<\vG{\J_0}\<\kappa_*E$}
 \bpic[xscale=3.25,yscale=1.5]
  \node(11) at (1,-1){\1};
  \node(12) at (2,-1){\3};
  \node(13) at (3,-1){\5};
  
  \node(21) at (1,-2){\2};
  \node(22) at (2,-2){\4};
  \node(23) at (3,-2){\6};
  
   \draw[->] (11) -- (12) node[above, midway, scale=.75]{$\Iso$}
                                    node[below, midway, scale=.75]{$\via\xi$}  ;
   \draw[->] (12) -- (13) node[above, midway, scale=.75]{$\Iso$}
                                    node[below, midway, scale=.75]{$\via\xi$}  ;
   
   \draw[->] (21) -- (22) node[above, midway, scale=.75]{$\Iso$}
                                    node[below, midway, scale=.75]{$\via\bar\xi$}  ;
   \draw[->] (22) -- (23) node[above, midway, scale=.75]{$\Iso$}
                                    node[below, midway, scale=.75]{$\via\bar\xi$}  ;
   
   \draw[->] (11) -- (21) node[left=1, midway, scale=.75]{$\kappa_*\zeta$} ;   
   \draw[->] (12) -- (22) ;
   \draw[->] (13) -- (23) node[left, midway, scale=.75]{\ref{injective}(ii)}
                                    node[right, midway, scale=.75]{$\simeq$}  ;
   
 \epic
 \]
 \vskip-4pt\noindent
shows then that $\kappa_*\zeta$ is an isomorphism, whence so is $\zeta$, because $\kappa$ is, topologically, a closed immersion. Thus \Lref{formal Gam preserves injective} holds.
\end{proof}

\begin{sprop}\label{derived intersect5}
Let\/ $\I$ and\/ $\J$ be\/ $\OX\<$-bases and $E$ an\/ $\OX$-complex. The natural map\va{-5} is an isomorphism
\[
\gamma^{}_{\I\<,\J}\colon\R\vG{\I\cap\>\J}E\iso\R\vG{\I}\R\vG{\J}E\>.
\]
such that the following natural diagram commutes:
\[
 \def\1{$\R\vG{\I\cap\>\J}E$}
\def\2{$\R\vG{\J}E$}
\def\3{$\R\vG{\I}E$}
\def\4{$\R\vG{\I}\R\vG{\J}E$}
 \bpic[xscale=3.5,yscale=1.5]
  \node(11) at (1,-1){\1};
  \node(12) at (2,-1){\2};
  
  \node(21) at (1,-2){\3};
  \node(22) at (2,-2){\4};
 
   \draw[->] (11) -- (12) ;
   \draw[<-] (21) -- (22) node[below=1, midway, scale=0.75]{$\R\vG{\I}(\iota^{}_{\J}E)$} ;
   
   \draw[->] (11) -- (21) ;   
   \draw[<-] (12) -- (22) node[right=1, midway, scale=0.75]{$\iota^{}_\I(\R\vG{\J}E)$} ;
   
    \draw[->] (11) -- (22) node[below=-2, midway, scale=0.75]
                          {$\gamma^{}_{\I\<,\J}\mkern45mu$} ;

 \epic
\]  
\end{sprop}

\begin{proof} By \ref{formal Gam preserves injective}, for $\gamma^{}_{\I\<,\J}$ to be an isomorphism
it suffices that $\vG{\I}\vG{\J}=\vG{\I\cap\>\J}\>$, which one can show by imitating the argument used to establish
\eqref{local intersect sos}. 

The commutativity can be shown by arguing  just as in the proof of~\ref{symm}.
\end{proof}

As in \ref{exams-mon}, derived tensor product
makes $\D\set\D(X)$  into a symmetric monoidal category, with unit object~$\OX$.

Deriving the inclusion $\vG{\I}\hookrightarrow\b1_\A$ (the identity functor of~$\A$) produces a functorial
map $\iota^{}_\I\colon \R\vG{\I}\to {\bf1_{\D}}$.  
 \va1  

\begin{scor}\label{Gamma_I}
Set \va{-2}
\begin{equation*}\label{Gama'}
(\R\varGamma'\!, \>\iota')\set(\R\vG{\I'},\iota^{}_{\I'}\<) 
\qquad(\I' \textup{ as in \ref{qct}}).
\tag*{(\ref{Gamma_I}.1)}
\end{equation*}
Then\/$:$

\textup{(i)} $\R\varGamma'\Dqc\subset\Dqct\mkern.5mu.$ 

\textup{(ii)} $\Dqct$ is the essential image of\/ $\R\varGamma'\colon\Dqct\to\D$. 
\end{scor}

\begin{proof}
By \cite[p.\,49, 5.2.1(a)]{DFS}, a complex $E\in\Dqc$ lies in~ $\Dqct$ 
if and only if $\iota'(E)\colon\R\varGamma'E\to E$ is an isomorphism. In particular, $\Dqct$ is contained in the essential image 
of $\R\varGamma'\colon\Dqct\to\D$. Moreover,
if $E\cong \R\varGamma'F\ (F\in\Dqc)$ then $E\in\Dqct$, since by \ref{derived intersect5},
\[
\R\varGamma'E\cong\R\varGamma'\R\varGamma'F\mkern-5mu\underset{\iota'(\R\vG{\I'}F)}\cong\R\varGamma'F\cong E.
\\[-2pt]
\]
Thus $\Dqct$ contains the essential image of $\R\varGamma'\colon\Dqc\to\D\>$; and \ref{Gamma_I} follows.
\end{proof}
 
\begin{sprop}\label{vGI tensor-coreflects}
For any\/ $\OX$-base\/ $\I,$ the pair $(\R\vG{\I},\iota^{}_\I\>)$ is a\/  $\otimes$-coreflection of\/~$\D;$ and so 
$(\R\vG{\I}\OX,\iota^{}_\I(\OX))$ is\/ $\D(X)$-idempotent.\looseness=-1
\end{sprop}

\begin{proof}
That $(\R\vG{\I},\iota^{}_\I\>)$ is coreflecting\va1 in $\D$ results from~\ref{derived intersect5}.
For compatibility with~$\OT{}$, argue as in the proof of~\ref{coreflections}(i), except for\va{-3} replacing 
$\Spec(R)$ there by an affine formal scheme $X\set\textup {Spf}\>(\<R)$ with $R$ an admissible noetherian ring, and taking~$J$ to be a coherent open ideal. The last assertion results from ~\ref{idpt3}.
\end{proof}

\begin{subcosa}\label{DRGamI}
By~\ref{monoidal D_A}, the essential image\va1 
$\D_{\R\vG{\I}}$ of~$\R\vG{\I}$  is a monoidal category, with product\va{1.5} \smash{$\OT{\,X}$} and unit $\R\vG{\I}\OX$.  By \ref{vGI tensor-coreflects} and \ref{D_idem}, an $\OX$-complex $E$ lies in 
$\D_{\R\vG{\I}}$ if and only if the natural map is an isomorphism $\R\vG{\I}E\iso E$. 
\va1 

The plumpness of $\A^{}_{\>\I}$ in $\A$ implies that $\D_{\R\vG{\I}}\subset\D_\I$, the full subcategory of $\D$ spanned by the complexes whose homology sheaves are all in~$\A_\I$.
The converse holds if for some open coherent $\OX$-ideal~$I$,
\begin{equation*}
\I\set \{\,\textup{open coherent $\OX$-ideals }G\mid \sqrt{G}\supset  I \,\},
\end{equation*}
as can be seen, via Koszul complexes,\va1 just as in the proof of \cite[p.\,49, 5.2.1(a)]{DFS} (with the ideal $\J$ there replaced by~$I\>$). 

\begin{sprop}\label{Dqct to itself'}
For any\/ $\OX$-base\/ $\I,$ it holds that\/ $\R\vG{\I}\Dqc\subset\Dqct\>.$
\end{sprop}

\begin{proof} 
For any $E\in\Dqc\>,$ one has
\[
\R\vG{\I}E
\underset{\textup{\ref{derived intersect5}$^{\mathstrut}$}}\cong \R\vG{\I}\R\varGamma'\<E
\underset{\textup{\ref{Gamma_I}$^{\mathstrut}$}}\in \R\vG{\I}\Dqct,
\]
so one need only see that $\R\vG{\I}\Dqct\subset\Dqct$.  
Hence, one may assume that $E\in\Dqct$ and $E$ is K-injective. Then for any open immersion 
$u\colon U\hookrightarrow X$, $u^*E$~is K-injective and $u^*\vG{\I}E\cong\vG{\I|_U}u^*E$, 
so\va{-.5} one can assume $X=\Spf(S)$  for an adic noetherian ring~$S$. 
Since $\Aqct$ is closed under $\dirlm{}$,\va{-1.3} \eqref{HGam as lim.2} allows one to assume 
that $\vG{\I}=\vG{I}$ with $\I$ an open coherent $\OX$-ideal. 

Let $\kappa\colon X\to X_0\set\Spec(S)$ be the (flat) completion map 
(see proof of \ref{formal Gam preserves injective}). By~\cite[p.\,47, 5.1.2]{DFS},
$I_0\set\kappa_*I$ is a coherent $\cO_{\<\>X_0}$-ideal, and $I=\kappa^*I_0=I_0\OX$. By~\cite[p.\,50, 5.2.4]{DFS},  $E_0\set\kappa_*E\in\D_{\qc Z}(X_0)$ and $E\cong\kappa^*E_0$. Hence by \cite[p.\,53, 5.2.8(b)]{DFS},
$\R\vG{\I}E\cong\kappa^*\R\vG{\I_0}E_0$, which, by \cite[p.\,50, 5.2.4]{DFS}, lies in $\Dqct$ since 
by \ref{RGam respects Dqc},  $E_0$ being exact outside~$Z$, one has $\R\vG{\I_0}E_0\in\Dqc{}_Z(X_0)$.
\end{proof}
\end{subcosa}

Recalling that \ref{ex:supports} holds for open coherent $\OX$-ideals~$I$, let $\I$ be the $\OX$-base 
that corresponds to $\Phi_{\<\<Z}$. Then for any $\OX$-module~$M$, $\vG{\I}M\subset\vG{\<Z}M$, \emph{with equality if} $M\in\Aqct$. The proof is the same as that of \ref{2Gammas}, modulo the observation that for any 
open $U\subset X$ and $s\in\Gamma(U,M)$, \textup{ann}$_U(s)$ is an open coherent $\cO_U$-ideal.
The following more general result comes from \cite[\S\S2.1--2.2]{AJS2}.

\begin{sprop}\label{2RGams formal}
Let\/ $\I$ be the\/ $\OX$-base 
that corresponds to\/ $\Phi_{\<\<Z}$. The natural map $\theta_{\<Z,E}$ is 
an \emph{isomorphism} $\R\vG{\I}E\iso\R\vG{\<Z}E$ for all\/ $E\in\Dqct$.
\end{sprop}

\begin{proof}
As in the proof of \ref{Dqct to itself'}, one may assume $X=\Spf(S)$ for an adic noetherian ring~$S$ 
with ideal of definition, say,~$L$. Let  $\kappa\colon X= \Spf(S)\to\Spec(S)=:X_0$ and let $\I_0$ be as in the 
proof of~\ref{formal Gam preserves injective}.

By \ref{Dqct to itself}, $\R\vG{\<Z}E\in\Dqct$, and by \ref{Dqct to itself'}, $\R\vG{\I}E\in\Dqct$.
Since $\kappa_*$ sends $\Dqct$ fully faithfully into $\Dqc(X_0)$ \cite[p.\,50, 5.2.4(a)]{DFS}, 
it suffices for \ref{2RGams formal}  to show that $\kappa_*\theta_{\<Z,E}$ is an isomorphism. 

One may assume that the $\OX$-complex $E$ is K-injective, whence, $\kappa$ being flat, the $\cO_{\<\<X_{\<0}}\<$-complex $\kappa_*E$ is K-injective. 
Thus one need only verify that the following $\D(X_0)$-diagram commutes:\va{-6}
\[\mkern25mu
\def\1{$\kappa_*\vG{\I}E$}
\def\2{$\kappa_*\vG{\<Z}E$}
\def\3{$\vG{\I_0}\<\kappa_*E$}
\def\4{$\vG{\<Z}\kappa_*E$}
 \bpic[xscale=2.5,yscale=1.5]
  \node(11) at (1,-1){\1};
  \node(12) at (2,-1){\2};
  
  \node(21) at (1,-2){\3};
  \node(22) at (2,-2){\4};
 
   \draw[->] (11) -- (12) node[above ,midway, scale=.75]{$\kappa_*\theta_{\<Z\<\<,E}$} ;
   \draw[->] (21) -- (22) node[above =1,midway, scale=.75]{$\Iso$} 
                                    node[below=1,midway, scale=.75]{\ref{2RGams}} ;
   
   \draw[->] (11) -- (21) node[left=1,midway, scale=.75]{\ref{zeta}} 
                                    node[right =1,midway, scale=.75]{$\simeq$} ;   
   \draw[double distance=2] (12) -- (22) node[right =1,midway, scale=.75]{\eqref{vGam and f}, \ref{PhitoW}} ;
   
  \epic
\]
\vskip-6pt

For this it's enough, by \ref{coconds},\va{.5} to show that the natural map\va2
$\kappa_*\vG{\I}E\to \kappa_*E$ factors\- naturally as $\kappa_*\vG{\I}E\underset{\lift1.2,\ref{zeta},}\iso \vG{\I_0}\<\kappa_*E\lra \kappa_*E$, a task that\va{.5} comes down easily to verifying that the natural composite isomorphism\va{-1}
\[
\kappa_*E\iso \kappa_*\sHom_\OX\<(\kappa^*\<\cO_{\<\<X_{\<0}},E)
\iso \sHom_{\cO_{\<\<X_{\<0}}}\<\<(\cO_{\<\<X_{\<0}},\>\kappa_*E)\iso \kappa_*E\\[-1pt]
\]
is the identity map of $\kappa_*E$---which results from \cite[p.\,117, 3.5.6(e)]{Dercat}, or from an explicit description of the isomorphisms involved. Details are left to the reader. 
\end{proof}

\begin{small}
Another way to prove \Pref{2RGams formal} is by upgrading the proof of \Pref{2RGams}. This means, ultimately, to adapt the proof of \cite[p.\,25, Lemma (3.2.3)]{AJL} to the formal-scheme context. For this, two points 
have to be addressed.  
 
First, if \kf $V$ is an affine formal scheme and $g\colon V\to W$ is a~separated---hence affine---map of formal schemes then the natural map is a $\D(V)$-isomorphism $g_*\OV\iso\R g_*\OV$. 
In view of \cite[p.\,68, (13.3.1)]{GD3}, this follows from the well-known case where $V$ and $W$ are ordinary schemes. (For greater generality, see \cite[p.\,39, 3.4.2]{DFS}.)

Second, one needs to extend the \emph{projection isomorphism} to the formal-scheme context. This is done in \Pref{P:projn} below.\par
\end{small}

\vskip2pt
$\Dqct$ has a monoidal structure with product~\smash{$\OT{}$} and unit object 
$\cO'\set\R\varGamma'\OX$ (see \ref{monoid1}). For any $\OX$-base~$\I$, 
$\R\vG{\I}\cO'\cong\R\vG{\I}\OX$ (\ref{derived intersect5} with $\J\set\I'$).

\begin{sprop}\label{Dqct to itself2} 
Let\/ $Z\subset X$ be specialization-stable, and let\/ $\I$ be the\/ $\OX$-base 
corresponding to\/ $\Phi_{\<\<Z}$. Then\/ $(\R\vG{\I},\iota^{}_\I)$ and\/ 
$(\R\vG{\<Z},\iota^{}_{\<Z})$\va{.6} restrict to naturally isomorphic\/ $\ot$-coreflections of\/ $\Dqct,$ whence\/
 $(\R\vG{\I}\cO'\<\<,\>\iota^{}_\I(\cO'))$  and\/
$(\R\vG{\<Z}\cO'\<\<,\>\iota^{}_{\<Z}(\cO'))$ are naturally isomorphic\/ $\Dqct$-idempotent pairs.
\end{sprop}

\begin{proof}
By \ref{Dqct to itself'}, $\R\vG{\I}\Dqct\subset\Dqct$.
So by \ref{vGI tensor-coreflects} and \ref{coreflect exams}(a), $(\R\vG{\I},\iota^{}_\I)$ \va{.6} restricts to a 
coreflection of\/ $\Dqct$, in fact a $\ot$-coreflection because by~\ref{vGI tensor-coreflects}, subdiagram~\circled4 in the following diagram (with $\psi$ as in \ref{coreflections}, 
$\theta$ as in \ref{2RGams formal} and $F\in\Dqct(X)$) commutes: \va{-7}
\[
\def\1{\raisebox{-13pt}{$\R\vG{\I}\cO'\<\OT{} F$}}
\def\2{$\R\vG{\I}F$}
\def\3{\raisebox{-13pt}{$\R\vG{\<Z}\cO'\<\OT{} F$}}
\def\4{$\R\vG{\<Z}F$}
\def\5{\raisebox{-13pt}{$\cO'\<\OT{} F$}}
\def\6{\raisebox{-13pt}{$\OX\<\OT{X}F$}}
\def\7{$F$}
 \bpic[xscale=3.3,yscale=1.8]
  \node(11) at (1,-1){\3};
  \node(12) at (2,-1){\4};
  
  \node(21) at (1,-2){\1};
  \node(22) at (2,-2){\2};
 
  \node(31) at (0,-2.75){\5};
  \node(32) at (1.55,-2.75){\6};
  \node(33) at (3,-2.75){\7};
  
   \draw[->] (11) -- (12) node[above=1, midway, scale=.85]{$\Iso$} 
                                    node[below=1, midway, scale=.75]{$\psi^{}_Z(\cO'\!,F\>)$} ;
   \draw[->] (21) -- (22) node[above=1, midway, scale=.85]{$\Iso$} 
                                    node[below=1, midway, scale=.75]{$\psi^{}_\I(\cO'\!,F\>)$} ;
   \draw[->] (31) -- (32) node[above=1, midway, scale=.85]{$\Iso$} 
                                    node[below=1, midway, scale=.75]{$\iota'(\cO')\OT{}\b1$} ;                                 
   \draw[->] (32) -- (2.9,-2.75) node[above=1, midway, scale=.85]{$\Iso$}
                                               node[below=1, midway, scale=.75]{$\ml_F$} ;
   
   \draw[->] (1,-1.2) -- (1,-1.84) node[below=-7, midway, scale=.75]{$\mkern88mu\theta_{\<Z,\cO'}^{-\<1}\<\OT{}\b1$}
                                    node[left=1,midway, scale=.75]{$\simeq$} ;   
   \draw[->] (.67,-2.1) -- (.21,-2.62) node[below=-3, midway, scale=.75]{$\mkern90mu\iota^{}_{\I}(\cO')\<\OT{}\<\b1$} ;
   
   \draw[->] (2,-1.2) -- (22) node[left, midway, scale=.75]{$\theta_{\<Z,F}^{-\<1}$} 
                                    node[right=1,midway, scale=.75]{$\simeq$} ;   
   \draw[->] (2.2,-2.1) -- (2.91,-2.685) node[below=-3, midway, scale=.75]{$\iota^{}_{\I}(F)\mkern75mu$} ;
   
   \draw[->] (.7,-1.2) -- (.075,-2.58)  node[above, midway, scale=.75]{$\iota^{}_{\<Z}(\cO')\<\OT{}\<\b1\mkern65mu$} ;   
   \draw[->] (2.175,-1.2) -- (33) node[above=2, midway, scale=.75]{$\mkern45mu\iota^{}_{\<Z}(F)$} ;
   
  \node at (.55,-2)[scale=.9]{\circled1} ;
  \node at (1.55,-1.55)[scale=.9]{\circled2} ;
  \node at (2.38,-2)[scale=.9]{\circled3} ;
  \node at (1.55,-2.46)[scale=.9]{\circled4} ;
 \epic
\]
\vskip-6pt

It holds that 
$\iota^{}_{\<Z}(F)\smcirc\theta_{\<Z,F}=\iota^{}_\I (F)$---clearly for K-injective $F\<$, hence for all~$F\<$. 
It~follows easily that the restriction of~$(\R\vG{\<Z},\iota^{}_{\<Z})$ to $\Dqct$ is coreflecting.
Moreover, for~$F\in\Dqct$, subdiagrams \circled1 and \circled3 in the above diagram commute.
\looseness=1

 The commutativity of \circled2 follows readily from the definitions of $\psi^{}_\I$ 
and~$\psi^{}_{\<Z}$ (details left to the reader). Thus the border of the diagram commutes, and~so
$(\R\vG{\<Z},\iota^{}_{\<Z})|_{\Dqct}$~is a $\ot$-coreflection.

The last assertion results from \ref{idpt3}. (For its converse, see \ref{A and Z} below.)\va{6}
\end{proof}

\begin{sprop}\label{RG_Z} If\/ $Z$ and\/ $Z'$ are specialization-stable  subsets of\/ $X$ then for all\/ $E\in\Dqct,$ the natural map is an isomorphism
\[
\R\vG{\<Z\cap Z'}E\iso\R\vG{\<Z}\R\vG{\<Z'}E.
\]
\end{sprop}

\begin{proof} 
This follows from \ref{derived intersect5}, \ref{2RGams formal} and \ref{Dqct to itself}.\va1
\end{proof}

\begin{sprop}\label{RGam and colim5}
Let  $I$ be an\/  $\OX\<$-base, $A$ a small filtered category,  $n\in\mathbb Z,$ and\/
$\mathscr M\colon A\to \{\OX$-complexes\} a functor.\va1
The natural map is an isomorphism\va{-2}
\[
\dirlm{A} \<H^n_{\<\I}\smcirc\mathscr M\iso 
H^n_{\<\I}\big(\>\dirlm{A} \mathscr M\big).\]
In particular, $\R\vG{\I}$ commutes with small direct sums in\/ $\D(X)$.\va1
\end{sprop}

\begin{proof}
Essentially the same as the (first) proof of \ref{RGam and colim3}.
\end{proof}

\begin{sprop} \label{RGam and colim4}
Let\/ $Z$ be a specialization-stable subset of\/ $X,$ $A$ a small filtered\- category, $\mathscr M$ a functor from\/~$A$ to the category of\/ $\OX$-complexes all of whose homology modules are in\/ $\Aqct(X),$ and\/
$n\in\mathbb Z$.
The natural map is an isomorphism
\begin{equation*}
\dirlm{A} \<H^n_{\<Z}\smcirc\mathscr M\iso 
H^n_{\<Z}\<\big(\mkern1.5mu\dirlm{A}\< \mathscr M\big).
\end{equation*}
In particular, $\R\vG{\<Z}$ commutes with small direct sums in\/ $\Dqct(X)$.\va2
\end{sprop}

\begin{proof}
This follows from \ref{RGam and colim5} and \ref{2RGams formal}.
\end{proof}

\vskip5pt

\centerline{* * * * *}

\begin{scor}\label{cor:GZ+tensor}
The objects of the full subcategory\/ 
$(\Dqct)_{\R\varGamma^{}_{\<\<\!Z}\cO'}\subset\Dqct$  \(\textup{\ref{D_A}}\kf\)
are those\/ $E\in\Dqct$ with\/  \mbox{$\Supp(E)\subset Z$.}
\end{scor}

\begin{proof} For any $E\in\Dqct$,  \ref{Supp in Y} gives
\[
\Supp(E) \subset Z \!\iff\! E\in\D_{\<\Phi_{\<\<Z}}(X)\cap\Dqct\>.
\]
In view of \ref{Dqct to itself2}, \ref{ringed space D_Phi}(i) with $\D\set\Dqct$  (same proof) shows that 
$\D_{\<\Phi_{\<\<Z}}(X)\cap\Dqct$ is the essential image $(\Dqct)_{\R\varGamma^{}_{\<\<\!Z}\cO'}$
of $\>\R\vG{\<Z}\colon \Dqct\to \Dqct\mkern.5mu$.
\end{proof}

\begin{sprop}\label{SuppRGam} 
$\Supp(\R\vG{\<Z}\cO')=Z$.
\end{sprop}
 
\begin{proof}
That $\Supp(\R\vG{\<Z}\cO')\subset Z$ is given by \ref{Supp in cup}.\va1

For the opposite inclusion,\va1 let $x\in Z$, let $\widehat{\cO_{\<\<X\<,\>x}}$ be the maximal-ideal completion of~  
$\cO_{\<\<X\<,\>x}\>$, $X_x$  the one\kf-point formal scheme\va1 $\Spf(\widehat{\cO_{\<\<X\<,\>x}})$ 
($\widehat{\cO_{\<\<X\<,\>x}}$ being topologized in the usual way), \va1 $\cK_x$ the residue field of 
$\widehat{\cO_{\<\<X\<,\>x}}$ (= residue field of $\cO_{\<\<X\<,\>x}$) viewed as an object of~$\Aqct(X_x)$,
$\iota_x\colon X_x\to X$ the 
canonical map,\va{.5} and $\cK(x)\set\iota_{x*}\cK_x\in\Aqct(X)$ (see \cite[p.\,47,~5.1.1]{DFS}.)\va1

Then $\cK(x)$ has support $\overline{\{x\}}\subset Z$ and is flabby, hence K-flabby (section~\ref{flabby}), hence
$\vG{\<\<Z}$-acyclic (section~\ref{K-flabby}), so there are natural isomorphisms
\[
\cK(x)=\varGamma^{}_{\<\<\!Z}\cK(x)\iso\R\varGamma^{}_{\<\<\!Z}\cK(x)
\underset{\ref{Dqct to itself2}}\iso\smash{\R\varGamma^{}_{\<\<\!Z}\cO'\OT{X}\cK(x).}
\]  
Since the stalk $\cK(x)_x\ne0$, therefore $(\R\varGamma^{}_{\<\<\!Z}\cO')_x\ne 0$, that is,   \mbox{$x\in\Supp(\R\varGamma^{}_{\<\<\!Z}\cO')$.} 
Thus $Z\subset\Supp(\R\varGamma^{}_{\<\<\!Z}\cO')$.
\end{proof}

\begin{scor}\label{corRG_Z}
$\R\varGamma^{}_{\<\<\!Z}\cO'\preccurlyeq\R\varGamma^{}_{\!Z'}\cO'\!\iff\! Z\subset Z'$.
\end{scor}

\begin{proof}
 One has
\begin{align*}
[\,\R\varGamma^{}_{\<\<\!Z}\cO'\preccurlyeq \R\varGamma^{}_{\<\<\!Z'}\cO'\,]
&\underset{\under{.7}{\textup{\ref{preccurly}}}}{\!\iff\!}
[\,\R\varGamma^{}_{\<\<\!Z}\cO'\in(\Dqct)_{\R\varGamma^{}_{\mkern-5mu \lift1.2,\sst Z\<',}\<\cO'}]\\[2pt]
&\underset{\under{.7}{\textup{\ref{cor:GZ+tensor}}}}{\!\iff\!}
[\,\Supp(\R\varGamma^{}_{\<\<\!Z}\cO')\subset Z'\,]
\underset{\under{.7}{\textup{\ref{SuppRGam}}}}{\!\iff\!} 
Z\subset Z'.\\[-34pt]
\end{align*}
\end{proof}
\vskip10pt
As in \ref{Dqct to itself2}, for any specialization-stable $Z\subset X$, the pair 
$\bigl(\R\varGamma^{}_{\<\<\!Z}\cO'\<,\>\iota^{}_{\<Z}(\cO')\bigr)$ is $\Dqct$-idempotent. The  converse is, essentially, \cite[p.\,604, Corollary 5.4]{AJS2}:

\begin{sprop}\label{A and Z}
Every\/ $\Dqct$-idempotent pair\/ $(A,\>\alpha)$ is isomorphic to a pair\/ 
$\bigl(\R\varGamma^{}_{\<\<\!Z}\cO'\<\<,\>\iota^{}_{\<Z}(\cO')\<\bigr)$ 
for some specialization-stable\/ $Z$---necessarily\/ $\Supp(A)$ \(\kern-1pt see~\textup{\ref{SuppRGam}}\).
\end{sprop}

\begin{proof}
The idea is to reduce, via localization and completion, to the known case where $X=\Spec(S)$ for a noetherian ring $S$. 

Since $H^i\<A\in\Aqct\subset \Avc\>\>$ is the union of all its coherent submodules, each of whose support 
is closed (see \S\ref{suppdef}), therefore $\Supp(A)=\cup_{i\in\mathbb Z}\, \Supp(H^i\<A)$ \emph{is specialization-stable}.
So by  \ref{SuppRGam} and \ref{isos}, it's enough to show that if $(A,\alpha)$ and~$(B,\beta)$ are
$\Dqct$-idempotent pairs with $\Supp(A)=\Supp(B)$, then $A\cong B$, i.e., by \ref{one map}, the maps 
$A\OT{X}B\to B$ and $A\OT{X}B\to A$ induced\va{-4} by $\alpha$ and~$\beta$, respectively,
induce homology isomorphisms. In view of \ref{f*idpt2} with $f$ an open immersion, the question is local, so one may assume that $X=\Spf(S)$ where $S$ is an adic noetherian ring with ideal of definition, say, $L$.

Let $\kappa\colon X\set\Spf(S)\to\Spec(S)\set X_0$ be the completion map
(as in the proof~of~\ref{formal Gam preserves injective}).  This map is flat,  so
the functor $\kappa^*\colon\A(X_0)\to\A(X)$ is exact, therefore extends\- to $\kappa^*\colon\D(X_0)\to\D(X)$. 

Since $X$ and $Y\set\Spec(S/L)$ are homeomorphic, and ~$\kappa$ is, topologically, the inclusion 
$Y\hookrightarrow X_0\>$, one can regard~$Y$ as a closed subset of~$X_0\>$, whence the 
functor $\kappa_*\colon\A(X)\to\A(X_0)$ is exact, therefore extends to $\kappa_*\colon\D(X)\to\D(X_0)$.
The category $\Dqc(X_0)$ has a monoidal structure \va1 
with product \smash{$\OT{\mkern2.5muX_0}\!\<\<$} and unit object $\cO_{\<\<X_0}$ (see~\ref{exams-mon}(a)). 
For any $W\subset X_0$, the pair $(\R\vG{W}\cO_{\<\<X_0},\> \iota^{}_W(\cO_{\<\<X_0}))$\va{.6} is $\Dqc(X_0)$-idempotent (see~ \ref{RGamidem}). So~by \ref{monoidal D_A},
$\D_{\qc W}(X_0)\set(\Dqc(X_0))_{\R\vG{W}\cO_{\<\<X_0}}\!\!$ has a monoidal
structure with product~\smash{$\OT{\>\>X_0}\!\<\<$} and unit object $\R\vG{W}\cO_{\<\<X_0}$. 
Also,\va1 if $W'\subset W$ then $\R\vG{W'}\cO_{\<\<X_0}\preccurlyeq \R\vG{W}\cO_{\<\<X_0}$, so by \ref{preccurly}(ii),
$\R\vG{W'}\cO_{\<\<X_0}$ is $\D_{\qc W}(X_0)$-idempotent.\va2

By \ref{cor:GZ+tensor} for the discrete formal scheme~$X_0$, 
$\D_{\qc Y}(X_0)$ is the full subcategory spanned~by the 
$\cO_{\<\<X_0}\<$-complexes whose homology modules are quasi-coherent and have support contained in~$Y\<$. Thus the notation here agrees with that in \cite[p.\,50, 5.2.4(a)]{DFS},  which shows that the functors $\kappa^*$ and~$\kappa_*$ 
give inverse iso\-morphisms between $\D_{\qc Y}(X_0)$ and $\Dqct(X)$. Also, one has
the usual isomorphism 
\[ 
\smash{v(E\<,F\>)\colon\kappa^*\<\<E\OT{\<\<X} \kappa^*\<\<F\iso\kappa^*\<(E\OT{\<\<X_0}F)
\qquad(E, F\in\D(X_0)).}
\]
Taking $E\set \R\vG{Y}\cO_{\<\<X_0}$, one deduces that $\kappa^*\R\vG{Y}\cO_{\<\<X_0}$ is a unit object in the monoidal category $\Dqct$ (see \ref{monoid1}); and one checks (directly,  
or via~\ref{f*idpt} with $\cO_2\set\R\vG{Y}\cO_{\<\<X_0}$, $\cO_1\set\kappa^*\cO_2$, and $u\set\b1_{\cO_1}$)
that $\kappa^*$ and $\kappa_*$ induce inverse bijections between the sets of
$\D_{\qc Y}(X_0)$-idempotents and $\Dqct(X)$-idempotents.\va1
 
Let $Z\subset Y$ be specialization-stable.
Over $X_0$, set $\I_0\set\I_{\<\Phi_{\<\<Z}}$ (see~\ref{PhisubI}, \ref{supports}). The $\OX$-base corresponding to the s.o.s.~$\Phi_{\<\<Z}$ in $X$ is 
\[
\I=\I_0\OX\set\{\,I_0\OX\mid I_0\in\I_0\,\}
\]
(see \ref{I and J}). 
For $F\in\D(X)$, one has, as in the lines following \ref{zeta}, the iso\-morphism
$\bar \xi_{\kappa\<,\I\<,F}\colon\kappa_*\R\vG{\I}F\iso \R\vG{\I_0}\<\kappa_*F$,
whence the natural composite map
\begin{equation*}\label{A and Z1}
\kappa^*\R\vG{\I_0}F
\lra 
\kappa^*\R\vG{\I_0}\kappa_*\kappa^*\<\<F
\xto[\lift1.2,\kappa^*\bar \xi^{-\<1},]{\Iso}
\kappa^*\kappa_*\R\vG{\I}\kappa^*\<\<F
\lra
\R\vG{\I}\kappa^*\<\<F.\\[-1pt]
\tag*{(\ref{A and Z}.1)}
\end{equation*}
This is an \emph{isomorphism}:
apply cohomology \mbox{$H^n\ (n\in\mathbb Z)$,} then use~ \eqref{HGam as lim.2} to reduce\- to where 
there is an open coherent $\cO_{\<\<X_0}$-ideal~$I_0$ such that 
\[
\I_0\set \{\,\textup{open coherent $\OX$-ideals }J\mid \sqrt{J}\supset  I_0 \,\};
\] 
then identify \ref{A and Z1}, via \cite[p.\,18, 3.1.1]{AJL},
with the natural isomorphism $v(\KK(\bt), F)$,
where $\bt$ is a finite sequence in $S$ that generates $I_0$, and $\KK(\bt)$ is as in the proof of \ref{injective}. 
(Details left to the reader.)%
\footnote{Alternatively, it is an instance of the isomorphism in \ref{lift vGI}. }

One has then the composite isomorphism
\[
\kappa^*\R\vG{Z}\cO_{\<\<X_0}
\<\underset{\ref{2RGams}}\iso\<
\kappa^*\R\vG{\I_0}\cO_{\<\<X_0}\>
\underset{\ref{A and Z1}}\iso\>
\R\vG{\I}\OX\<\<
\underset{\ref{derived intersect5}}\iso
\<\R\vG{\I}\R\varGamma'\OX\<=\\R\vG{\I}\cO'\<
\<\<\underset{\ref{Dqct to itself2}}\iso\<
\R\vG{Z}\cO'\<.
\]

Accordingly, it suffices that any $\D_{\qc Y}(X_0)$-idempotent be 
isomorphic to $\R\vG{Z}\cO_{\<\<X_0}$ for some specialization-stable $Z\subset Y$. But in view of the monoidal equivalence between $\Dqc(X_0)$ and $\D(S)$ 
\cite[p.\,225, Theorem 5.1]{BN}, that results, essentially, from 
\cite[p.\,526, Theorem 2.8]{chrom} and its proof---cf.~\cite[Proposition 3.5.7\kf]{Lectures} and the remarks following it, keeping in mind the bijection between specialization-stable subsets of~$\Spec(S)$ and decent topologies on $S$, as in \ref{toprings} above.\va1 
\end{proof}

\begin{scor}\label{idem and specstable}
The mapping that takes the isomorphism class of\/ 
$A$ to\/~$\Supp(A)$ induces an order-preserving bijection
\[
\textup{\{isomorphism classes of\/ 
$\Dqct$-idempotents\} \!$\leftrightarrow$\!
\{specialization-stable subsets of~$X$\}.}
\]
\end{scor}

\begin{proof} 
This follows from \ref{Dqct to itself2}, \ref{corRG_Z} and \ref{A and Z}.\va1
\end{proof}

\begin{scor}\label{idem and bases}
 There is an order-reversing bijection
\[
\textup{\{$\OX$-bases\} $\leftrightarrow$
\{isomorphism classes of\/ 
$\Dqct$-idempotents\}}
\]
that sends an\/ $\OX\<$-base\/ $\I$ to the isomorphism class of\/ $\R\vG{\I}\OX$.
\end{scor}

\begin{proof} The order-reversing bijection arising from \ref{ex:supports} takes
each $\OX$-base~$\I$ to the specialization-stable set $Z=\cup_{I\in\I}\>Z(I)$; and the order-preserving bijection in~\ref{idem and specstable} takes $Z$ to the isomorphism class of
$\R\vG{\<Z}\cO'\cong \R\vG{\I}\cO'\cong\R\vG{\I}\OX$ (see \ref{2RGams formal} and~\ref{derived intersect5}).\va1
\end{proof}

\begin{scor}\label{Supp in Z}
If\/ $A$ is a\/ $\Dqct(X)$-idempotent and\/ $E\in\Dqct(X),$ then
\[
E\in (\Dqct)_{\<\<A}\!\iff\!\Supp(E)\subset\Supp(A).
\]
\end{scor}

\begin{proof}
With $Z\set\Supp(A)$, \Pref{A and Z} allows one to assume $A=\R\varGamma^{}_{\<\<\!Z}\cO'\mkern-1.5mu$,
in which case the assertion is just  Corollary~\ref{cor:GZ+tensor}.\va2
\end{proof}

\begin{scor}\label{x in Supp} Let\/ $x\in X$ and let\/ $\cK(x)\in\Dqct(X)$ be as in the proof 
of~\textup{\ref{SuppRGam}}. For any\/ $\Dqct(X)$-idempotent\/ $(A,\alpha),$ 
\[
\cK(x)\!\in\<\<(\Dqct)_{\<\<A}\!\iff\! A\OT{\<X}\<\cK(x)\<\<\underset{\via\alpha}\cong\cK(x)\!\iff\! A\OT{\<X}\<\cK(x)\ne0
\!\iff\! x\in\Supp(A).
\]
\end{scor}

\textsc{Proof.}
$\Supp(A)$ is specialization-stable,\va{1.5} so if $x\in\Supp(A)$ then
$\overline{\{x\}}\subset\Supp(A)$,
whence by \ref{Supp in Z}, $\cK(x)\in(\Dqct)_{\<\<A}\>$, that is, by \ref{D_idem}, 
$A\OT{X}\cK(x)\underset{\lift1.2,\via\alpha,}\cong\cK(x)$, whence \smash{$A\OT{X}\cK(x)\ne0$}, whence, by \ref{A and Z},
if $Z\set\Supp(A)$ then\va3 \smash{$\R\varGamma^{}_{\<\<\!Z}\cO'\OT{X}\cK(x)\ne0$.} 

Conversely, by 
\ref{Dqct to itself2} and since $\cK(x)$ is flabby, hence K-flabby,\va1 and since $x$ lies in the support of any nonzero section of $\cK(x)$ over any open set $U\<$, therefore
\[
0\ne\R\vG{\<\<Z}\cO'\OT{X}\cK(x)
\cong\R\vG{\<\<Z}\cK(x)\cong\vG{\<\<Z}\cK(x)
\implies x\in Z=\Supp(A).\makebox[0pt]{$\mkern100mu\square$}
\]

\centerline{* * * * *}
\medskip

The next Proposition enhances \Cref{idem and specstable}. 

Set $\cO'_{\!X}\set\R\varGamma'_{\!\<\<X}\OX$ and\/ $\cO'_{\<\<W}\set\R\varGamma'_{\!W}\OW$.

\begin{sprop}\label{liftidem}
Let\/ $f\colon W\to X$ be a map of noetherian formal schemes. Set\/~$\bL'\!f^*\set\R\varGamma'_{\!W}\bL f^*\<\<$. Then the natural map is an isomorphism 
\[
\bL'\!f^*\<\cO'_{\!X}=\bL'\!f^*\R\varGamma'_{\!\<\<X}\OX
\iso\bL'\!f^*\OX=\cO'_W\>;
\]
and if\/ $\alpha\colon A\to \cO'_{\!X}$~is a\/ 
$\Dqct(X)$-idempotent pair with\/ $\Supp(A)=Z,$ then\/ 
\[
\bL'\!f^*\!\alpha\colon \bL'\!f^*\!\<A\lra \bL'\!f^*\<\cO'_{\!X}=\cO'_W   
\]
is a $\Dqct(W)$-idempotent pair with\/ $\Supp(B)=f^{-1}Z\<$. 
\end{sprop}

\begin{proof}
For the first assertion, see \cite[p.\,53, 5.2.8(c)]{DFS}.

The pair $(A,\alpha)$ is clearly $\D(X)_{\<\cO'}$-idempotent, whence
$A\xto{\,\alpha\,}\cO'_{\!X}\xto{\!\iota'\<(\OX\<)\<}\OX$ is $\D(X)$-idempotent, see 
\ref{vGI tensor-coreflects},  \ref{preccurly}. By \ref{f*idpt2} and
 \ref{Gam preserves idem} with $(\Gamma,\iota)\set(\R\varGamma'_{\!W}, \iota'_W)$ (see 
again \ref{vGI tensor-coreflects}), the composition
\[
\bL'\!f^*\!A\xto{\,\bL'\!f^*\!\alpha\,}\bL'\!f^*\<\cO'_{\!X}
\xto{\!\bL'\!f^*\!\iota'_{\<\<X}\<\<(\OX\<\<)\<}\bL'\!f^*\<\OX
\xto{\iota'(\bL f^*\<\<\OX)}\bL  f^*\<\<\OX=\OW
\]
is $\D(W)$-idempotent. In particular, this holds when $\alpha=\b1_{\cO'_{\!X}}$. 

Note that $\bL f^*\!A\in\Dqc(W)$: the question being local on $X\<$, one can assume\- that $A\in\Dvc(X)$ and apply
\cite[p.\,37, 3.3.5]{DFS}. So by \ref{Dqct to itself'}, $\bL'\!f^*\!\<A$ in $\Dqct(W)$.
Hence, as in the proof of \ref{preccurly}(ii), 
with $B\set\bL'\!f^*\!A$,  $A\set\bL'\!f^*\<\<\OX$ and $\lambda\set\bL'\!f^*\alpha$,
\[
\bL'\!f^*\!\alpha\colon \bL'\!f^*\!\<A\to \bL'\!f^*\<\cO'_{\!X}=\cO'_W
\]
is $\D(W)_{\cO'_{\<\<W}}$-idempotent, and thus $\Dqct(W)$-idempotent.

It remains to be shown that  $B\set\bL'\!f^*\!\<A$ has support $f^{-1}Z$.

Assuming, as one may, that $A$ is K-flat, one has for $w\in W$ and $x\set f(w)$ that 
$$
(\bL f^*\!\<A)_w= (f^*\!\<A)_w=\cO_{W\!,\>w}\ot_{\cO_{\<\<X\!,x}}\< A_{x}.
$$
If $x\notin Z$ then $A_{x}$ is exact and K-flat,%
\footnote{For any exact $\cO_{\<\<X\!,x}$-complex~$C$,\vs{-1} the extension by 0 of the
constant sheaf~$C$ on $\bar x$
is exact, as is its tensor product
with the K-flat $\OX$-complex~$A$, whence $C\ot_{\cO_{\<\<X\!,x}}\< \<A_{x}$
is exact.}
and therefore 
$(\bL f^*\!\<A)_w$ is exact---as one sees upon replacing  $\cO_{W\!,\>w}$ by a quasi-isomorphic 
K-flat $\cO_{\<\<X\!,x}\>$-complex; in other words, 
$w\notin\Supp(\bL f^*\!\<A)$. Hence $\Supp(B)\subset f^{-1}Z$.

For the opposite inclusion, suppose $w\in f^{-1} Z\setminus \Supp(B)$. Let \mbox{$\cK(w)\in\Aqct(W)$}
be as in the proof of~\ref{SuppRGam}. This sheaf is K-flabby, hence \mbox{$\fst$-acyclic} (see section~\ref{K-flabby}). 
One has
$\R\fst\cK(w)\cong \fst\cK(w)\in\Aqct(X)$ \cite[p.\,47, 5.1.1]{DFS}, 
the stalk $(\fst\cK(w))_{x}$ is the residue field of $\cO_{\<\<W\<\<,\>w}\>$, and 
$\fst\cK(w)$ vanishes outside
$\,\overline{\{f(w)\}}\,\subset Z\<$, so that $\Supp(\R\fst\cK(w))\subset \Supp(A)$. Hence
$\smash{0\ne\R\fst\cK(w)\cong A\OT{X}\R\fst\cK(w),}$
where the last isomorphism comes from \ref{Supp in Z} and \ref{D_idem}.

Since $A\<\in\<\Dqct(X)\<\subset\Dvc(X)$ and 
$\cK(w)\<\in\<\Dqct(W)\<\subset\Dvc(W)$, 
\Pref{P:projn} below gives a natural ``projection"  isomorphism
$$
\smash{0\ne A\OT{X}\R\fst\cK(w)\cong\R\fst\big(\bL f^*\!\<A\OT{W}\cK(w) \big);}
$$
so one has, via \ref{Psi for coref}, \ref{vGI tensor-coreflects} and \ref{Gama'}, natural isomorphisms 
\begin{align*}
\smash{0\ne \bL f^*\!\<A\OT{W}\cK(w)} &\cong
\smash{\bL f^*\!\<A\OT{W}\R\varGamma'_{\!W}\cK(w)}\\
&\cong
 \smash{\R\varGamma'_{\!W}\big(\bL f^*\!\<A\OT{W}\cK(w)\big)}\\
&\cong
\smash{\R\varGamma'_{\!W}\bL f^*\!\<A\OT{W}\cK(w)=B\OT{W}\cK(w)}=0,
\end{align*}
where the last equality comes from \ref{x in Supp}.
This contradiction shows $w$ can't exist. 
Thus $\Supp(B)=f^{-1} Z.$
\end{proof}

\begin{sprop}\label{lift vGI} 
Let\/ $f\colon W\to X$ be a map of noetherian formal schemes, let\/~$\I$ be an\/ $\OX$-basis, and let\/ $\I_{\!f}$ be the\/ $\OW$-basis
\[
\I_{\!f}\set\{\,\textup{open coherent $\OW$-ideals }J\mid \sqrt{J}\supset  I\OW \textup{ for some }I\in\I\,\}.
\]
There is a unique functorial isomorphism
\[
\xi(\I, E)\colon\bL f^*\R\vG{\I}E\iso\R\vG{\I_{\!\<f}}\bL f^*\<\<E \qquad(E\in\D(X))
\]
whose composition with the natural map\/ $s\colon\R\vG{\I_{\!\<f}}\bL f^*\<\<E\to \bL f^*\<\<E$ is 
the natural map\/ $q\colon\bL f^*\R\vG{\I}E\to  \bL f^*\<\<E$.
\end{sprop}

\begin{proof} 
Set $\D\set\D(W)$. First of all,\va{.5} 
it holds that $\bL f^*\R\vG{\I}E$ \emph{is in the essential\va{-1} image\/~$\D_{\R\vG{\I_{\!\<f}}}$ of\/}~
$\R\vG{\I_{\!f}}$---which implies the existence\va{-1} and uniqueness of $\xi(\I,E)$ as a $\D(W)$-\emph{map} 
(see~\ref{vGI tensor-coreflects} and \ref{coconds}). \va1

To see this, assume without loss of generality that $E$ is K-injective. Regard the ordered set $\I$ as a category in the usual way (with containments as morphisms), and let~$P$ be a functor from $\I$ to 
the category of maps of $\OX$-complexes such that for each
 $I\in\I$, $P(I):P_I\to \vG{I}E$ is a K-flat resolution, and for each $\I$-morphism $I'\supset I$, the resulting map 
 $\vG{I'}E\to\vG{I}E$ is the natural one. The existence of such a~$P$ is given, for instance,  by \cite[p.\,61, 2.5.5]{Dercat}.
Then with $\dirlm{}\!\<\set\dirlm{I\in\I^{}\,}\!$, $\dirlm{}\!P_I$ is a K-flat resolution of $\dirlm{}\!\vG{I}E=\vG{\I}E$; and  so
\[
\bL f^*\<\<\R\vG{\I}E \cong f^*\dirlm{}\!P_I\cong \dirlm{}\<\<f^*\<\<P_I.
\]

If for an open coherent $\OX$-ideal~$I$,  $\I$ is the $\OX$-basis
\[\label{I_I}
\I_{\<\<I}\set \{\,\textup{open coherent $\OX$-ideals }G\mid \sqrt{G}\supset  I \,\},\tag*{(\ref{lift vGI}.1)}
\]
then
\[
\I_{\!f}=\I_{\<\<I\<\OW\!}=\{\,\textup{open coherent $\OW$-ideals }J\mid \sqrt{J}\supset  I\OW\,\}.\\[2pt]
\] 
So in this case, $\xi(\I,E)$ is the isomorphism given by \cite[p.\,53, 5.2.8(b)]{DFS}, 
whence 
\[\label{f^*PI}
f^*\<\<P_I\cong\bL f^*\R\vG{I}E\cong\R\vG{I\OW}\bL f^*\<\<E\!
=\R\vG{\I_{\!f}}\bL f^*\<\<E\in\D_{\R\vG{\I_{\!\<f}}}.\tag*{(\ref{lift vGI}.2)}
\]

The category $\D_{\R\vG{\I_{\!\<f}}}$  is a triangulated subcategory of $\D(X)$: it is clearly closed under translation, and if $\>T$ is a $\D(X)$-triangle with two vertices in~$\D_{\R\vG{\I_{\!\<f}}}\mkern-1.5mu$, then since
$\R\vG{\I_{\!\<f}}$ is coreflecting (see \ref{vGI tensor-coreflects}), the natural map is an isomorphism\va2 $\R\vG{\I_{\!\<f}}T\iso T$, so the third vertex is also in $\D_{\R\vG{\I_{\!\<f}}}\<$. 

Moreover,  \ref{RGam and colim5} implies that 
$\D_{\R\vG{\I_{\!\<f}}}$ is closed under small direct sums, that~is, $\D_{\R\vG{\I_{\!\<f}}}$ is a localizing subcategory of $\D$.
So \cite[pp.\,232--233, Theorems 2.2 and~3.1]{AJS1} give that
\[
\bL f^*\R\vG{\I}E\cong\bL f^*\<\<\vG{\I}E \cong\dirlm{\lift .25,{\halfsize{$I$}\in{\lift1,\I,}},\,}\<f^*\<\<P_I\in\D_{\R\vG{\I_{\!\<f}}}\<,
\]
\vskip-4pt\noindent
as desired.\va2

Note next that the natural map 
\[
\nu^{}_{\mkern-1.5muI}\colon
\dirlm{\lift .25,{\halfsize{$I$}\in{\lift1,\I,}},\,}\<\<H^n\bL f^* \R\vG{I}E
\lra 
H^n\bL f^*\R\vG{\I}E\\[3pt]
\]
is the natural composite \emph{isomorphism}
\[
\dirlm{\lift .25,{\halfsize{$I$}\in{\lift1,\I,}},\,}\<\<H^n\bL f^* \R\vG{I}E
\iso
 \dirlm{\lift .25,{\halfsize{$I$}\in{\lift1,\I,}},\,}\<\<H^n \<\<f^*\<\< P_I
 \iso
 H^n\<\< f^*\> \dirlm{\lift .25,{\halfsize{$I$}\in{\lift1,\I,}},\,}\!P_I
 \iso
 H^n\bL f^*\R\vG{\I}E.\\[4pt]
 \]
So for $\xi(\I,E)$to be an isomorphism it's enough that for all $n\in\mathbb Z$,
subdiagram~\circled1 of the following natural diagram commutes, or equivalently, that the border of the whole diagram commutes.
\[
\def\1{\raisebox{-9pt}{$H^n\bL f^* \R\vG{I}E$}}
\def\2{\raisebox{-12pt}{$\dirlm{\lift .25,{\halfsize{$I$}\in{\lift1,\I,}},\,}\<\<H^n\bL f^* \R\vG{I}E$}}
\def\3{\raisebox{-9pt}{$H^n\bL f^*\R\vG{\I}E$}}
\def\4{\raisebox{-12pt}{$H^n\R\vG{I\OW}\bL f^*\<\<E$}}
\def\5{\raisebox{-12pt}{$\dirlm{\lift .25,{\halfsize{$I$}\in{\lift1,\I,}},\,}\<\<H^n \R\vG{I\OW}\bL f^*\<\<E$}}
\def\6{\raisebox{-12pt}{$H^n\R\vG{\I_{\!f}}\bL f^*\<\<E$}}
 \bpic[xscale=3.75,yscale=1.65]
  \node(11) at (1,-1){\1} ;
  \node(12) at (2,-1){\2} ;
  \node(13) at (3,-1){\3} ;
  
  \node(21) at (1,-2){\4} ;
  \node(22) at (2,-2){\5} ;
  \node(23) at (3,-2){\6} ;
  
   \draw[->] (11) -- (12) ;
   \draw[->] (12) -- (13) node[below,midway, scale=.75]{$\via\>\nu^{}_{\<\<I}$} 
                                    node[above=.7,midway, scale=.75]{$\Iso$} ;
                                    
   \draw[->] (21) -- (22) ;
   \draw[->] (22) -- (23) node[below=1,midway, scale=.75]{\eqref{HGam as lim.2}}
                                     node[above=.7,midway, scale=.75]{$\Iso$} ;   
   \draw[->] (11) -- (1,-1.8) node[left,midway, scale=.75]{$H^n\xi(\I_{\<\<I},E)$}
                                         node[right, midway, scale=.75]{$\simeq$}  ;
   \draw[->] (2,-1.25) -- (2,-1.8) node[left,midway, scale=.75]{$\dirlm{\lift .25,{\halfsize{$I$}\in{\lift1,\I,}},\,}\<\<H^n\xi(\I_{\<\<I},E)^{\mathstrut}$}
                                                 node[right, midway, scale=.75]{$\simeq$} ;
   \draw[->] (3,-1.25) -- (3,-1.8) node[right,midway, scale=.75]{$H^n\xi(\I,E)$}; 
   
  \node at (2.55,-1.55){\circled1} ;   
 \epic
\]  
Thus one has to see that in the following natural diagram, $ca=db$. 
\[
\def\0{$\bL f^*\<\<E$}
\def\1{\raisebox{-9pt}{$\bL f^* \R\vG{I}E$}}
\def\3{\raisebox{-9pt}{$\bL f^*\R\vG{\I}E$}}
\def\4{\raisebox{-9pt}{$\R\vG{I\OW}\bL f^*\<\<E$}}
\def\6{\raisebox{-12pt}{$\R\vG{\I_{\!f}}\bL f^*\<\<E$}}
 \bpic[xscale=2.5,yscale=1.8]
  \node(11) at (1,-1){\1} ;
  \node(13) at (3,-1){\3} ;
  
  \node(21) at (1,-2){\4} ;
  \node(23) at (3,-2){\6} ;
  
  \node(00) at (2,-1.5){\0} ;
  
   \draw[->] (11) -- (13) node[above, midway, scale=.75]{$a$} ;
   
   \draw[->] (21) -- (23) node[below, midway, scale=.75]{$d$} ;
   
   \draw[->] (11) -- (1,-1.8) node[left=1, midway, scale=.75]{$b=\xi(\I_{\<\<I},E)$} ;
   \draw[->] (3,-1.25) -- (3,-1.8) node[right=1, midway, scale=.75]{$c=\xi(\I,E)$} ; 
   
  \draw[->] (11) -- (00) node[above, midway, scale=.75]{$p$} ; 
  \draw[->] (13) -- (00) node[above, midway, scale=.75]{$q$} ;
  \draw[->] (21) -- (00) node[below, midway, scale=.75]{$r$} ;
  \draw[->] (23) -- (00) node[below, midway, scale=.75]{$s$} ; 
  
 \epic
\]  
\vskip-3pt\noindent
But  
$\bL f^*\R\vG{I}E\in\D_{\R\vG{\I_{\!\<f}}}\mkern-1.5mu$ (see \ref{f^*PI}) so by \ref{coconds},\va{-1} 
it suffices to note that since 
(clearly) all the subdiagrams commute, therefore $sca=qa=p=rb=sdb$.
 \end{proof}

\begin{small}
As an exercise, show that  for $E=\OX$  and
$A\in\Dqct(X),$ \ref{lift vGI} implies \ref{liftidem}.
\end{small}   
   
\medskip
Recall that $\cO_{\!X}'\set\cO'=\R\varGamma'\OX$, that $(\cO_{\!X}',\iota'_{\<\<X}\<(\OX))$ is $\D(X)$-idempotent (see \ref{vGI tensor-coreflects}), and that for any ringed-space map $f\colon W\to X$, 
$(\bL f^*\cO_{\!X}',\bL f^*\<\iota'_{\<\<X}\<(\OX))$ is $\D(W)$-idempotent (see \ref{f*idpt} and the remarks following it). Recall also the meaning of $B\preccurlyeq A$ for $\OW$-idempotents $B$ and $A$ (see \ref{defcurly}).

\begin{sprop}\label{rho and psi2}  
Let\/ $f\colon W\to X$ be a map of noetherian formal schemes. Then $\cO_{W}'\preccurlyeq \bL f^*\cO_{\!X}'\>;$ 
and $\cO_{W}'\cong \bL f^*\cO_{\!X}'$ if and only if\/ $f$ is adic.
\end{sprop}

\textsc{Proof.}
Let $I$ (respectively~$J$) be an ideal of definition of $X$ (respectively~$W$). 
Then $I\OW\subset\sqrt{J}$ (see \cite[p.\,416, (10.6.10)(i)]{GD}\kf), or equivalently, $\I_{\!f}=\I_{\<I\OW}\supset\I_{\<\<J}$
(see~\ref{I_I}\,\emph{ff.}), 
or equivalently (by \ref{idem and bases}), 
\[
\mspace{100mu}\cO'_W=\R\vG{\<\<J}\>\OW\preccurlyeq \R\vG{I\OW}\OW\!\<
\underset{\ref{lift vGI}}\cong\bL f^*\R\vG{\<I}\OX =\bL f^*\cO_{\!X}'.
\mspace{87mu}\square
\]

If, in addition, $\bL f^*\cO_{\!X}'\preccurlyeq \cO_{W}'$, i.e., $J\subset\sqrt{I\OW}$,\va{.6} then $I\OW$ is an ideal of definition of~$W\<$, and so $f$ is adic \cite[p.\,436, (10.12.1)]{GD}.\va5

\centerline{* * * * *}
\vskip5pt

The following basic facts, \ref{monoid1} and ~\ref{P:projn}, were referred to before.

\begin{sprop}\label{monoid1}
\textup{(i)} For an ordinary scheme\/ $X,$ the usual monoidal structure on\/ $\D(X)$ restricts to one on\/ $\Dqc(X)$.

\textup{(ii)} For a noetherian formal scheme\/ $X,$  the usual monoidal structure on\/ $\D(X)$ restricts to one on\/~ $\Dvc(X)$. 

\textup{(iii)} For a noetherian formal scheme\/ $X,$ there is a monoidal structure on \/~ $\Dqct(X)$ with product map and associativity and symmetry isomorphisms inherited from the usual monoidal structure on\/ $\D(X),$
with unit element\/  $\cO'\set\R\varGamma'\OX,$
and with unit isomorphisms\/ \smash{$\ml'_E\set\ml_E\smcirc (\iota'(\OX)\OT{}\b1_E)$} and\/~ 
$\mr'_E\set\mr_E\smcirc (\b1_E\OT{}\iota'(\OX))\ (E\in\Dqct)$\va{-2}.
\end{sprop}

\begin{proof}
One needs to show that $\Dqc(X)$ (resp.~$\Dvc(X)$, $\Dqct(X)$)\va1 is \smash{$\OT{\mkern2.5mu X}$}-closed,
i.e.,~if all the cohomology sheaves of $\OX$-complexes $E$ and $F$ lie in $\Aqc\set\Aqc(X)$ 
(resp.~$\Avc\set\Avc(X)$, $\Aqc\set\Aqct(X)$) then the same holds for \smash{$E\OT{X}F$}.
After that, one applies \ref{monoidal D_A} with $\D_{\<*}\set \Dqc$ and $\alpha\set\b1_\OX\!$ or $\Dvc\>$, or (keeping in mind \ref{Gamma_I}) with $\D_{\<*}\set\Dqct$ and 
$\alpha\set\iota'(\OX)\colon\cO'\to\OX$.\va1

Let $P\to E$, $P'\to F$ be K-flat $\OX$-resolutions, so that, with $P^{\lle
u}$ the complex obtained from $P$ by replacing $P^n$ by 0 for all
$n>u$ and $P^u$ by the kernel of $P^u\to P^{u+1}\<$, one has, for all $i\in\mathbb Z$,\va{-2}
$$
H^i(E\OT{\<\<X}\< F)\cong H^i(P\otimes\sX\< P')=\dirlm{u}H^i(P^{\lle u}\otimes\sX \<P').\\[-3pt]
$$
Since $\Aqct$ is closed under $\dirlm{}$ \cite[p.\,48, 5.1.3]{DFS},\vspace{-1pt} and the same clearly holds 
for~$\Avc\>$ and for $\Aqc$, therefore 
$P$ can be replaced by a bounded-above\va{.5} flat resolution
of~$P^{\lle u}\<$. Then one can do likewise with~$P'\<$. 
So one may assume $E$ and~$F$ bounded-above.
Since $\Dqct(X)$ (resp.~$\Dvc(X)$, $\Dqc(X)$) is a triangulated subcategory of $\D(X)$ 
(\cite[p.\,48, 5.1.3]{DFS}, \cite[p.\,34, 3.2.2]{DFS}, \cite[p.\,217, (2.2.2)\kf(iii)]{GD}),
\cite[p.\,73, Proposition 7.3(ii)]{H} (dualized, and
for whose terminology see \cite[p.\,38, Definition]{H})
yields a further reduction to where $E$ and $F$ are single
sheaves (complexes vanishing in nonzero degrees). To be shown then is
that ${\T}or_i^X(E,F\>)\in\Aqct$ (resp.~$\Avc\>$, $\Aqc$).

For $\Aqc$ the problem is local, say $X=\Spec\>(R)$, and is easily
disposed of via the standard equivalence of categories between $\Aqc$
and the category of $R$-modules (an equivalence which preserves free
resolutions). For $\Aqct$, one has more generally that if 
$E\in\Aqc$ and $F\in\Aqct$ then  ${\T}or_i^X(E,F\>)\in\Aqct\>$: one localizes to the case
$X=\Spf(S)$ where $S$ is a noetherian ring complete with respect to
the topology defined by powers of an ideal $I$, such that $E$ is a cokernel of a map of
free $\cO\sX$-modules, so that $E\in\Avc(X)$ \cite[p.\,32, 3.1.4]{DFS}; and one uses the
equivalences of categories described\vs1 in \cite[p.\,31, 3.1.1]{DFS} and 
\cite[p.\,47, 5.1.2]{DFS}\va1 to reduce to consideration of two 
$S$-modules $E_0$, $F_0$, such that
$
F_0=\dirlm{\<\<n\<>0\,\,\>}\Hom_S(S/I^n\<,\>F_0).
$
Since $\textup{Tor}_i^S\<$ commutes with $\dirlm{}\!\!$\vspace{-3pt} one may assume\vs{.6}
that  $F_0$ is annihilated by some fixed power $I^n\<$,
whence so is $\textup{Tor}_i^S(E_0,\>F_0)$, whence the conclusion.\vs1

\pagebreak[3]
As for $\Avc$, some caution must be taken because being in~$\Avc$ is
not a local property. But since ${\T}or_i$ commutes with
$\dirlm{}\!\!$\vspace{-2.4pt} one may assume that $E$ and $F$ have
coherent homology, and then the problem is to show\vspace{.4pt} that
so does ${\T}or_i(E,F\>)$. This problem \emph{is} local, and so one can
use the equivalence of categories described in
\mbox{\cite[p.\,31, Proposition 3.1.1]{DFS}} to reduce to the analogous---and
easily handled---problem for finitely-generated modules over a
noetherian ring.
\end{proof}

\begin{sprop}\label{P:projn}  Let\/ $\psi\colon X\to Y$ be a map of 
noetherian formal schemes. For all\/ $F\in\Dvc(X),$ $G\in\Dvc(Y),$  the
projection map is an\/ \emph{isomorphism}
$$
\smash{\R\psi_*F\OT{X} G\iso \R\psi_*(F\OT{X} \bL\psi^*G).}
$$
\end{sprop}

\begin{proof}
Once the necessary preliminaries are in place the proof is
essentially that of \cite[Proposition 3.9.4]{Dercat}.  These
preliminaries are as follows.\vspace{1.4pt}

1) The question is local on $Y$ (cf.~e.g., \emph{loc.\,cit.}), so one
   can assume that $Y$ is affine. Then \cite[p.\,37, Prop.\,3.3.5]{DFS} gives
   $\bL\psi^*G\in\Dvc(X)$, and so, by \ref{monoid1}, 
   $F\otimes\bL\psi^*G\in\Dvc(X)$.\vspace{1pt}

2) The functor $\R\psi_*$ is bounded-above on $\Dvc(X)$ \cite[p.\,39,
   Prop.\,3.4.3\kern.4pt (b)]{DFS}.\va{1.5}

3) The functors $\R\psi_*$, $\bL\psi^*$ and~$\OT{\,X}$ all\va{-3}
   commute with direct sums: for the first, see \cite[p.\,41, Prop.\,3.5.2]{DFS}, 
   and for the last two see \cite[Prop.\,(3.8.2)]{Dercat}.\va{1.5}

4) For any noetherian formal scheme $Z$, $\Avc\>(Z)$ is a
    \emph{plump} subcategory of $\A(Z)$ \cite[p.\,34,
    Prop.\,3.2.2]{DFS}.\vspace{1pt}

5) Over an affine noetherian formal scheme $Z$, every object in
   $\Avc\>(Z)$ is a homomorphic image of a free $\cO_{\<\<Z}$-module
   \cite[p.\,32, Corollary 3.1.4]{DFS}. \vspace{2pt}

These facts enable a ``way-out" reduction of the proof of
 \Pref{P:projn} for bounded-above $\Dvc\,$-\kf complexes to the simple
 case where $G= \cO_Y$ (cf.~proof of \cite[Proposition
 3.9.4]{Dercat}).  Then for the unbounded case, one uses that
 $\Avc\>(X)$ is stable under $\dirlm{}\!\!$.

The rest is left to the reader.
\end{proof}


\newpage
\begin{thebibliography}{AJL3}

\bibitem[AJL1]{AJL} L.\:Alonso Tarr{\'{\i}}o, A.\:Jerem{\'{\i}}as L{\'o}pez
and J.\:Lipman, Local homology and cohomology on schemes,
{\it Ann.\:Scient.\:{\'E}c.\:Norm.\:Sup.\;}{\bf 30}\,(1997), 1--39, plus
CORRECTIONS at {\tt www.math.purdue.edu/\~{}lipman}, or p.\,879, 
vol.\,2 of  {\it Collected Papers of
\mbox{Joseph Lipman},} Queen's Papers in Pure and Applied Math., Vol.\;{\bf 117},
Queen's University, Kingston, \mbox{Ontario,} Canada, 2000.

\bibitem[AJL2]{DFS}\bysame, Duality and flat base change on formal schemes, 
{\it Contemporary Math.,} Vol.\;{\bf 244}, Amer. Math. Soc., Providence, R.I.
(1999), 3--90.

\bibitem[AJS1]{AJS1} L.\:Alonso Tarr{\'{\i}}o, A.\:Jerem{\'{\i}}as L{\'o}pez
and M.\,J.\:Souto Salorio, Localization in
categories of complexes and unbounded resolutions,
{\it Canadian J.~Math.,} {\bf 52}  (2000),  225--247.

\bibitem[AJS2]{AJS2} \bysame, Bousfield localization on formal schemes,
{\it J. Algebra,}  {\bf 278}  (2004),  585--610.

\bibitem[\kf BN]{BN} M.\:B\"okstedt and A.\:Neeman, Homotopy limits in
triangulated categories, {\it Compositio Math.\:}{\bf 86}~(1993),
209--234.

\bibitem[\kf Brb]{Bou} N.\:Bourbaki, {\it Alg\`ebre Commutative,}
Actualit\'es Sci. et Industrielles, nos.~1290, 1293, Hermann, Paris, 1961.

\bibitem[Gdm]{Go} R.\:Godement, \!{\it Topologie Alg\'ebrique et Th\'eorie des Faisceaux,}
Act.\:Sci.\,et
\kern-1pt Industrielles \kern-1pt 1252, Herrmann, Paris, 1973.

\bibitem[GD]{GD} A.\:Grothendieck and J.\:Dieudonn\'{e}, {\it \'Elements
de G\'{e}om\'{e}trie Alg\'{e}brique\/} {\bf I}, Springer-Verlag, New
York, 1971.

\bibitem[GD3]{GD3} \bysame, 
{\it \'Elements de G\'{e}om\'{e}trie Alg\'{e}brique\/} {\bf III},
Publications Math. IHES {\bf 11}, Paris, 1961.

\bibitem[GR2]{SGA2} A.\:Grothendieck, 
{\it Cohomologie locale des faisceaux coh\'erents et th\'eor\`emes de Lefschetz locaux et globaux (SGA 2),} 
{\tt arxiv.org/pdf/math/0511279v1.pdf.}

\bibitem[GS]{GS} 
J.\,A.\:Navarro Gonz\'alez  and J.\,B.\:Sancho de Salas,
Sections of quasi-coherent sheaves,
Comm. Algebra 20 (1992), no. 8, 2289--2293. 

\bibitem[\kf Hg]{Hg} M.\:Hogancamp, \mbox{Idempotents in triangulated monoidal
categories,\hskip88pt{}} {\tt arxiv.org/pdf/1703.01001.pdf}.

\bibitem[\kf Hrt]{H} R.\:Hartshorne, {\it Residues and Duality,} Lecture
Notes in Math., no.\,{\bf 20}, Springer-Verlag, New York, 1966.

\bibitem[Kf\kern1pt]{Kf} G.\,R.\:Kempf, Some elementary proofs of basic
theorems in the cohomology of quasi-coherent sheaves, {\it Rocky
Mountain J.~Math.} {\bf 10} (1980), 637--645.

\bibitem[\kf Lp1]{Dercat}  J.\:Lipman, {\it Notes on Derived Functors and Grothendieck Duality,}
Lecture
Notes in Math., no.\,{\bf 1960}, Springer-Verlag, New York, 2009.

\bibitem[\kf Lp2]{Lectures} \bysame,  Lectures on local cohomology and
duality, in {\it Local Cohomology and Its Applications,}
 (ed.~G.\:Lyubeznik), Marcel Dekker, New York, 2001, 39--89. 

\bibitem[Lu]{Lu} J.\:Lurie,  \textit{Higher Algebra} (Sept.~18, 2017). 
\texttt {https://www.math.ias.edu/~lurie/ }

\bibitem[\kf Mc]{Ma2} S.\:Mac\:Lane, Categories for the Working
 Mathematician, second edition, Springer, New York, 1998.

\bibitem[\kf Nm1]{chrom} A.\:Neeman, The chromatic tower for $\D(R)$, {\it Topology}
{\bf 31} (1992), 519--532.

\bibitem[\kf Nm2]{Nmn} \bysame, The Grothendieck duality theorem via Bousfield's
techniques and Brown representability, {\it J.\:Amer.\:Math.\:Soc.}\;{\bf
9}\,(1996), 205--236.

\bibitem[Spn]{Sp} N.\:Spaltenstein, Resolutions of unbounded complexes,  {\it
Comp.~Math.}\;{\bf 65}$\mspace{-.1mu}$(1988), \mbox{121--154.}

\bibitem[St\kf]{Stacks} Stacks Project, {\tt https://stacks.math.columbia.edu}

\end{thebibliography}
\end{document}